\documentclass[oneside,11pt]{article}
\usepackage[top=22mm,bottom=22mm,left=25mm,right=25mm]{geometry}
\usepackage{etoolbox}
\usepackage{setspace}
\usepackage[T1]{fontenc}
\usepackage[utf8]{inputenc}
\usepackage{lmodern}

\makeatletter
\patchcmd{\thebibliography}
{\list}
{\small\setstretch{1.2}\list}
{}{}
\makeatother
\usepackage{amsmath,amssymb,amsthm,bm}
\usepackage{graphicx,array,tabularx}
\usepackage[font=small,labelfont=bf]{caption}
\usepackage{subcaption}
\usepackage{cite,needspace,float}
\usepackage{xurl}
\usepackage{xcolor}
\usepackage[colorlinks,citecolor=blue,
linkcolor=blue,hypertexnames=false]{hyperref} 
\usepackage[section]{placeins}
\allowdisplaybreaks

\title{{\sl Extrapolating escape dynamics from non-escaping trajectories in Hamiltonian systems}}

\author{Jianming Liu\thanks{School of Mathematics, China University of Mining and Technology, Xuzhou, 221116, P. R. China. Email: \url{liujianming1026@outlook.com}} 
	\and 
	Zhicheng Tong\thanks{Corresponding author. School of Mathematics, Jilin University, Changchun 130012, P. R. China. Email: \url{tongzc25@jlu.edu.cn}}
}

\date{}

\begin{document}
	\maketitle
	
	\begin{abstract}
		Hamiltonian systems exhibit rich dynamical behaviors, among which escape is a typical phenomenon. Below the escape threshold, trajectories initiated inside the central potential well remain bounded. Once the energy exceeds this threshold, and the initial conditions lie within the escape basin, the system eventually evolves along unbounded orbits toward infinity. Predicting the escape dynamics is important for understanding the evolution of these systems. In this study, we employ the next-generation reservoir computing framework, a powerful machine learning architecture widely used in dynamical prediction, to address this challenging issue. To capture the intrinsic features of the system, we train the model using data from multiple energy values. Numerical experiments on two Hamiltonian systems show that, with only three time series as training data, each corresponding to an energy value below the escape threshold, the model is capable of not only predicting the escape threshold itself, but also reliably predicting the structure of the escape channels. As a primary contribution of this work, these findings demonstrate that escape dynamics can be inferred from subcritical data alone in the two systems studied, offering a promising paradigm for data-driven studies of Hamiltonian systems.
		\\
		\\
		{\bf Keywords:} Hamiltonian systems, escape dynamics, escape threshold, escape channels, next-generation reservoir computing\vspace{2mm}
		\\
		\noindent{\bf 2020 Mathematics Subject Classification:} 37M05, 68T07, 37M10, 37J06
	\end{abstract}
	
	\tableofcontents
	
	\section{Introduction}\label{sec:introduction}
	Hamiltonian systems form a class of dynamical systems that govern the evolution of a physical system through its total energy. The phase space of such a system possesses a symplectic structure, which endows the flow with two fundamental invariants: conservation of energy (for autonomous systems) and conservation of phase-space volume. Depending on the energy level and the topology of the potential energy surface, the dynamics can be classified into bounded and unbounded regimes. Escape dynamics describes the process by which a trajectory, initially confined to a bounded region, eventually escapes to infinity. This phenomenon is prevalent in open Hamiltonian systems, where the potential energy function has saddle-type stationary points at which exit channels open in the energy landscape. For energies above the escape threshold, not all trajectories escape; rather, the energy surface is partitioned into basins of escape and basins of bounded motion.
	
	Beyond its theoretical significance, escape dynamics also plays a key role in a variety of practical contexts. In celestial mechanics, it governs the gravitational ejection of stars from globular clusters \cite{fukushige2000-20}, the scattering and ejection of planets in multiplanetary systems \cite{ford2008-18}, and the long-term diffusion of comets into interstellar space \cite{duncan1987-16}; the escape outcome determines whether the body remains bound or is ejected. In molecular physics, classical reaction dynamics, including dissociation, can be described as transport out of a reactant region through phase-space bottlenecks near saddle points \cite{waalkens2008-44}; the relative populations of competing product channels determine branching ratios \cite{ezra2009-17}. In atomic physics, the microwave-driven ionization of a Rydberg atom provides a paradigmatic example of Hamiltonian escape dynamics, where the electron escapes from the Coulomb well through chaotic diffusion \cite{casati1987-6}; the ionization threshold and probability depend on the driving field and the initial state \cite{casati1987-6,leeuwen1985-43}. Across these phenomena, the common thread is transport out of an initially bound interaction region, together with sensitivity to initial conditions and system parameters. Predicting whether, when, and through which channel escape will occur is of great practical importance, yet remains highly challenging due to the intricate phase-space mechanisms governing chaotic transport. In the autonomous systems studied here, the phenomenon is escape from a potential well, investigated in the framework of chaotic escape or chaotic scattering~\cite{aguirre2001-1,contopoulos1990-11,contopoulos2004-12}, where unstable periodic orbits and their stable and unstable manifolds organize transport through the exits~\cite{contopoulos2004-12,waalkens2008-44}. Furthermore, escape is compatible with bounded motion on invariant tori, as regular islands can persist above the escape energy~\cite{aguirre2001-1,barrio2009-3}.
	
	The recurrent neural network (RNN) is a powerful computational framework for modeling time series. The traditional RNN is trained by gradient descent, incurring high computational costs. The echo state network (ESN) \cite{jaeger2001-26,jaeger2004-27,lukoeviius2009-34}, a special type of RNN, adopts a distinct approach to weight setting. Unlike the traditional RNN, in which all weights are trainable, an ESN randomly initializes and fixes its input and hidden layer weights. Only the output weights are learned, typically via linear regression. The architecture of the ESN incorporates a high-dimensional dynamical system, termed the reservoir; hence this paradigm is also known as reservoir computing (RC). A recent development of RC is the next-generation reservoir computing (NG-RC) \cite{barbosa2022-2,gauthier2021-21}. In contrast to classical RC, which relies on the intrinsic dynamics of randomly connected neurons to generate high-dimensional feature vectors, the NG-RC framework constructs the feature vector explicitly by concatenating the delayed input states with their polynomial terms. This design significantly reduces the number of hyperparameters to be tuned, thereby lowering the computational cost and simplifying implementation. NG-RC has been widely applied to nonlinear system modeling tasks, such as chaotic dynamics prediction \cite{chepuri2024-10,cox2024-15,gauthier2025-23,grigoryeva2025-24,liu2023-32,liu2024-33,pyragas2026-35,santos2025-37,wang2025-45,wang2026-46,yang2025-47,zhang2025-50}, bifurcation diagram reconstruction \cite{kglmayr2024-31,ratas2024-36,zhang2025-51}, system control \cite{cheng2026-9,kent2024-30}, coexisting attractors learning \cite{gauthier2022-22}, dynamical inference \cite{cestnik2026-8}, quantum dynamics prediction \cite{sornsaeng2024-41}, and superconducting qubit readout \cite{kent2026-29}.
	
	The data-driven approach is an important means of studying Hamiltonian systems. Despite some progress \cite{bousakreltayar2026-4,bramburger2026-5,celletti2022-7,hasmi2026-25,jimnezlpez2025-28,saltas2024-38,tong2026-42,zhang2021-48,zhang2022-49}, many problems remain to be further explored. In this study, we address the prediction of escape dynamics in Hamiltonian systems. Owing to its unique design, NG-RC can effectively extract the dynamical features hidden in data. Trained with data from only a few energy values, the model can predict the escape threshold with high accuracy and reliably reconstruct the structure of the escape channels. To the best of our knowledge, this represents the first result demonstrating that \textit{critical escape behaviors can be anticipated without ever observing actual escape events,} highlighting a key insight of this study.
	
	The remainder of this paper is organized as follows. Section~\ref{sec:method} briefly introduces the Hamiltonian system, formulates the problem, and details the NG-RC prediction method. Section~\ref{sec:results} reports the simulation results on two representative systems. Finally, Section~\ref{sec:conclusion} summarizes the main findings of this work and concludes the paper.
	
	\Needspace{4\baselineskip}\section{Problem statement and method}\label{sec:method}
	\Needspace{4\baselineskip}\subsection{Problem statement}\label{sec:problem}
	In the canonical formulation, the state of a system with $n$ degrees of freedom is described by a set of generalized coordinates $\boldsymbol{q}=(q_{1},\allowbreak q_{2},\allowbreak \ldots{},\allowbreak q_{n})$ and their conjugate momenta $\boldsymbol{p}=(p_{1},\allowbreak p_{2},\allowbreak \ldots{},\allowbreak p_{n})$, which together form a $2n$-dimensional phase space. The dynamics are prescribed by a scalar function $H(\boldsymbol{q},\allowbreak \boldsymbol{p},\allowbreak t)$, known as the Hamiltonian, which typically represents the total energy of the system. The equations of motion are given by Hamilton's canonical equations:
	
	\begin{equation*}
		\dot{\boldsymbol{q}}=\frac{\partial{}H}{\partial{}\boldsymbol{p}},\allowbreak \qquad \dot{\boldsymbol{p}}=-\frac{\partial{}H}{\partial{}\boldsymbol{q}},\allowbreak 
	\end{equation*}
	where the overdot denotes differentiation with respect to time $t$.
	
	For an autonomous system, i.e., one in which the Hamiltonian does not explicitly depend on time, the Hamiltonian itself is a conserved quantity. In this case, $H(\boldsymbol{q},\allowbreak \boldsymbol{p})=E$, where $E$ is the total energy of the system. For the natural Hamiltonians considered here, the Hamiltonian is expressed as the sum of kinetic energy $T(\boldsymbol{p})$ and potential energy $V(\boldsymbol{q})$:
	
	\begin{equation*}
		H(\boldsymbol{q},\allowbreak \boldsymbol{p})=T(\boldsymbol{p})+V(\boldsymbol{q})=\frac{1}{2}\boldsymbol{p}^{T}\boldsymbol{M}^{-1}\boldsymbol{p}+V(\boldsymbol{q}).
	\end{equation*}
	Here, $\boldsymbol{M}$ is a symmetric positive-definite mass matrix, often taken as the identity in simple systems. For a system with $n$ degrees of freedom, the total energy $E$ acts as a constant of motion; at regular energy levels, it constrains the system's evolution to a $(2n-1)$-dimensional energy hypersurface embedded in the phase space. Consequently, the value of $E$ determines the accessible regions of the configuration space. Since the total energy is the sum of kinetic and potential energy, and the kinetic energy is always non-negative, the motion can only take place where $V(\boldsymbol{q})\leqslant{}E$. The boundary $V(\boldsymbol{q})=E$ is the zero-kinetic-energy boundary; at regular boundary points, a trajectory comes to rest and reverses direction. As $E$ increases, the allowed region expands, and new parts of the configuration space become accessible. When $E$ exceeds a critical value, which is typically the energy of a saddle point, the central allowed component becomes connected to infinity, giving rise to escape channels. This is the central phenomenon of escape dynamics.
	
	The equations governing the trajectories of a Hamiltonian system are obtained by explicitly substituting the Hamiltonian into Hamilton's canonical equations. For a two-degree-of-freedom system with generalized coordinates $(x,\allowbreak y)$ and conjugate momenta $(p_{x},\allowbreak p_{y})$, the Hamiltonian typically takes the form:
	
	\begin{equation}
		H(x,\allowbreak y,\allowbreak p_{x},\allowbreak p_{y})=\frac{1}{2}(p_{x}^{2}+p_{y}^{2})+V(x,\allowbreak y),\allowbreak 
		\label{eq:hamiltonian}
	\end{equation}
	where we have assumed unit mass. Substituting the Hamiltonian into Hamilton's canonical equations yields the equations of motion:
	
	\begin{equation}
		\begin{aligned}\dot{x}&=\frac{\partial{}H}{\partial{}p_{x}}=p_{x},\allowbreak  \\ \dot{y}&=\frac{\partial{}H}{\partial{}p_{y}}=p_{y},\allowbreak  \\ \dot{p}_{x}&=-\frac{\partial{}H}{\partial{}x}=-\frac{\partial{}V}{\partial{}x},\allowbreak  \\ \dot{p}_{y}&=-\frac{\partial{}H}{\partial{}y}=-\frac{\partial{}V}{\partial{}y}.\end{aligned}
		\label{eq:hamilton-equations}
	\end{equation}
	For an initial value $\boldsymbol{h}(0)=(x_{0},\allowbreak y_{0},\allowbreak p_{x0},\allowbreak p_{y0})$, according to Eq. \eqref{eq:hamiltonian}, the corresponding energy is $E_{0}=H(\boldsymbol{h}(0))$. Starting from $\boldsymbol{h}(0)$, numerical integration of Eq. \eqref{eq:hamilton-equations} yields the trajectory data on the energy surface $E_{0}$. The goal of this study is to train an NG-RC model on trajectory data with energies below the escape threshold and then predict the escape dynamics of the system.
	
	\Needspace{4\baselineskip}\subsection{Method}\label{sec:ngrc}
	Assume that the training data are sampled from $m$ energy values $E_{j}$, $j=1,\allowbreak 2,\allowbreak \cdots{},\allowbreak m$. Suppose $\boldsymbol{h}_{j}(1),\allowbreak \boldsymbol{h}_{j}(2),\allowbreak \boldsymbol{h}_{j}(3),\allowbreak \ldots{},\allowbreak \boldsymbol{h}_{j}(L)\in{}R^{d}$ are the states of a $d$-dimensional Hamiltonian system with energy $E_{j}$. For NG-RC, the total feature vector $\boldsymbol{O}_{j}(t)$ associated with $\boldsymbol{h}_{j}(t)$ is given by \cite{gauthier2021-21}:
	
	\begin{equation}
		\boldsymbol{O}_{j}(t)=c\oplus{}\boldsymbol{O}_{j}(t)^{(1)}\oplus{}\boldsymbol{O}_{j}(t)^{(2)},\allowbreak 
		\label{eq:feature-vector}
	\end{equation}
	where $\oplus{}$ represents the vector concatenation operation, $c$ is the bias, $\boldsymbol{O}_{j}(t)^{(1)}$ is the linear vector, and $\boldsymbol{O}_{j}(t)^{(2)}$ is the nonlinear vector. We set $c=1$ in this study. The linear part is formed by concatenating the energy value with $k$ state vectors:
	
	\begin{equation}
		\boldsymbol{O}_{j}(t)^{(1)}=E_{j}\oplus{}\boldsymbol{h}_{j}(t)\oplus{}\boldsymbol{h}_{j}(t-1)\oplus{}\ldots{}\oplus{}\boldsymbol{h}_{j}(t-k+1).
		\label{eq:linear-features}
	\end{equation}
	In Eq. \eqref{eq:linear-features}, $k$ is a hyperparameter that needs to be specified in advance. We include $E_{j}$ as an additional dimension in $\boldsymbol{O}_{j}(t)^{(1)}$ to enhance the model's awareness of the energy. Obviously, the dimension of $\boldsymbol{O}_{j}(t)^{(1)}$ is $D_{1}=1+d\times{}k$. Writing $\boldsymbol{h}_{j}(t)=[h_{j1}(t),\allowbreak h_{j2}(t),\allowbreak \ldots{},\allowbreak h_{jd}(t)]^{T}\in{}R^{d}$, the elements of $\boldsymbol{O}_{j}(t)^{(2)}$ are obtained by taking all pairwise products of the elements of $\boldsymbol{O}_{j}(t)^{(1)}$:
	
	\begin{equation}
		\boldsymbol{O}_{j}(t)^{(2)}=[\boldsymbol{O}_{j,\allowbreak r}(t)^{(1)}\times{}\boldsymbol{O}_{j,\allowbreak s}(t)^{(1)}]_{1\leqslant{}r\leqslant{}s\leqslant{}D_{1}}^{T},\allowbreak 
		\label{eq:quadratic-features}
	\end{equation}
	where $\boldsymbol{O}_{j,\allowbreak r}(t)^{(1)}$ denotes the $r$-th component of $\boldsymbol{O}_{j}(t)^{(1)}$, and $D_{1}$ is the dimension of $\boldsymbol{O}_{j}(t)^{(1)}$. From Eq. \eqref{eq:quadratic-features}, the dimension of $\boldsymbol{O}_{j}(t)^{(2)}$ is $D_{2}=D_{1}(1+D_{1})/2$. Therefore, the dimension of the feature vector $\boldsymbol{O}_{j}(t)$ is $1+D_{1}+D_{2}$.
	
	Each time series is independently fed into the model through Eq. \eqref{eq:feature-vector}, generating the feature vectors $\boldsymbol{O}_{j}(k),\allowbreak \boldsymbol{O}_{j}(k+1),\allowbreak \ldots{},\allowbreak \boldsymbol{O}_{j}(L-1)$. For each feature vector $\boldsymbol{O}_{j}(t)$, the target output is given by $\boldsymbol{y}_{j}(t+1)=\boldsymbol{h}_{j}(t+1)-\boldsymbol{h}_{j}(t)$, where the time index satisfies $k\leqslant{}t\leqslant{}L-1$. Denote the feature vector matrix as $\boldsymbol{O}_{j} = (\boldsymbol{O}_{j}(k),\allowbreak  \boldsymbol{O}_{j}(k+1),\allowbreak  \ldots{},\allowbreak  \boldsymbol{O}_{j}(L-1))$, and the desired output matrix as $\boldsymbol{Y}_{j} = (\boldsymbol{y}_{j}(k+1),\allowbreak  \boldsymbol{y}_{j}(k+2),\allowbreak  \ldots{},\allowbreak  \boldsymbol{y}_{j}(L))$. The readout weight matrix $\boldsymbol{W}_{out}$ is fitted to the target relation $\boldsymbol{Y}_{j} = \boldsymbol{W}_{out}\boldsymbol{O}_{j}$, $j=1,\allowbreak 2,\allowbreak \ldots{},\allowbreak m$. Then $\boldsymbol{W}_{out}$ can be computed using ridge regression:
	
	\begin{equation}
		{\min}_{\boldsymbol{W}_{out}}(\sum{}_{j=1}^{m}\,||\boldsymbol{W}_{out}\boldsymbol{O}_{j}-\boldsymbol{Y}_{j}||_{\operatorname{F}}^{2}+\gamma{}||\boldsymbol{W}_{out}||_{\operatorname{F}}^{2}),\allowbreak 
		\label{eq:ridge-objective}
	\end{equation}
	where $\gamma{}>0$ is the Tikhonov regularization parameter, and F denotes the Frobenius norm. Let $\boldsymbol{O}=[\boldsymbol{O}_{1},\allowbreak \boldsymbol{O}_{2},\allowbreak \cdots{},\allowbreak \boldsymbol{O}_{m}]$ be the concatenated feature matrix, and $\boldsymbol{Y}=[\boldsymbol{Y}_{1},\allowbreak \boldsymbol{Y}_{2},\allowbreak \cdots{},\allowbreak \boldsymbol{Y}_{m}]$ be the concatenated target matrix. For Eq. \eqref{eq:ridge-objective}, we have the following closed-form solution:
	
	\begin{equation}
		\boldsymbol{W}_{out}=\boldsymbol{Y}\boldsymbol{O}^{T}(\boldsymbol{O}\boldsymbol{O}^{T}+\gamma{}\boldsymbol{I})^{-1},\allowbreak 
		\label{eq:ridge-solution}
	\end{equation}
	where $\boldsymbol{I}$ is the identity matrix.
	
	During prediction, feature indices label the next state. The learned $\boldsymbol{W}_{out}$ predicts evolution at different energies. For an energy value $E_{i}$, and the corresponding initial value $\boldsymbol{h}_{i}(0)=(x_{0,\allowbreak i},\allowbreak y_{0,\allowbreak i},\allowbreak p_{x0,\allowbreak i},\allowbreak p_{y0,\allowbreak i})$, the feature vector $\boldsymbol{O}_{i}(1)$ is constructed via Eq. \eqref{eq:feature-vector}. Here, as the initial condition $\boldsymbol{h}_{i}(0)$ is the only known state, we replicate it $k$ times to fill the delay window for the linear vector $\boldsymbol{O}_{i}(1)^{(1)}$. Then the first predicted state is computed as $\boldsymbol{h}_{i}(1)=\boldsymbol{W}_{out}\boldsymbol{O}_{i}(1)+\boldsymbol{h}_{i}(0)$. Similarly, after obtaining the feature vector $\boldsymbol{O}_{i}(t)$, the $t$-th predicted state is computed as $\boldsymbol{h}_{i}(t)=\boldsymbol{W}_{out}\boldsymbol{O}_{i}(t)+\boldsymbol{h}_{i}(t-1)$.
	
	Once the predictions are generated, the next step is to analyze the predicted data to reveal the underlying dynamical characteristics. This analysis allows us to identify the escape threshold and escape basins on the chosen initial-condition slices.
	
	\Needspace{4\baselineskip}\section{Results}\label{sec:results}
	\Needspace{4\baselineskip}\subsection{A Hamiltonian system with two exit channels}\label{sec:two-exit}
\begin{figure}[!htbp]
		\centering
		\begin{subfigure}[t]{.49\linewidth}
			\centering
			\includegraphics[width=\linewidth]{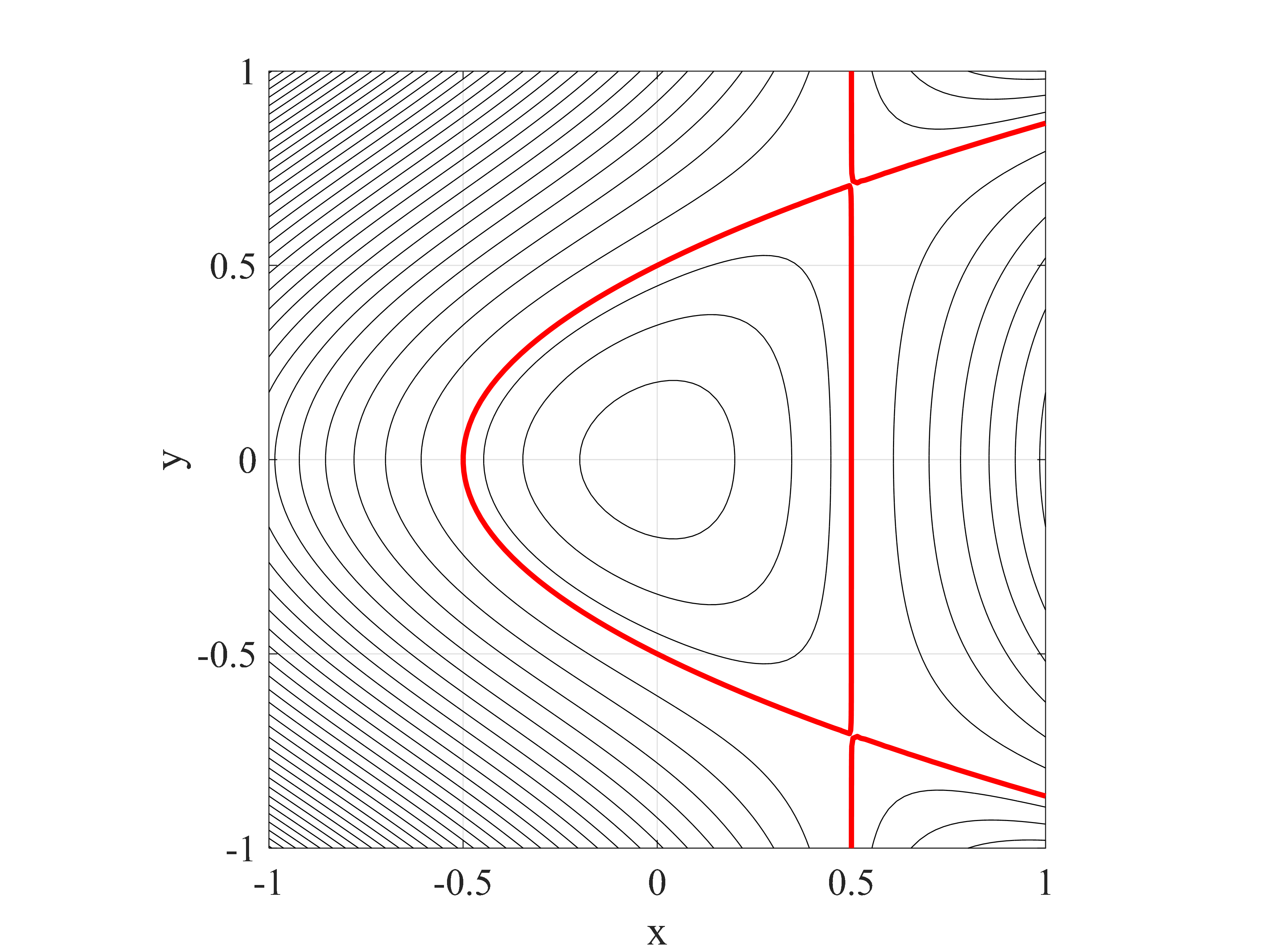}
			\caption{}\label{fig:two-exit-geometry:a}
		\end{subfigure}\hfill
		\begin{subfigure}[t]{.49\linewidth}
			\centering
			\includegraphics[width=\linewidth]{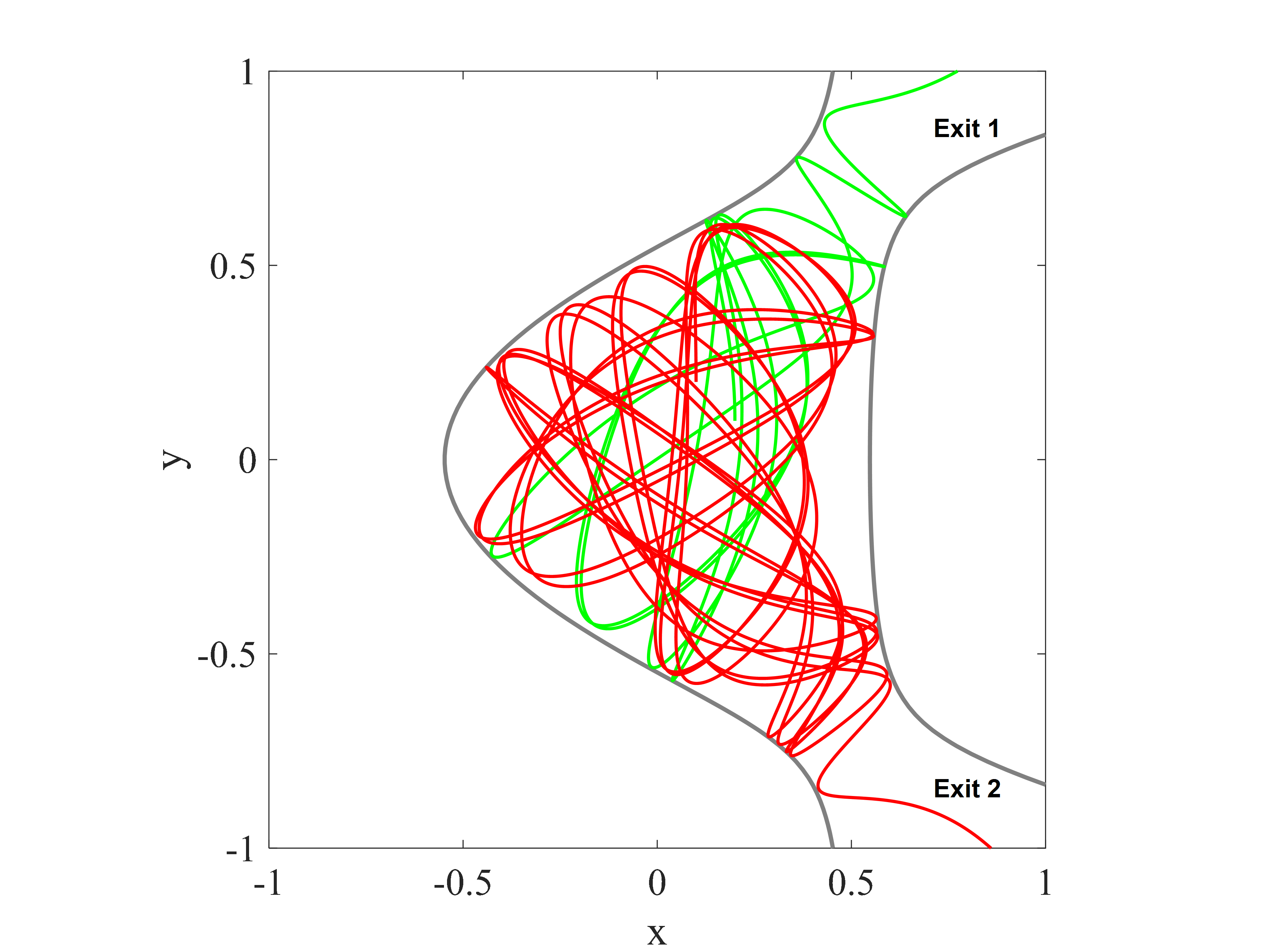}
			\caption{}\label{fig:two-exit-geometry:b}
		\end{subfigure}
		\caption{Equipotential curves and exit channels of the system. (\subref{fig:two-exit-geometry:a}) Equipotential curves for different energy values. The equipotential curve corresponding to the escape energy is shown in red. (\subref{fig:two-exit-geometry:b}) The two exit channels and two escaping trajectories at $E=0.15$. The trajectory escaping from exit 1 is plotted in green, while the trajectory escaping through exit 2 is plotted in red.}
		\label{fig:two-exit-geometry}
	\end{figure}

	Coupled-oscillator potentials of this type provide simple models for studying stellar escape from a galactic potential \cite{contopoulos1990-11,contopoulos2004-12}. Their unstable periodic orbits, chaotic scattering, and escape basins make them useful benchmark systems for investigating transport in open Hamiltonian systems \cite{contopoulos2004-12,contopoulos2012-13,zotos2014-52}. We consider a Hamiltonian system with the following potential function \cite{contopoulos1990-11,contopoulos2004-12,contopoulos2012-13,zotos2014-52}:
	
	\begin{equation}\notag
		V(x,\allowbreak y)=\frac{1}{2}(x^{2}+y^{2})-xy^{2}.
	\end{equation}
	With the conjugate momenta $p_{x}=\dot{x}$ and $p_{y}=\dot{y}$, the Hamiltonian reads
	
	\begin{equation}
		H=\frac{1}{2}(p_{x}^{2}+p_{y}^{2}+x^{2}+y^{2})-xy^{2}.
		\label{eq:two-exit-hamiltonian}
	\end{equation}
	From Eq. \eqref{eq:two-exit-hamiltonian}, we can obtain the equations of motion:
	
	\begin{equation}
		\begin{aligned}\dot{x}&=p_{x},\allowbreak  \\ \dot{y}&=p_{y},\allowbreak  \\ \dot{p}_{x}&=-x+y^{2},\allowbreak  \\ \dot{p}_{y}&=-y+2xy.\end{aligned}
		\label{eq:two-exit-motion}
	\end{equation}
	The equipotential curves and exit channels of system \eqref{eq:two-exit-motion} are illustrated in Fig.~\ref{fig:two-exit-geometry}. This Hamiltonian system has an escape energy $E_{esc}=0.125$. The escape threshold $E_{esc}$ marks the opening of exits from the central well. We use data generated at energies below the escape threshold (i.e., $E<E_{esc}$) to build a predictive model and then apply the model to predict escape dynamics for energies above $E_{esc}$.
	
	The system has two degrees of freedom, corresponding to a four-dimensional phase space. Once the energy $E$ is given, the motion is confined to a three-dimensional energy surface. With $p_{x0}=0$, for any initial value $(x_{0},\allowbreak y_{0})$, the potential energy is
	
	\begin{equation*}
		V(x_{0},\allowbreak y_{0})=\frac{1}{2}(x_{0}^{2}+y_{0}^{2})-x_{0}y_{0}^{2}.
	\end{equation*}
	Here, the initial value $(x_{0},\allowbreak y_{0})$ must satisfy the physical admissibility condition
	
	\begin{equation}
		V(x_{0},\allowbreak y_{0})\leqslant{}E,\allowbreak 
		\label{eq:two-exit-admissibility}
	\end{equation}
	which ensures that the initial kinetic energy $T_{0}=E-V(x_{0},\allowbreak y_{0})$ is nonnegative.
	
	We define the following function:
	
	\begin{equation*}
		g(x,\allowbreak y,\allowbreak E)=\sqrt{2(E-V(x,\allowbreak y))}.
	\end{equation*}
	For $p_{x0}=0$, and initial value $(x_{0},\allowbreak y_{0})$, setting $p_{y0}=g(x_{0},\allowbreak y_{0},\allowbreak E)$ can guarantee that the total energy of the system is exactly $E$.
	
	To learn the escape dynamics of this system, we build an NG-RC model using data from different energies. Specifically, we select three values $E=E_{1}=0.09$, $E=E_{2}=0.1$, and $E=E_{3}=0.11$, and define three initial values:
	
	\begin{equation*}
		\boldsymbol{z}_{1}=(0.1,\allowbreak  0.1,\allowbreak  0,\allowbreak  g(0.1,\allowbreak 0.1,\allowbreak E_{1})),\allowbreak 
	\end{equation*}
	\begin{equation*}
		\boldsymbol{z}_{2}=(0.1,\allowbreak  0.1,\allowbreak  0,\allowbreak  g(0.1,\allowbreak 0.1,\allowbreak E_{2})),\allowbreak 
	\end{equation*}
	\begin{equation*}
		\boldsymbol{z}_{3}=(0.1,\allowbreak  0.1,\allowbreak  0,\allowbreak  g(0.1,\allowbreak 0.1,\allowbreak E_{3})).
	\end{equation*}
	We generate data from Eq. \eqref{eq:two-exit-motion} with a step size of $\Delta{}t=0.03$ under $\boldsymbol{z}_{1}$, $\boldsymbol{z}_{2}$, and $\boldsymbol{z}_{3}$. The obtained three time series, each has a length of 8000, are employed to train the model. The two hyperparameters, the number of state vectors $k$ and the regularization parameter of ridge regression $\gamma{}$, are set as follows: $k=3$, $\gamma{}=5\times{}10^{-4}$. After training, the model is used for prediction. All predictions are made with this single model.
	
	To estimate the escape threshold, we predict the system's evolution under energy values ranging from 0.111 to 0.128 in increments of 0.001. For each energy value, $(x_{0},\allowbreak y_{0})$ is sampled over a uniform grid covering the domain $[-0.8,\allowbreak 0.49]\times{}[-0.8,\allowbreak 0.8]$ with a step size of 0.01, so that all grid points are visited. $p_{x0}$ is fixed at 0 for all trials, while $p_{y0}$ is determined by the function $g(x,\allowbreak y,\allowbreak E)$. Importantly, we only consider those points that satisfy the constraint given by condition \eqref{eq:two-exit-admissibility}, and any grid point that does not meet this criterion is excluded from the simulation. For each valid initial condition, we predict the trajectory and classify its escape behavior over the prediction horizon. The lowest sampled energy with a predicted escape is used as a numerical threshold estimate for this grid and prediction horizon.
	
	\begin{table}[!htbp]
		\centering
		\caption{Predictions of escape ($\checkmark$) and non-escape ($\times$) for energy values from 0.111 to 0.128.}
		\label{tab:threshold-1}
		\small
		\setlength{\tabcolsep}{0pt}
		\renewcommand{\arraystretch}{1.3}
		\begin{tabular}{|>{\centering\arraybackslash}m{16mm}|*{9}{>{\centering\arraybackslash}m{15.6mm}|}}
			\hline
			$E$ & 0.111 & 0.112 & 0.113 & 0.114 & 0.115 & 0.116 & 0.117 & 0.118 & 0.119 \\
			\hline
			Result & $\times$ & $\times$ & $\times$ & $\times$ & $\times$ & $\times$ & $\times$ & $\times$ & $\times$ \\
			\hline
			$E$ & 0.12 & 0.121 & 0.122 & 0.123 & 0.124 & 0.125 & 0.126 & 0.127 & 0.128 \\
			\hline
			Result & $\times$ & $\times$ & $\times$ & $\times$ & $\times$ & $\checkmark$ & $\checkmark$ & $\checkmark$ & $\checkmark$ \\
			\hline
		\end{tabular}
	\end{table}
	
	\begin{table}[!htbp]
		\centering
		\caption{Predictions of escape and non-escape for energy values from 0.12405 to 0.125.}
		\label{tab:threshold-2}
		\small
		\setlength{\tabcolsep}{0pt}
		\renewcommand{\arraystretch}{1.3}
		\begin{tabular}{|>{\centering\arraybackslash}m{16mm}|*{8}{>{\centering\arraybackslash}m{17.55mm}|}}
			\hline
			$E$ & 0.12405 & 0.1241 & 0.12415 & 0.1242 & 0.12425 & 0.1243 & 0.12435 & 0.1244 \\
			\hline
			Result & $\times$ & $\times$ & $\times$ & $\times$ & $\checkmark$ & $\checkmark$ & $\checkmark$ & $\checkmark$ \\
			\hline
			$E$ & 0.12445 & 0.1245 & 0.12455 & 0.1246 & 0.12465 & 0.1247 & 0.12475 & 0.1248 \\
			\hline
			Result & $\checkmark$ & $\checkmark$ & $\checkmark$ & $\checkmark$ & $\checkmark$ & $\checkmark$ & $\checkmark$ & $\checkmark$ \\
			\hline
			$E$ & 0.12485 & 0.1249 & 0.12495 & 0.125 &  &  &  &  \\
			\hline
			Result & $\checkmark$ & $\checkmark$ & $\checkmark$ & $\checkmark$ &  &  &  &  \\
			\hline
		\end{tabular}
	\end{table}
	
\begin{figure}[!htbp]
		\centering
		\begin{subfigure}[t]{.47\linewidth}
			\centering
			\includegraphics[width=\linewidth]{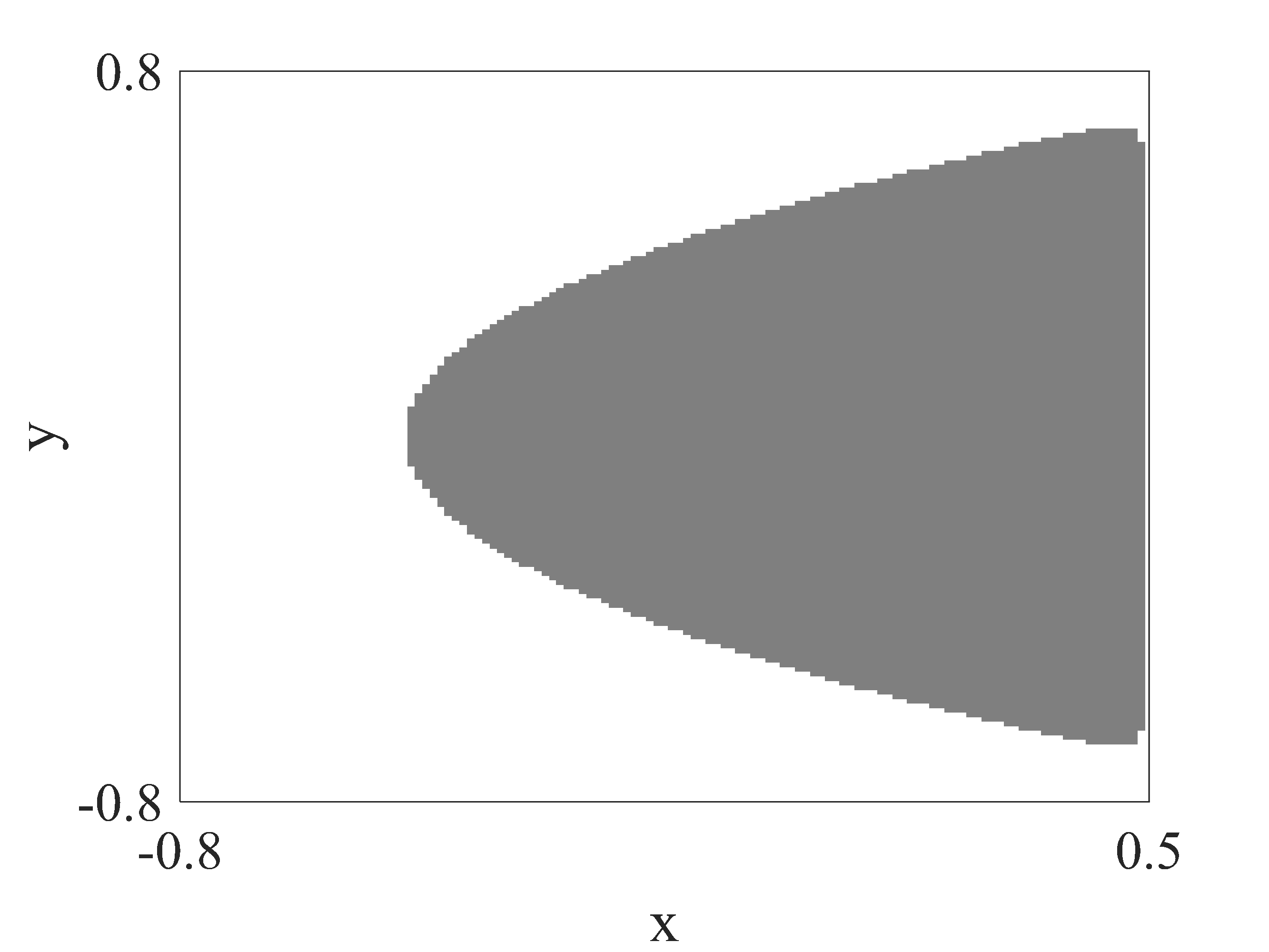}
			\caption{}\label{fig:two-exit-threshold-trajectories:a}
		\end{subfigure}\hfill
		\begin{subfigure}[t]{.47\linewidth}
			\centering
			\includegraphics[width=\linewidth]{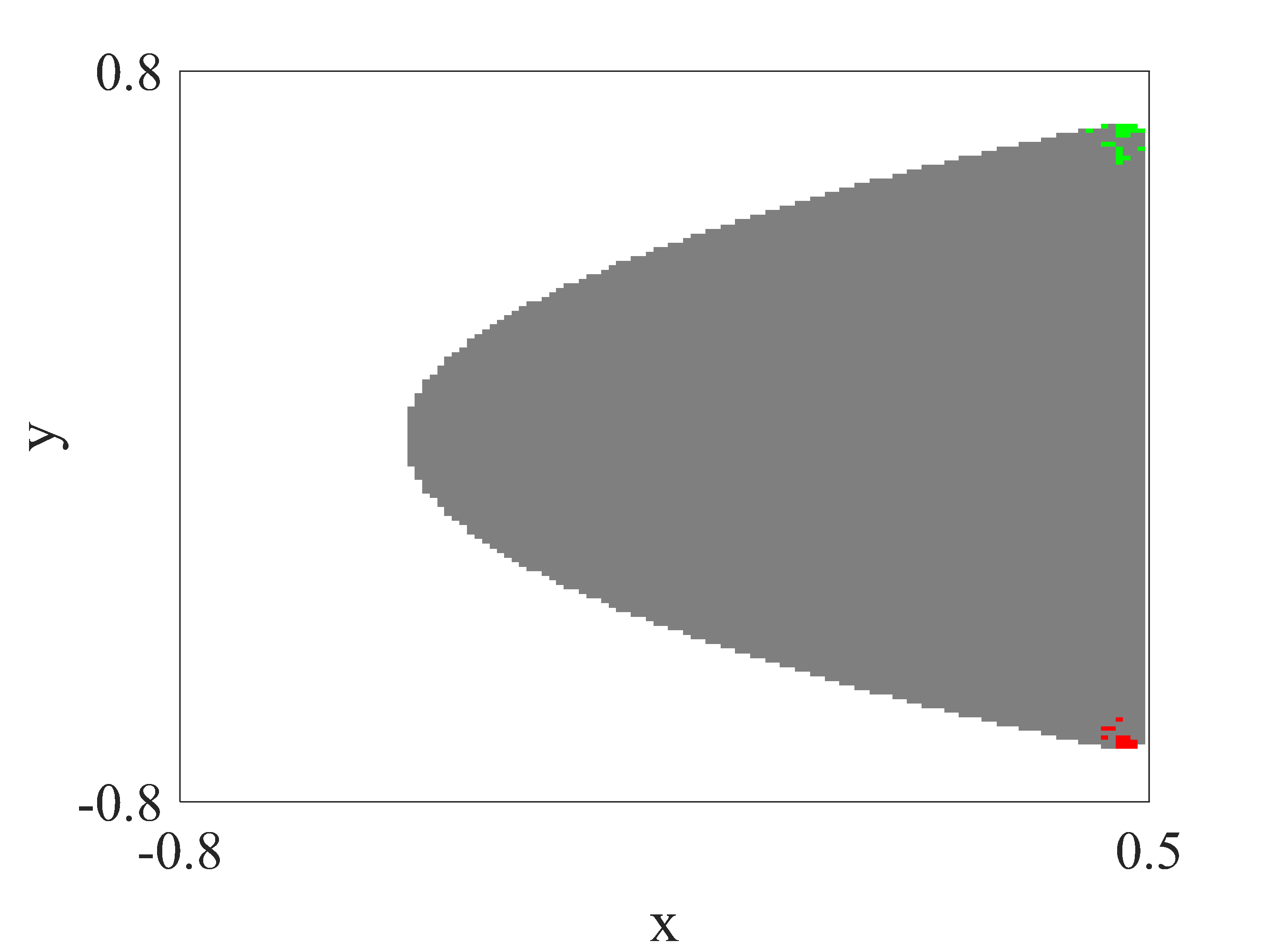}
			\caption{}\label{fig:two-exit-threshold-trajectories:b}
		\end{subfigure}\par\smallskip
		\begin{subfigure}[t]{.47\linewidth}
			\centering
			\includegraphics[width=\linewidth]{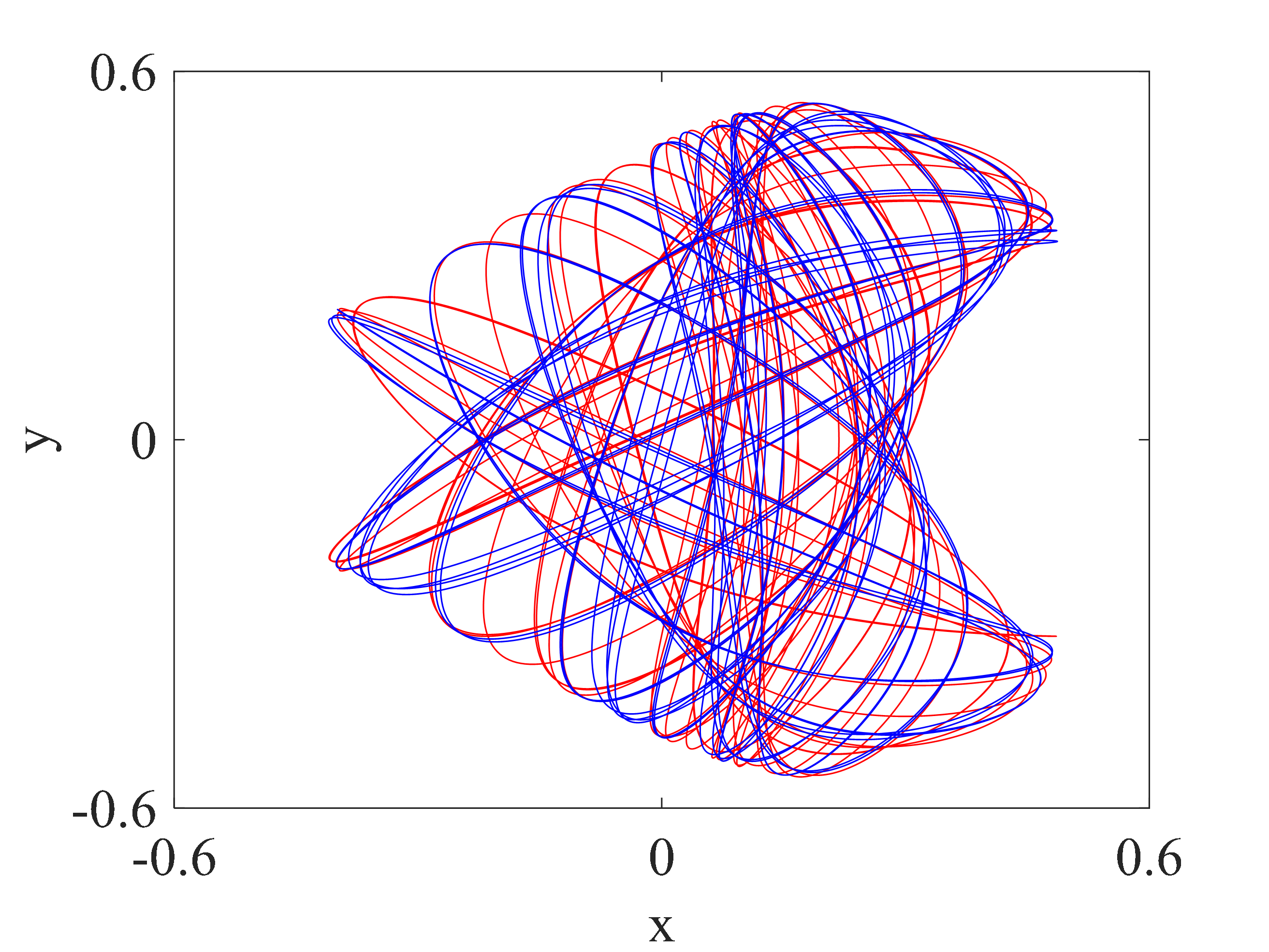}
			\caption{}\label{fig:two-exit-threshold-trajectories:c}
		\end{subfigure}\hfill
		\begin{subfigure}[t]{.47\linewidth}
			\centering
			\includegraphics[width=\linewidth]{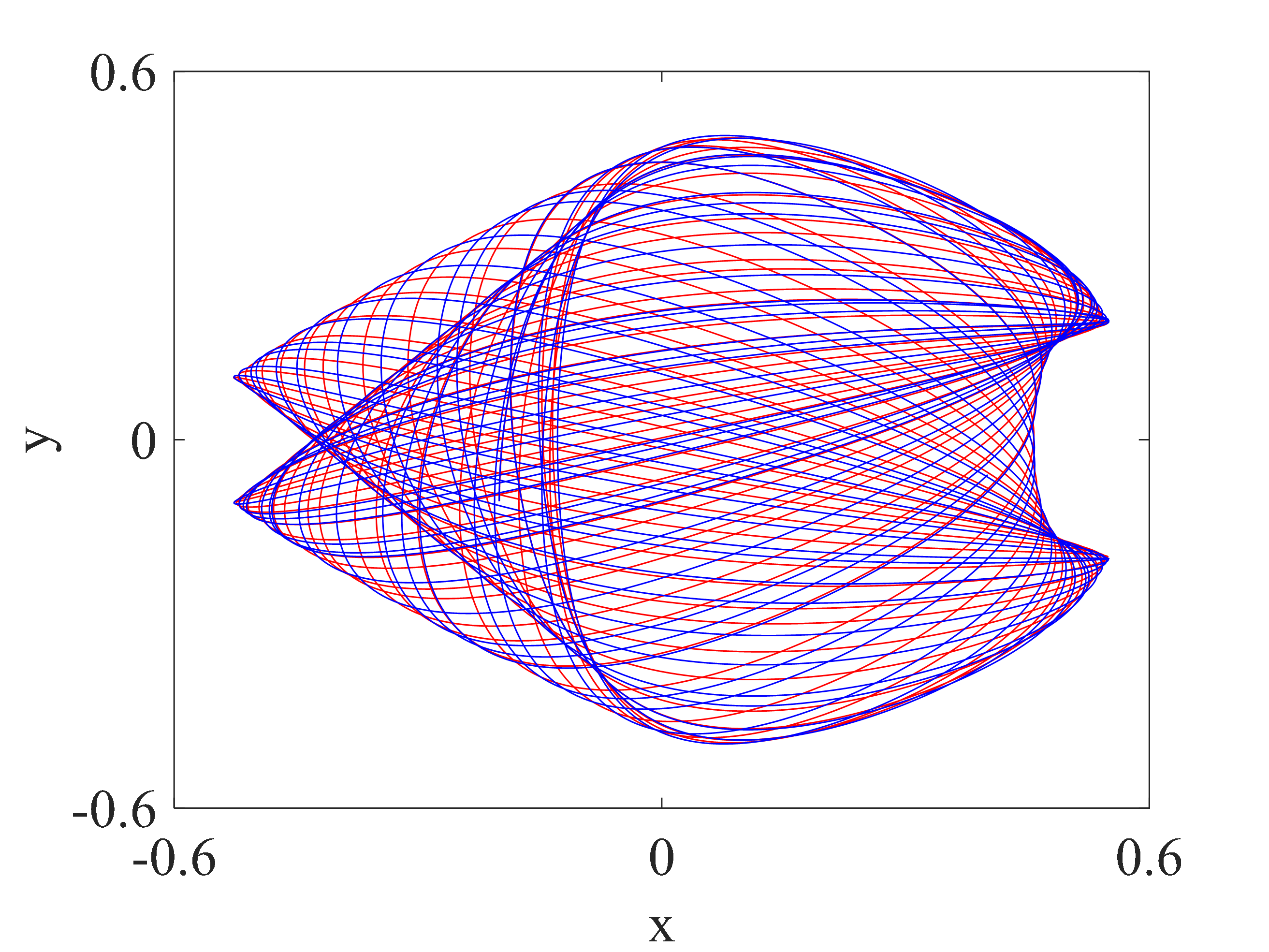}
			\caption{}\label{fig:two-exit-threshold-trajectories:d}
		\end{subfigure}\par\smallskip
		\begin{subfigure}[t]{.47\linewidth}
			\centering
			\includegraphics[width=\linewidth]{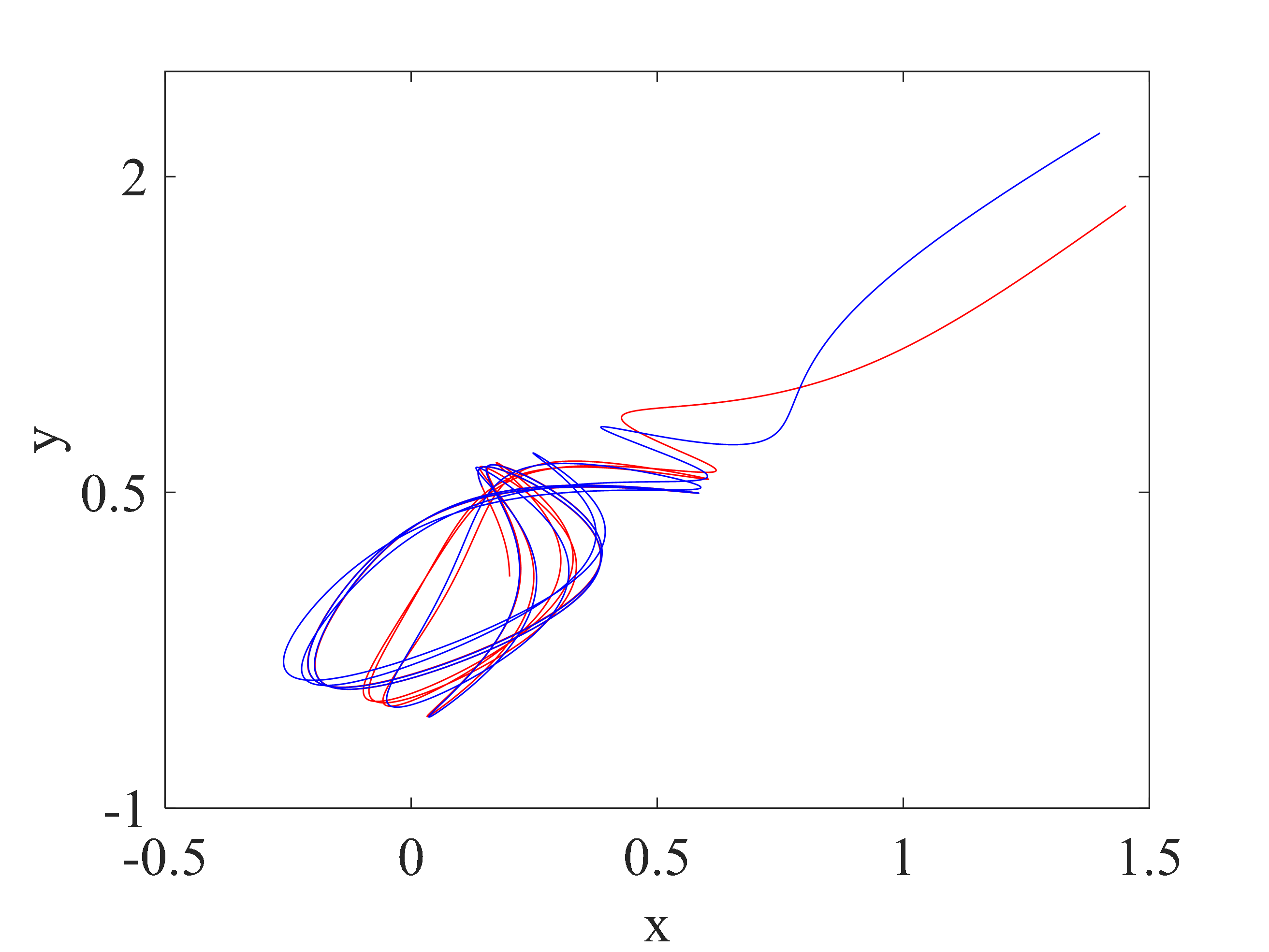}
			\caption{}\label{fig:two-exit-threshold-trajectories:e}
		\end{subfigure}\hfill
		\begin{subfigure}[t]{.47\linewidth}
			\centering
			\includegraphics[width=\linewidth]{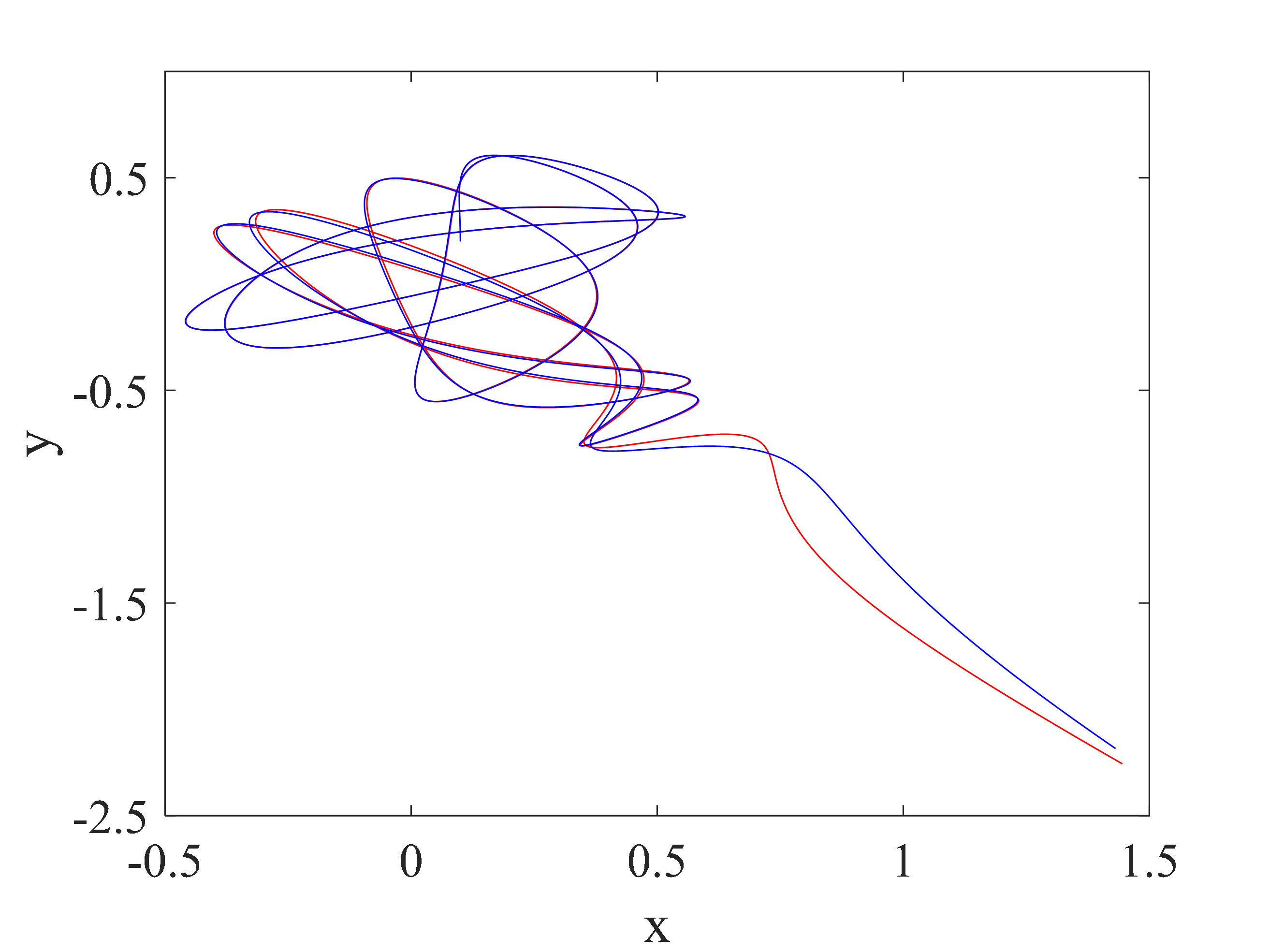}
			\caption{}\label{fig:two-exit-threshold-trajectories:f}
		\end{subfigure}
		\caption{Predicted escape and non-escape structures ((\subref{fig:two-exit-threshold-trajectories:a})--(\subref{fig:two-exit-threshold-trajectories:b}), gray: non-escape; green: escape through exit 1; red: escape through exit 2) and examples of predicted trajectories under different energy values and initial conditions ((\subref{fig:two-exit-threshold-trajectories:c})--(\subref{fig:two-exit-threshold-trajectories:f}), blue: true; red: predicted). (\subref{fig:two-exit-threshold-trajectories:a}) Predicted escape and non-escape structure at $E=0.1242$. (\subref{fig:two-exit-threshold-trajectories:b}) Predicted escape and non-escape structure at $E=0.12425$. (\subref{fig:two-exit-threshold-trajectories:c}) $E=E_{4}=0.12$, $\boldsymbol{z}_{4}=(0.1,\allowbreak  0.1,\allowbreak  0,\allowbreak  g(0.1,\allowbreak 0.1,\allowbreak E_{4}))$, trapped. (\subref{fig:two-exit-threshold-trajectories:d}) $E=E_{5}=0.15$, $\boldsymbol{z}_{5}=(-0.2,\allowbreak  -0.1,\allowbreak  0,\allowbreak  g(-0.2,\allowbreak -0.1,\allowbreak E_{5}))$, trapped. (\subref{fig:two-exit-threshold-trajectories:e}) $E=E_{5}=0.15$, $\boldsymbol{z}_{6}=(0.2,\allowbreak  0.1,\allowbreak  0,\allowbreak  g(0.2,\allowbreak 0.1,\allowbreak E_{5}))$, escape through exit 1. (\subref{fig:two-exit-threshold-trajectories:f}) $E=E_{5}=0.15$, $\boldsymbol{z}_{7}=(0.1,\allowbreak  0.2,\allowbreak  0,\allowbreak   g(0.1,\allowbreak 0.2,\allowbreak E_{5}))$, escape through exit 2.}
		\label{fig:two-exit-threshold-trajectories}
	\end{figure}

	Table~\ref{tab:threshold-1} presents the escape predictions. Here, a checkmark ($\checkmark$) denotes that at least one predicted trajectory is classified as escaping, whereas a cross ($\times$) indicates that no escape is detected within the prediction horizon. As shown in Table~\ref{tab:threshold-1}, for $E\leqslant{}0.124$, no predicted escape is detected. In contrast, for every tested energy value $E\geqslant{}0.125$, at least one predicted trajectory is classified as escaping, confirming that the model's estimated escape threshold lies between 0.124 and 0.125. To pinpoint the threshold more precisely, we conducted a finer scan over the interval $E\in{}[0.12405,\allowbreak 0.125]$ with a step size of $5\times{}10^{-5}$, and the results are summarized in Table~\ref{tab:threshold-2}. For the energy values up to and including 0.1242, no escaping trajectory is observed. Fig.~\ref{fig:two-exit-threshold-trajectories}(\subref{fig:two-exit-threshold-trajectories:a})--(\subref{fig:two-exit-threshold-trajectories:b}) present the predicted escape and non-escape structures at $E=0.1242$ and $E=0.12425$, respectively. Consequently, the model gives a grid-based threshold estimate of 0.12425. This estimated value differs from the theoretical threshold of 0.125 by only 0.00075.
	
	Denote $E=E_{4}=0.12$, and $E=E_{5}=0.15$. Define $\boldsymbol{z}_{4}=(0.1,\allowbreak  0.1,\allowbreak  0,\allowbreak  g(0.1,\allowbreak 0.1,\allowbreak E_{4}))$, $\boldsymbol{z}_{5}=(-0.2,\allowbreak  -0.1,\allowbreak  0,\allowbreak  g(-0.2,\allowbreak -0.1,\allowbreak E_{5}))$, $\boldsymbol{z}_{6}=(0.2,\allowbreak  0.1,\allowbreak  0,\allowbreak  g(0.2,\allowbreak 0.1,\allowbreak E_{5}))$, and $\boldsymbol{z}_{7}=(0.1,\allowbreak  0.2,\allowbreak  0,\allowbreak   g(0.1,\allowbreak 0.2,\allowbreak E_{5}))$. Fig.~\ref{fig:two-exit-threshold-trajectories}(\subref{fig:two-exit-threshold-trajectories:c})--(\subref{fig:two-exit-threshold-trajectories:f}) compare the predicted and true trajectories in the $(x,\allowbreak y)$ plane under initial values $\boldsymbol{z}_{4}$, $\boldsymbol{z}_{5}$, $\boldsymbol{z}_{6}$, and $\boldsymbol{z}_{7}$. For the cases $\boldsymbol{z}_{4}$ and $\boldsymbol{z}_{5}$, the system remains in a trapped state, and the predicted orbit projections agree with the reference projections. Meanwhile, for $\boldsymbol{z}_{6}$ and $\boldsymbol{z}_{7}$, the system escapes to infinity, specifically through exit 1 for $\boldsymbol{z}_{6}$ and through exit 2 for $\boldsymbol{z}_{7}$. The model accurately captures not only the occurrence of escape but also the correct escape channel.
	
	To further investigate the escape dynamics, we predict the structure of the escape channels in the $(x,\allowbreak y)$ plane. The predictions are performed over three distinct energy ranges, each with a corresponding domain for the initial value $(x_{0},\allowbreak y_{0})$. Specifically, for $0.13\leqslant{}E\leqslant{}0.21$, we sample $(x_{0},\allowbreak y_{0})$ from the square $[-1,\allowbreak 1]\times{}[-1,\allowbreak 1]$; for $0.22\leqslant{}E\leqslant{}0.3$, the sampling domain is enlarged to $[-1.2,\allowbreak 1.2]\times{}[-1.2,\allowbreak 1.2]$; and for $0.31\leqslant{}E\leqslant{}0.4$, the domain is further extended to $[-1.4,\allowbreak 1.4]\times{}[-1.4,\allowbreak 1.4]$. In all cases, the domain is discretized into a uniform grid with a step size of 0.01, and every admissible grid point is used as an initial value. The momentum $p_{x0}$ is fixed at 0 for all trajectories, while $p_{y0}$ is determined by the function $g(x,\allowbreak y,\allowbreak E)$, which ensures that the total energy equals the specified value $E$.
	
	\begin{figure}[!htbp]
		\centering
		\begin{subfigure}[t]{.31\linewidth}
			\centering
			\includegraphics[width=\linewidth]{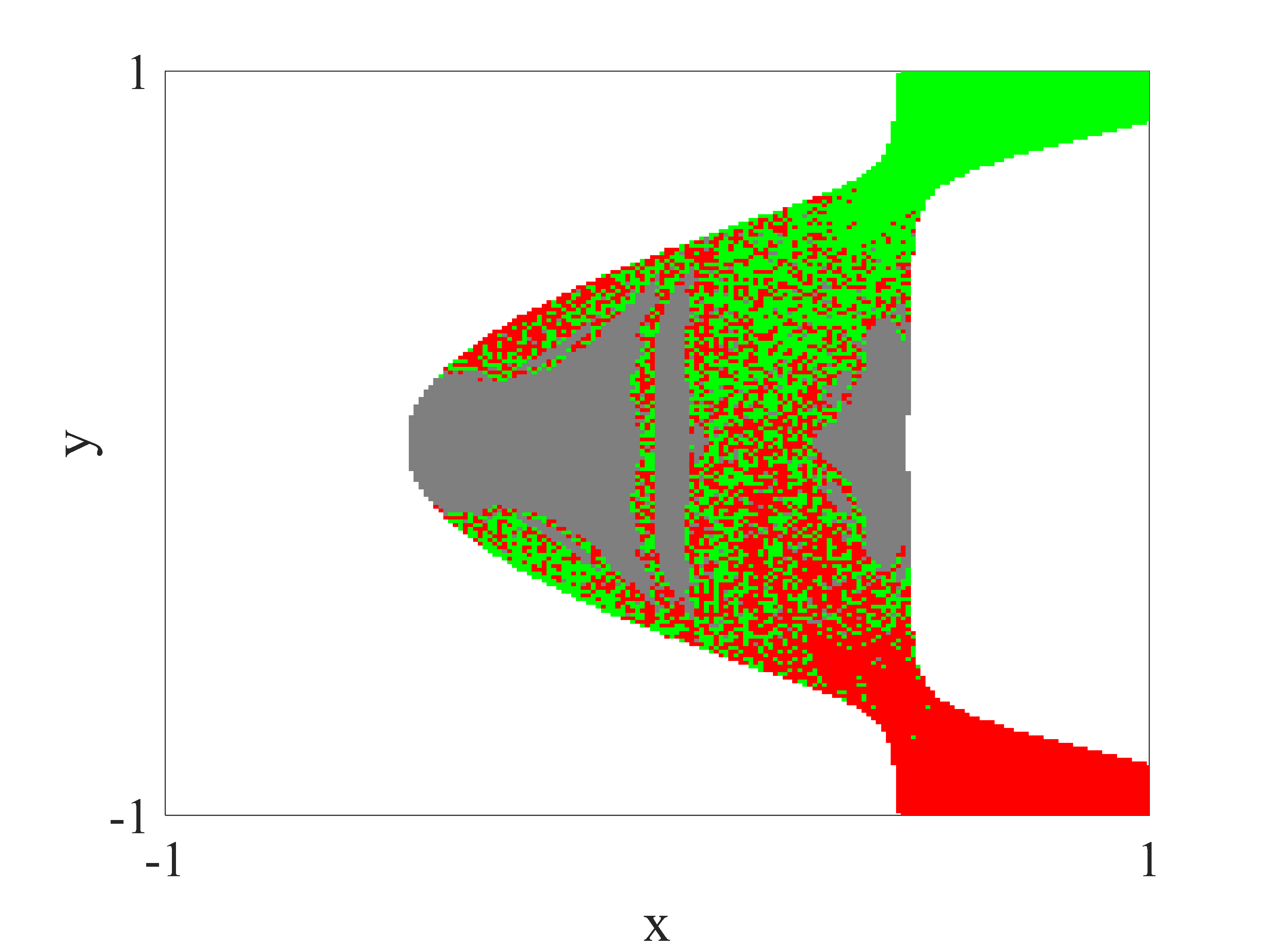}
			\caption{}\label{fig:two-exit-basins:a}
		\end{subfigure}\hfill
		\begin{subfigure}[t]{.31\linewidth}
			\centering
			\includegraphics[width=\linewidth]{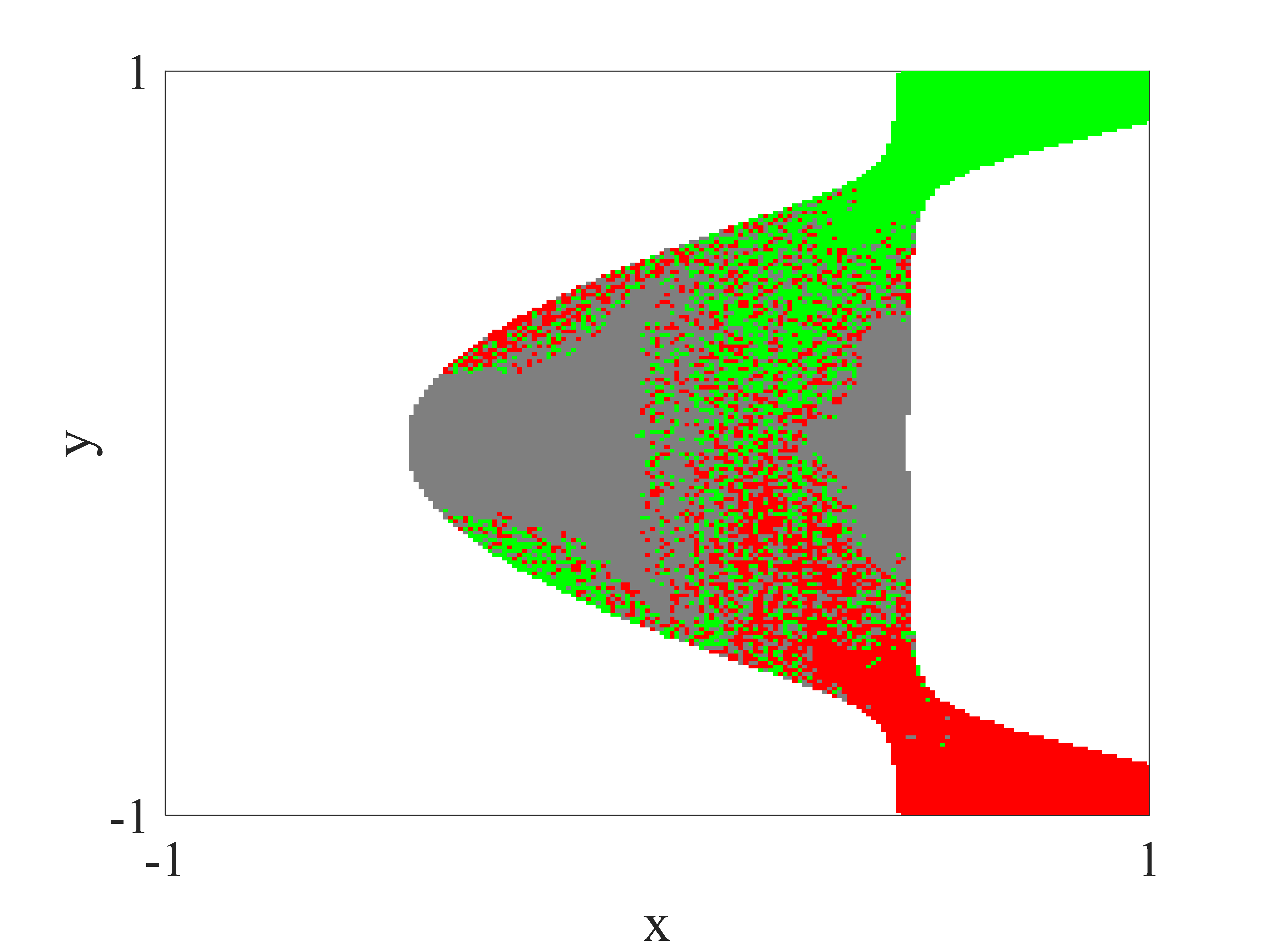}
			\caption{}\label{fig:two-exit-basins:b}
		\end{subfigure}\hfill
		\begin{subfigure}[t]{.31\linewidth}
			\centering
			\includegraphics[width=\linewidth]{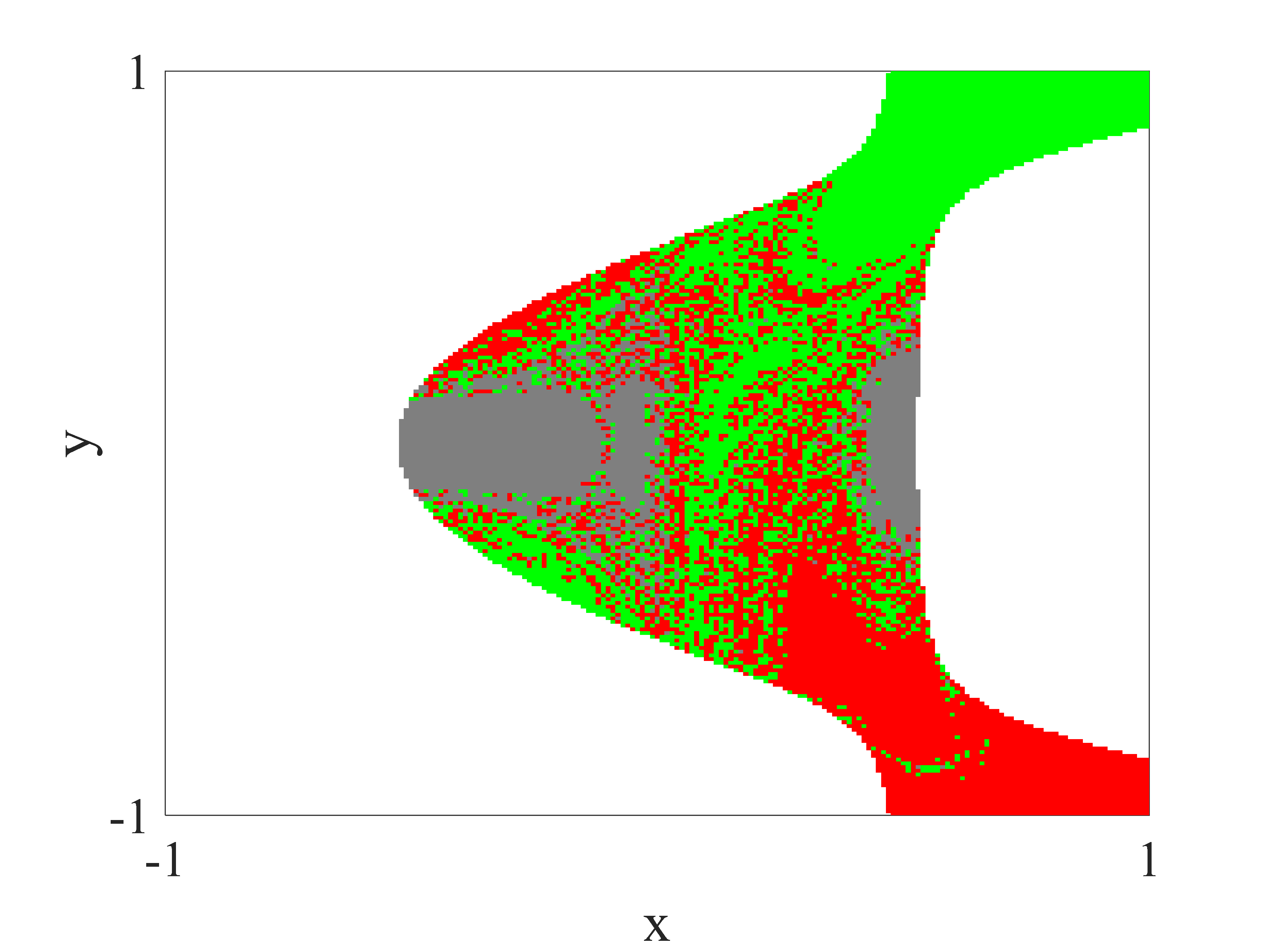}
			\caption{}\label{fig:two-exit-basins:c}
		\end{subfigure}\par\smallskip
		\begin{subfigure}[t]{.31\linewidth}
			\centering
			\includegraphics[width=\linewidth]{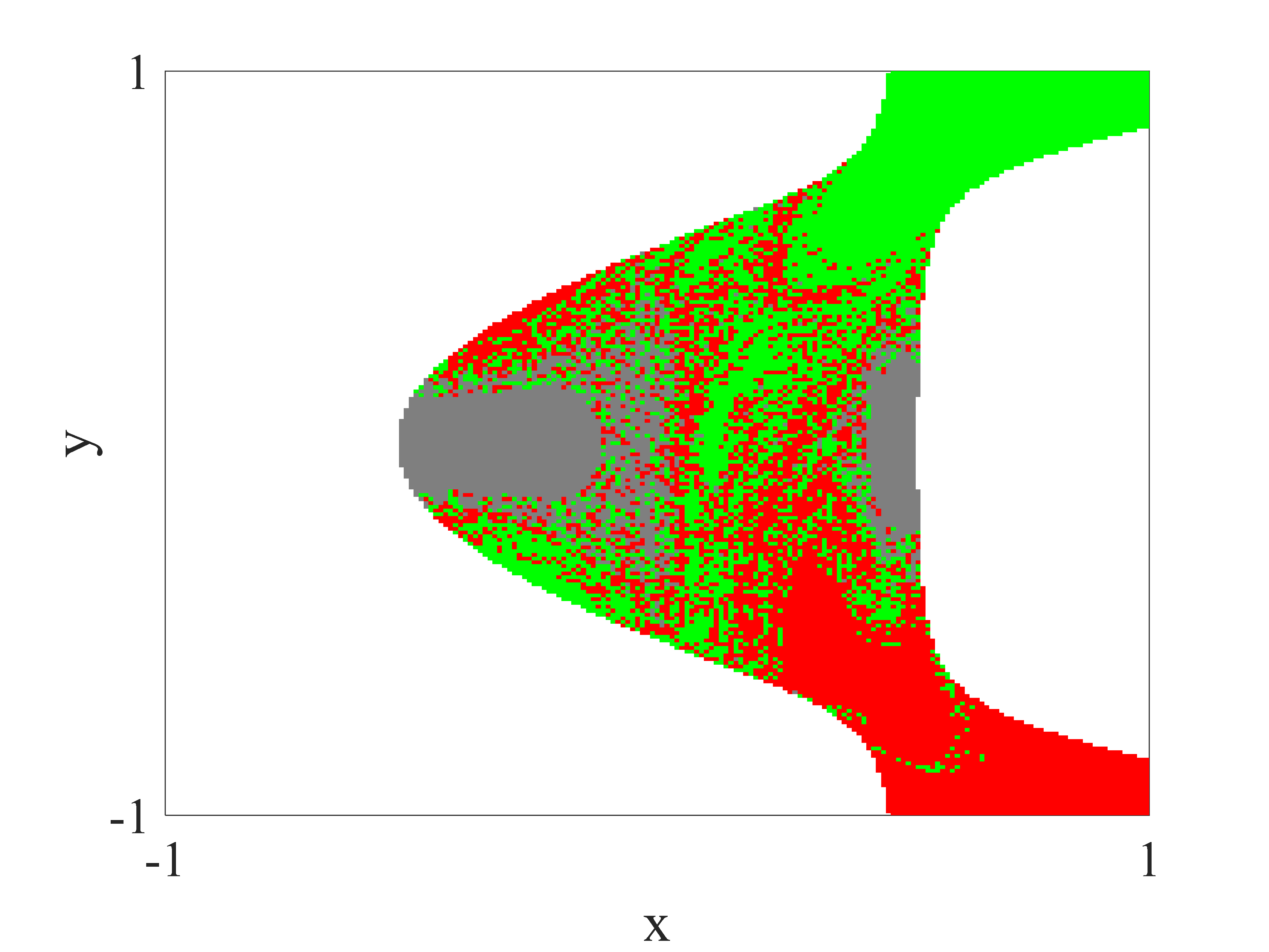}
			\caption{}\label{fig:two-exit-basins:d}
		\end{subfigure}\hfill
		\begin{subfigure}[t]{.31\linewidth}
			\centering
			\includegraphics[width=\linewidth]{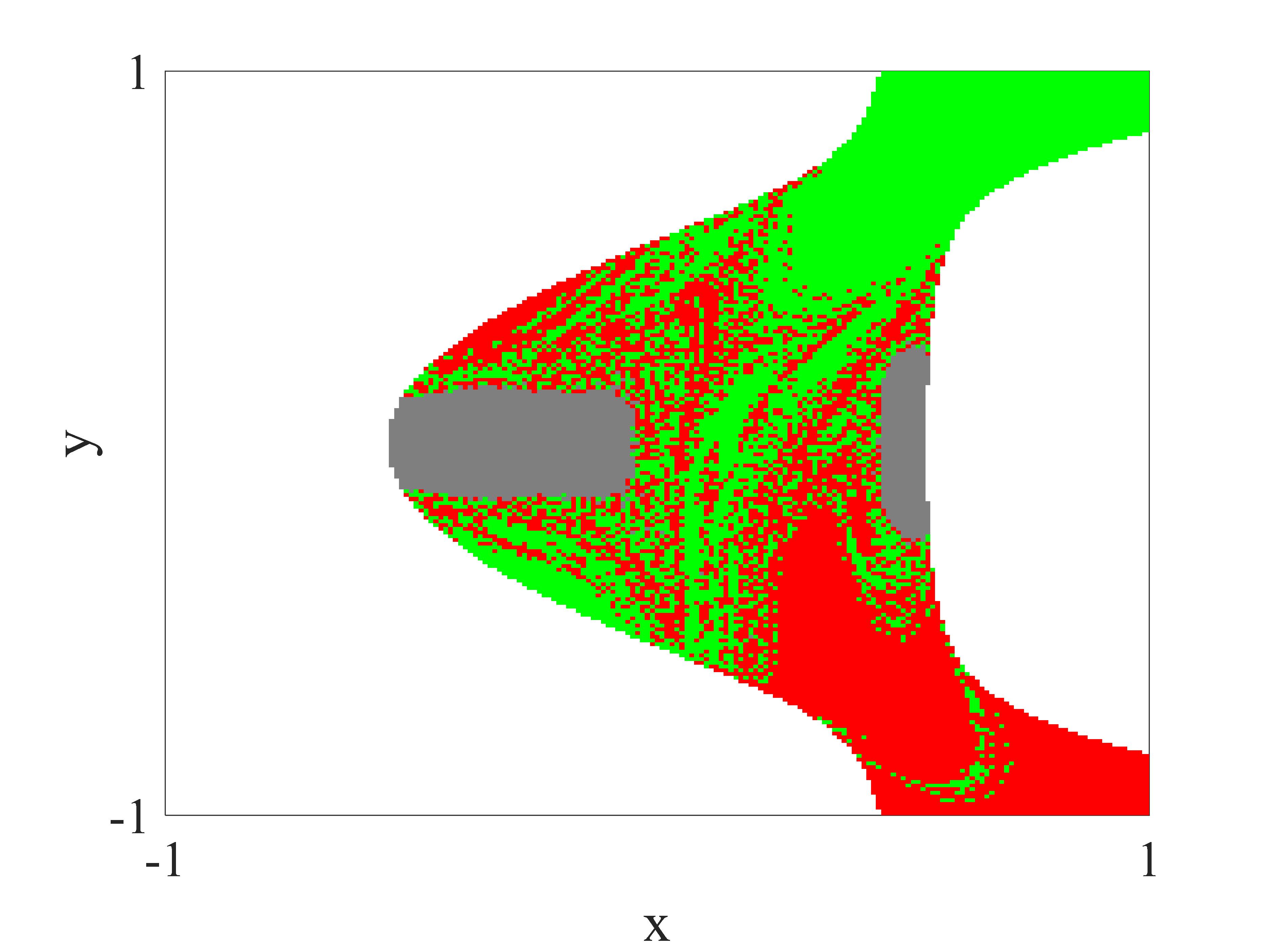}
			\caption{}\label{fig:two-exit-basins:e}
		\end{subfigure}\hfill
		\begin{subfigure}[t]{.31\linewidth}
			\centering
			\includegraphics[width=\linewidth]{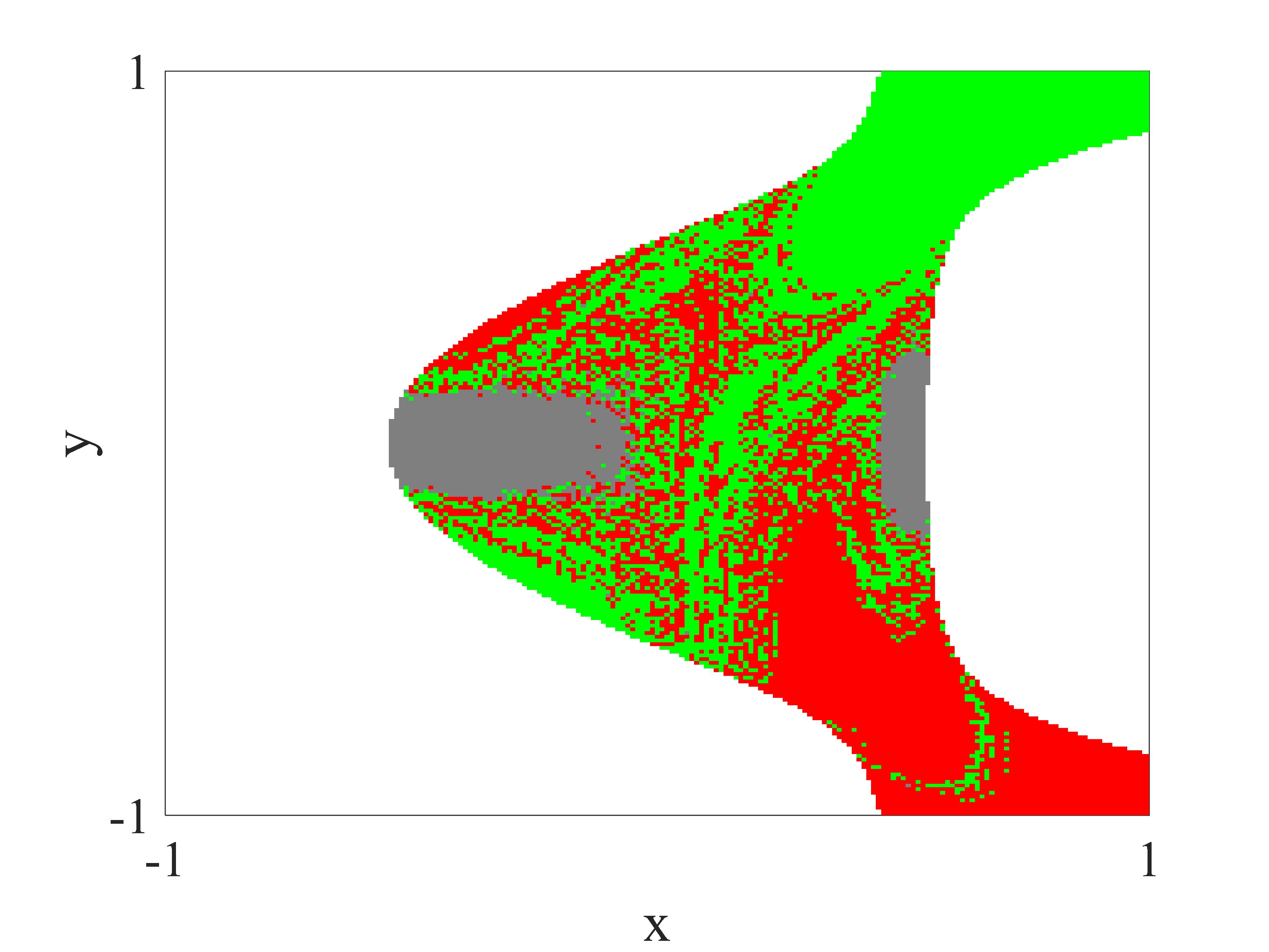}
			\caption{}\label{fig:two-exit-basins:f}
		\end{subfigure}
		\caption{Predictions of the structure of the escape channels in the $(x,\allowbreak y)$ plane for different energy values. Gray: trapped; green: escape through exit 1; red: escape through exit 2. (\subref{fig:two-exit-basins:a}) $E=0.13$, true. (\subref{fig:two-exit-basins:b}) $E=0.13$, predicted. (\subref{fig:two-exit-basins:c}) $E=0.14$, true. (\subref{fig:two-exit-basins:d}) $E=0.14$, predicted. (\subref{fig:two-exit-basins:e}) $E=0.15$, true. (\subref{fig:two-exit-basins:f}) $E=0.15$, predicted. (\subref{fig:two-exit-basins:g}) $E=0.16$, true. (\subref{fig:two-exit-basins:h}) $E=0.16$, predicted. (\subref{fig:two-exit-basins:i}) $E=0.17$, true. (\subref{fig:two-exit-basins:j}) $E=0.17$, predicted. (\subref{fig:two-exit-basins:k}) $E=0.18$, true. (\subref{fig:two-exit-basins:l}) $E=0.18$, predicted. (\subref{fig:two-exit-basins:m}) $E=0.20$, true. (\subref{fig:two-exit-basins:n}) $E=0.20$, predicted. (\subref{fig:two-exit-basins:o}) $E=0.3$, true. (\subref{fig:two-exit-basins:p}) $E=0.3$, predicted. (\subref{fig:two-exit-basins:q}) $E=0.4$, true. (\subref{fig:two-exit-basins:r}) $E=0.4$, predicted.}
		\label{fig:two-exit-basins}
	\end{figure}
	
	\begin{figure}[!htbp]
		\centering
		\ContinuedFloat
		\begin{subfigure}[t]{.31\linewidth}
			\centering
			\includegraphics[width=\linewidth]{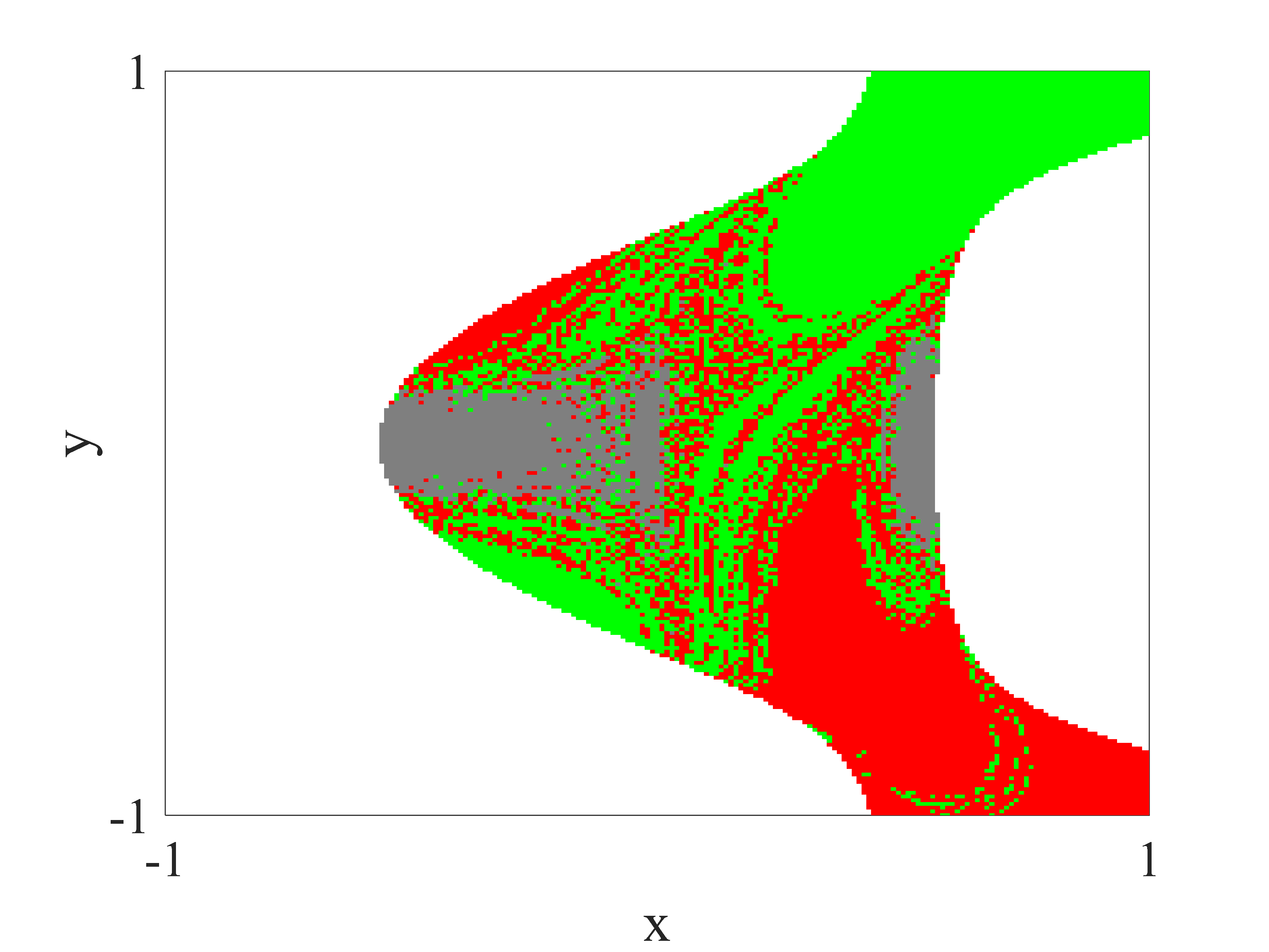}
			\caption{}\label{fig:two-exit-basins:g}
		\end{subfigure}\hfill
		\begin{subfigure}[t]{.31\linewidth}
			\centering
			\includegraphics[width=\linewidth]{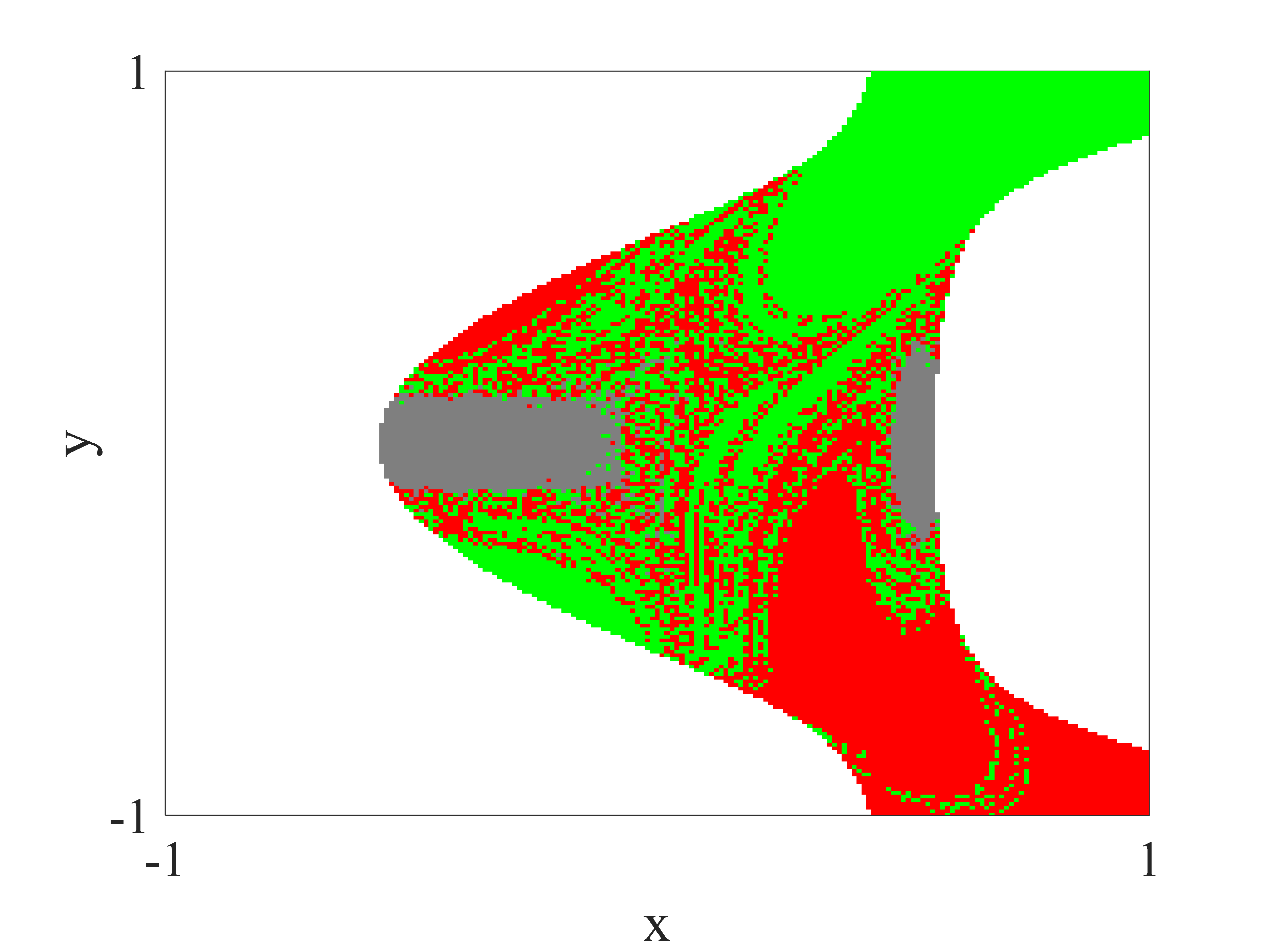}
			\caption{}\label{fig:two-exit-basins:h}
		\end{subfigure}\hfill
		\begin{subfigure}[t]{.31\linewidth}
			\centering
			\includegraphics[width=\linewidth]{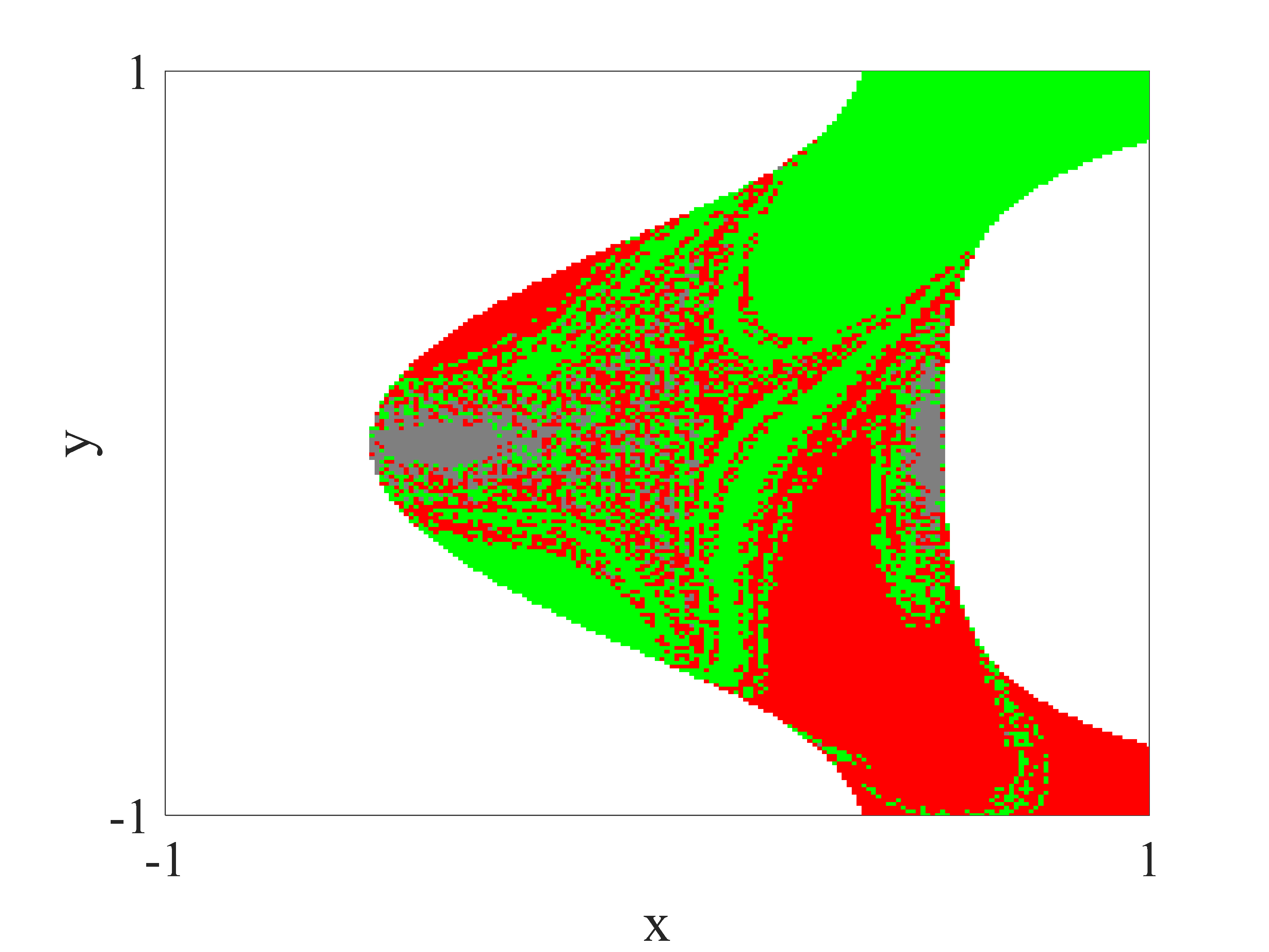}
			\caption{}\label{fig:two-exit-basins:i}
		\end{subfigure}\par\smallskip
		\begin{subfigure}[t]{.31\linewidth}
			\centering
			\includegraphics[width=\linewidth]{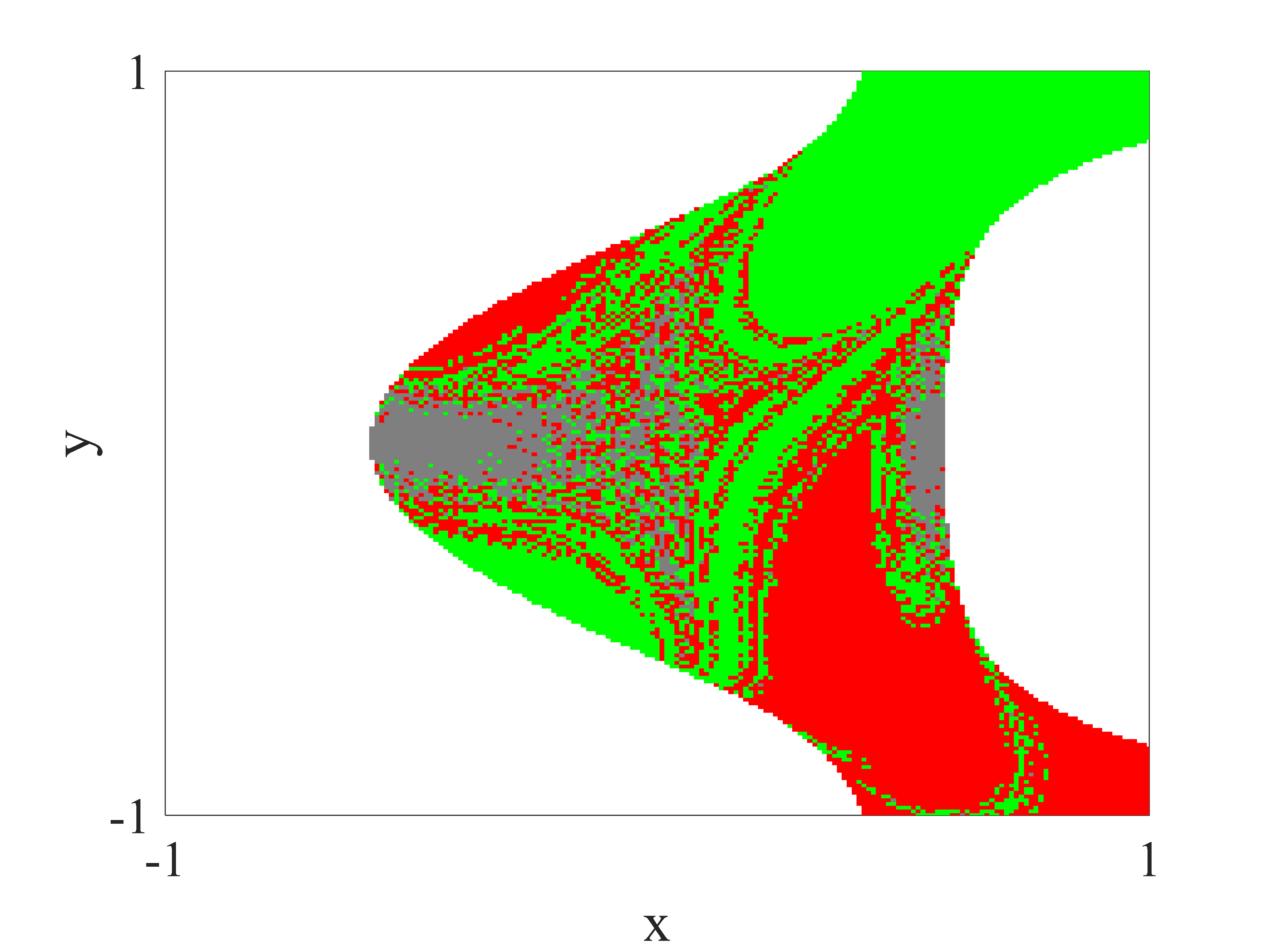}
			\caption{}\label{fig:two-exit-basins:j}
		\end{subfigure}\hfill
		\begin{subfigure}[t]{.31\linewidth}
			\centering
			\includegraphics[width=\linewidth]{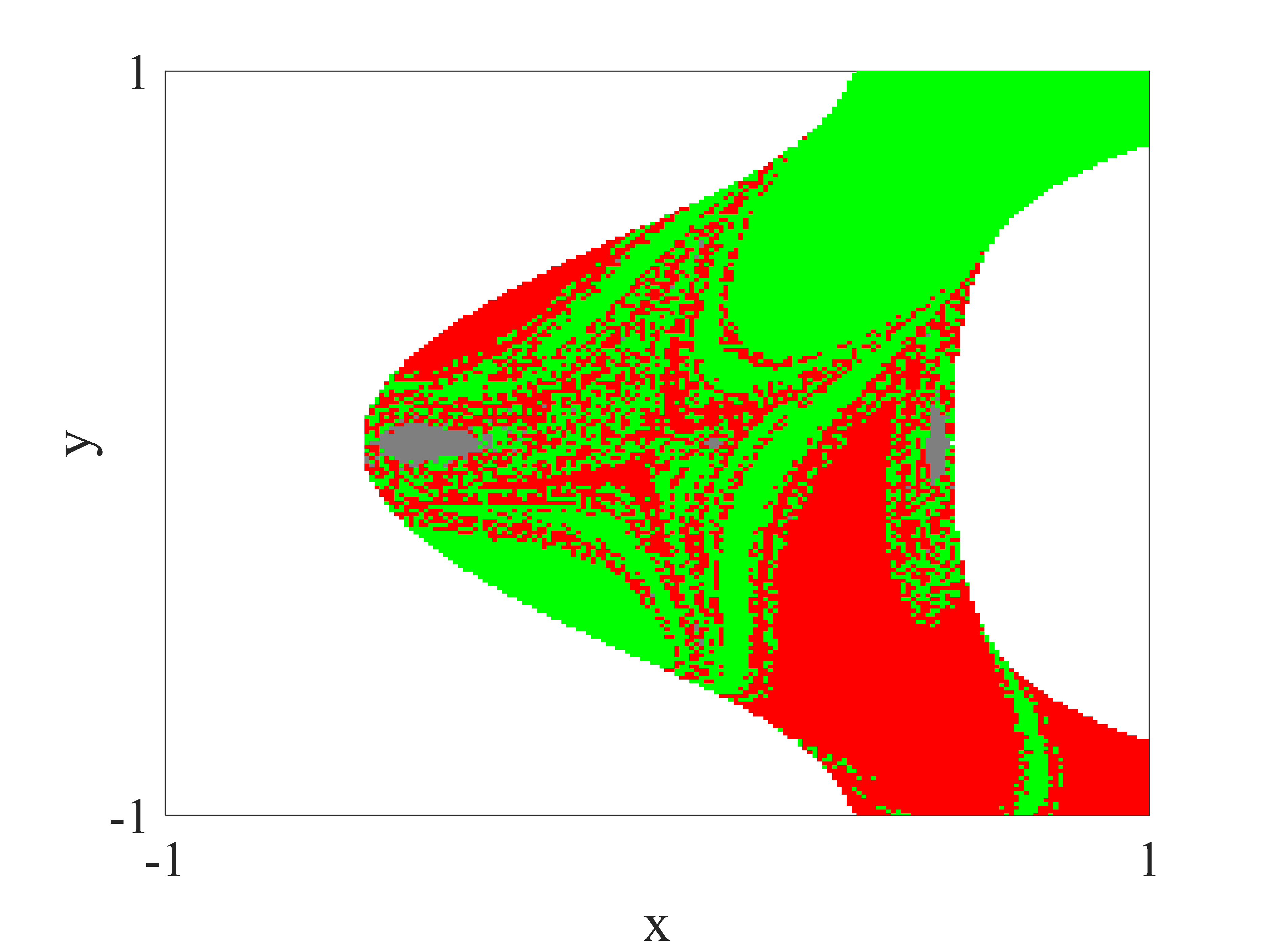}
			\caption{}\label{fig:two-exit-basins:k}
		\end{subfigure}\hfill
		\begin{subfigure}[t]{.31\linewidth}
			\centering
			\includegraphics[width=\linewidth]{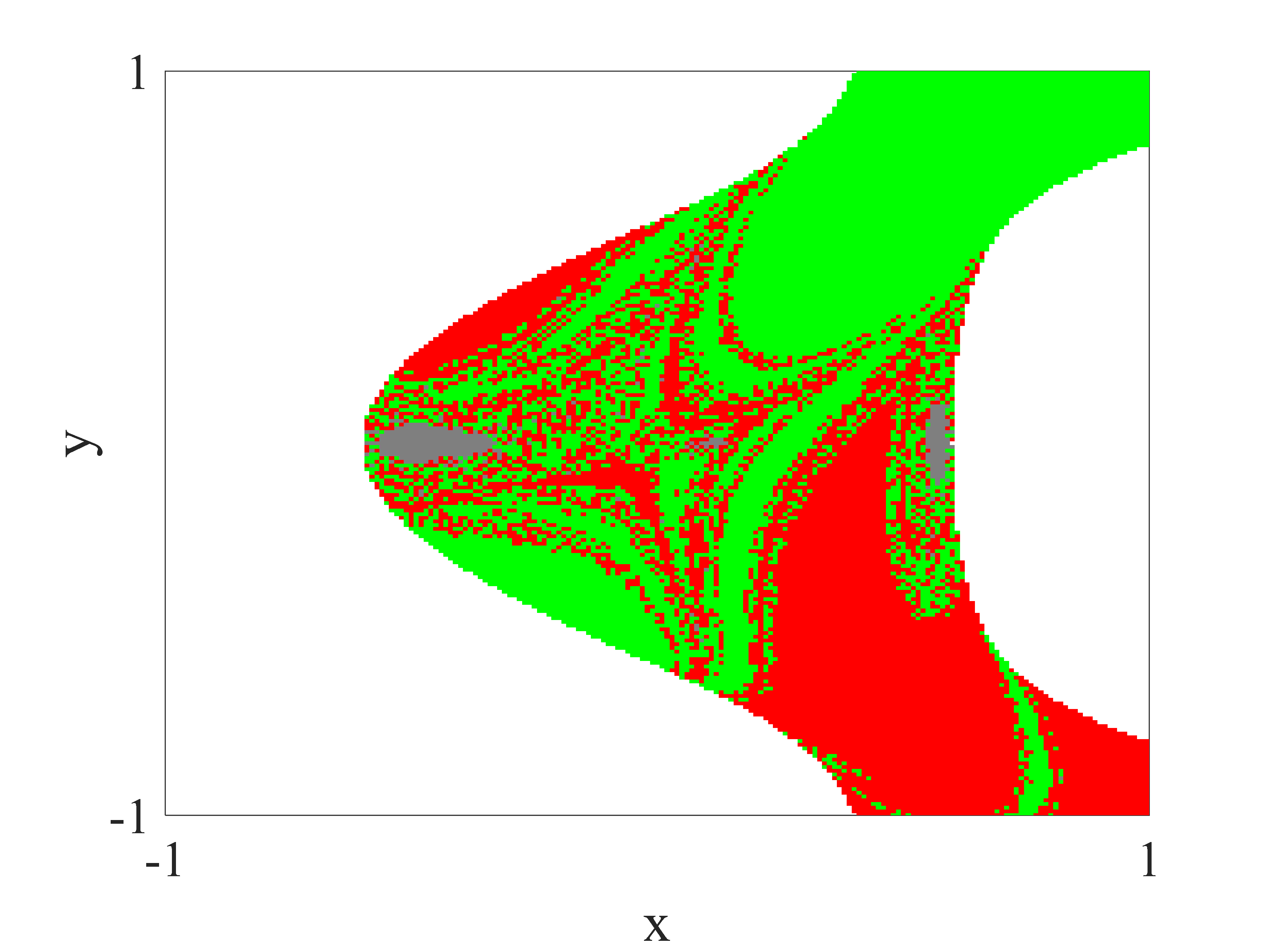}
			\caption{}\label{fig:two-exit-basins:l}
		\end{subfigure}\par\smallskip
		\begin{subfigure}[t]{.31\linewidth}
			\centering
			\includegraphics[width=\linewidth]{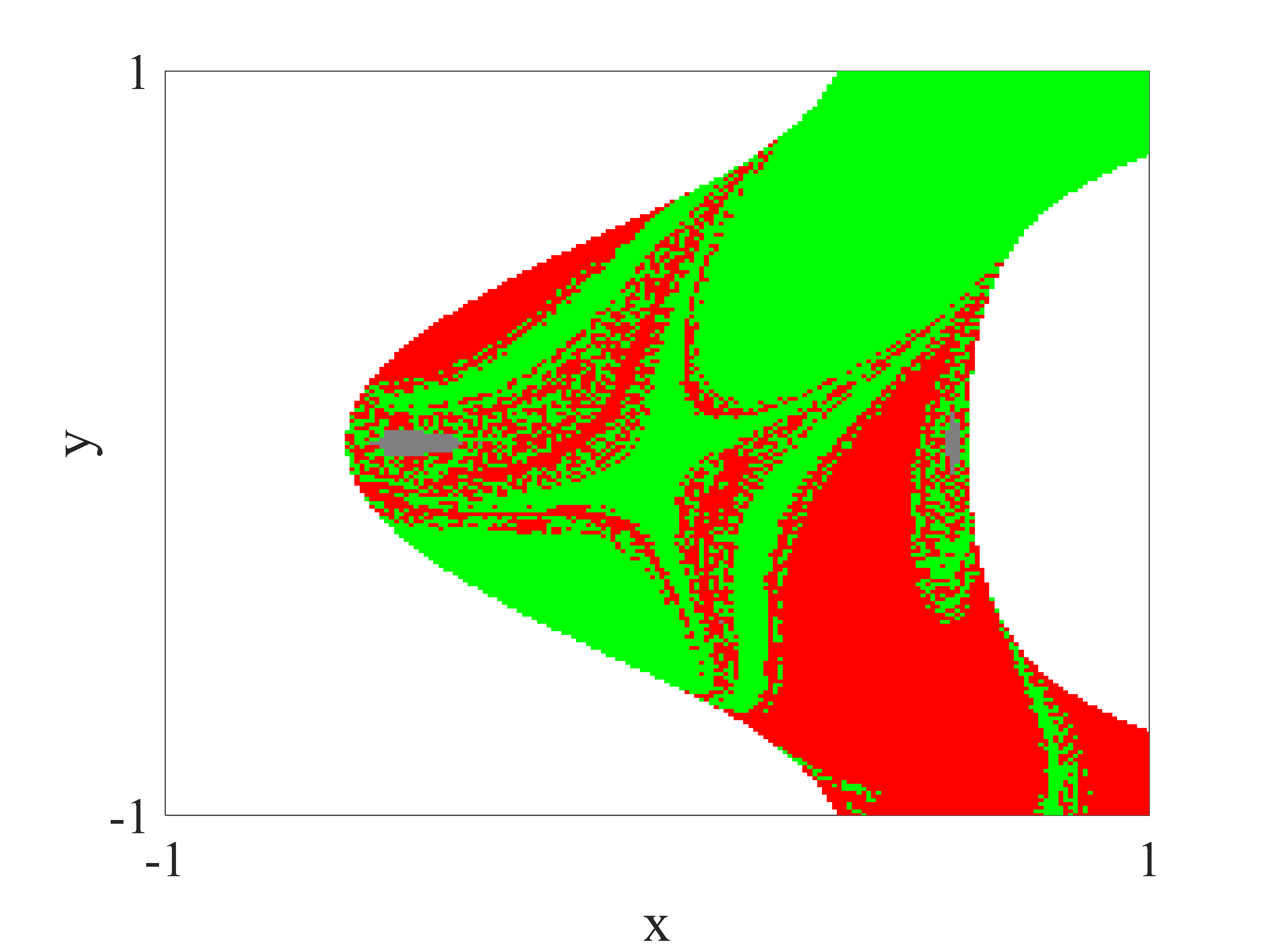}
			\caption{}\label{fig:two-exit-basins:m}
		\end{subfigure}\hfill
		\begin{subfigure}[t]{.31\linewidth}
			\centering
			\includegraphics[width=\linewidth]{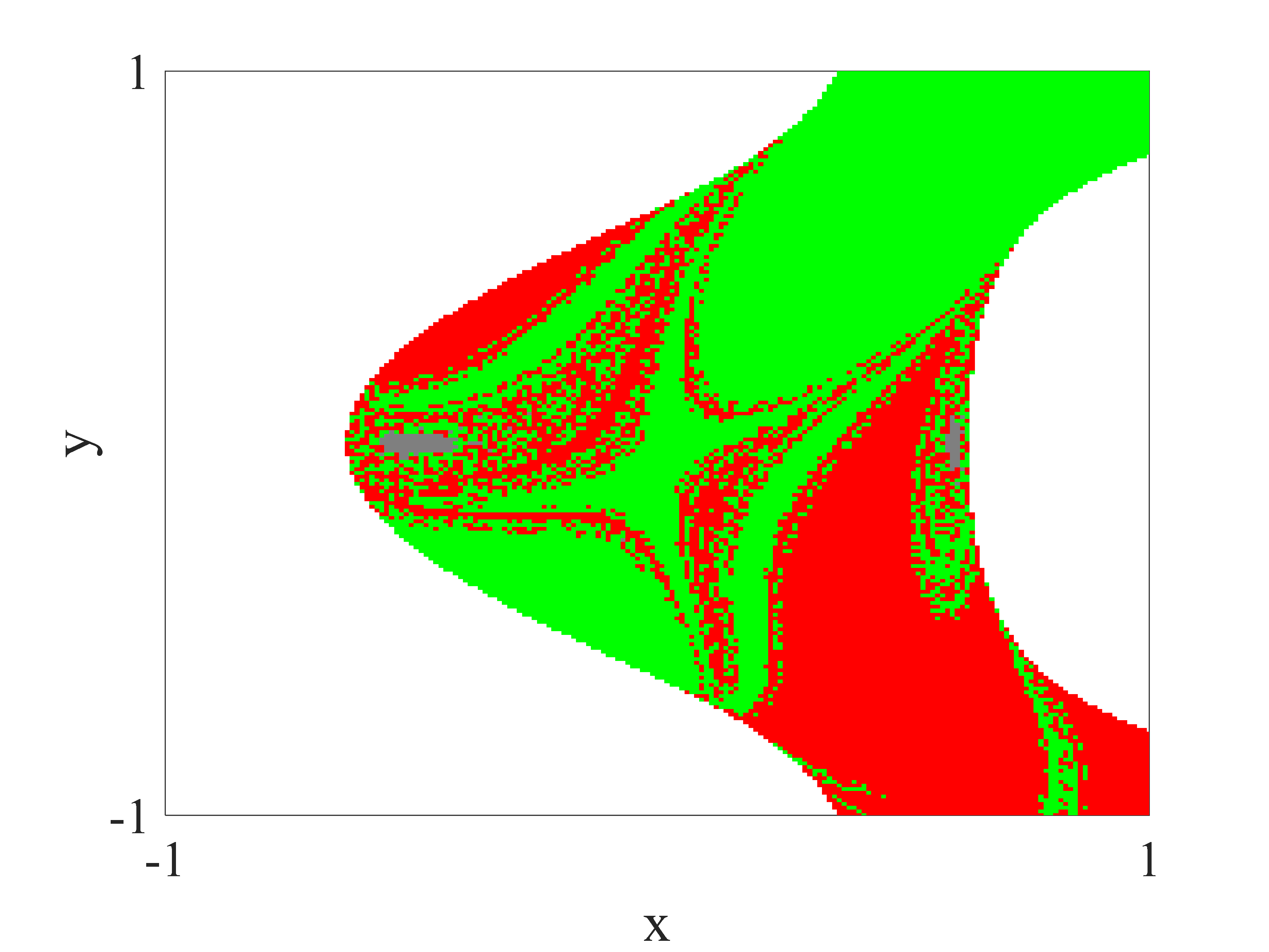}
			\caption{}\label{fig:two-exit-basins:n}
		\end{subfigure}\hfill
		\begin{subfigure}[t]{.31\linewidth}
			\centering
			\includegraphics[width=\linewidth]{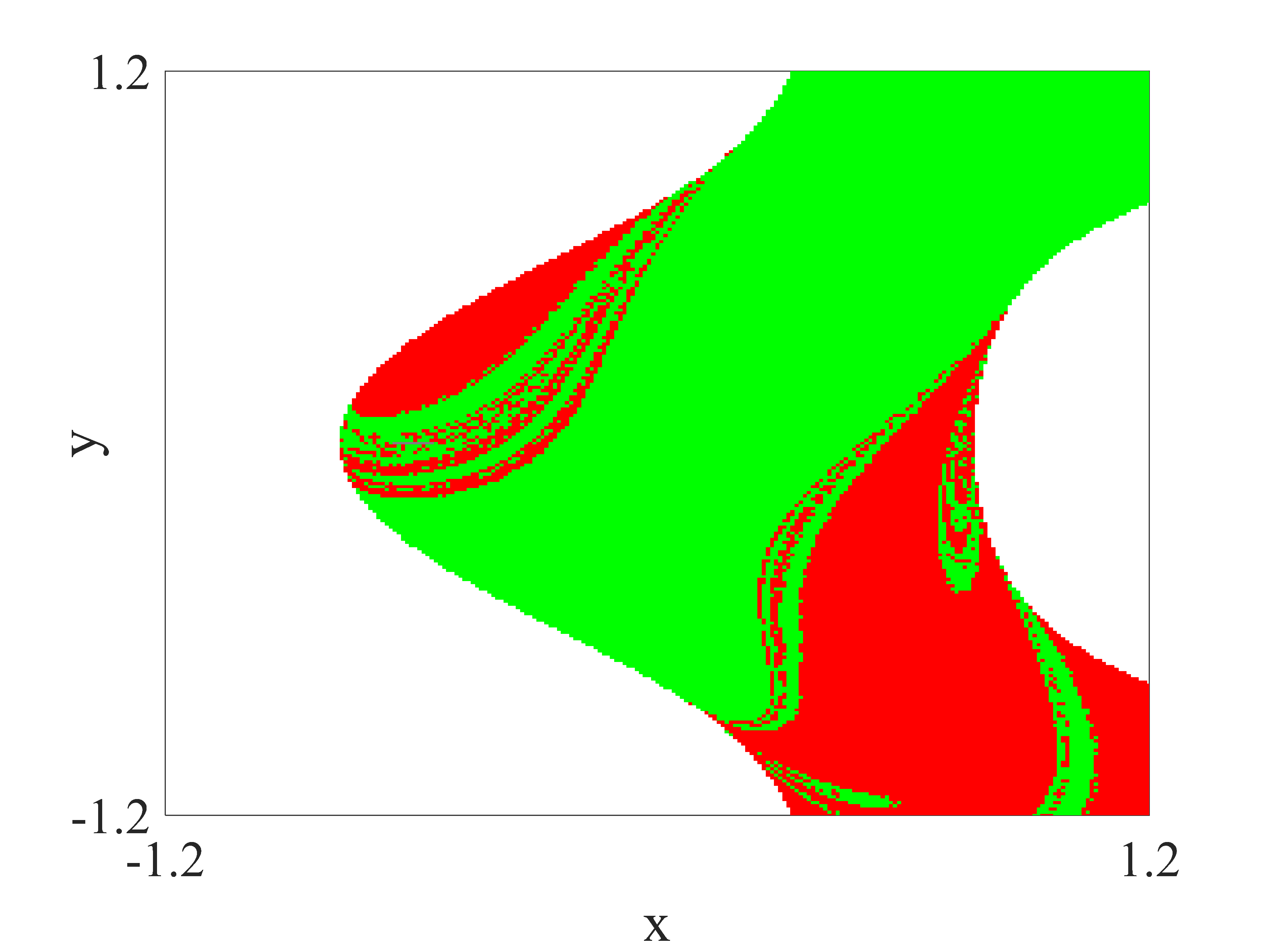}
			\caption{}\label{fig:two-exit-basins:o}
		\end{subfigure}\par\smallskip
		\begin{subfigure}[t]{.31\linewidth}
			\centering
			\includegraphics[width=\linewidth]{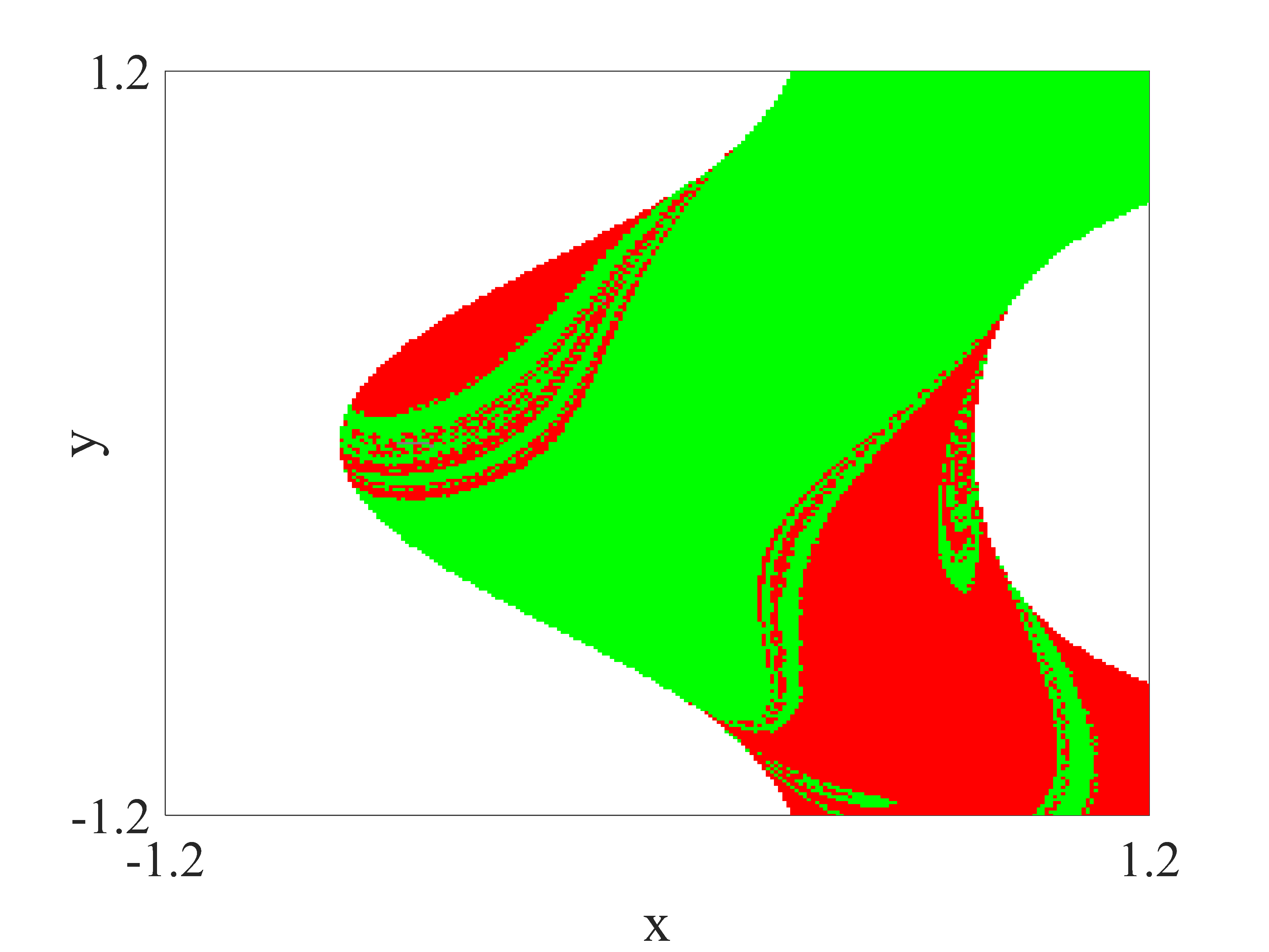}
			\caption{}\label{fig:two-exit-basins:p}
		\end{subfigure}\hfill
		\begin{subfigure}[t]{.31\linewidth}
			\centering
			\includegraphics[width=\linewidth]{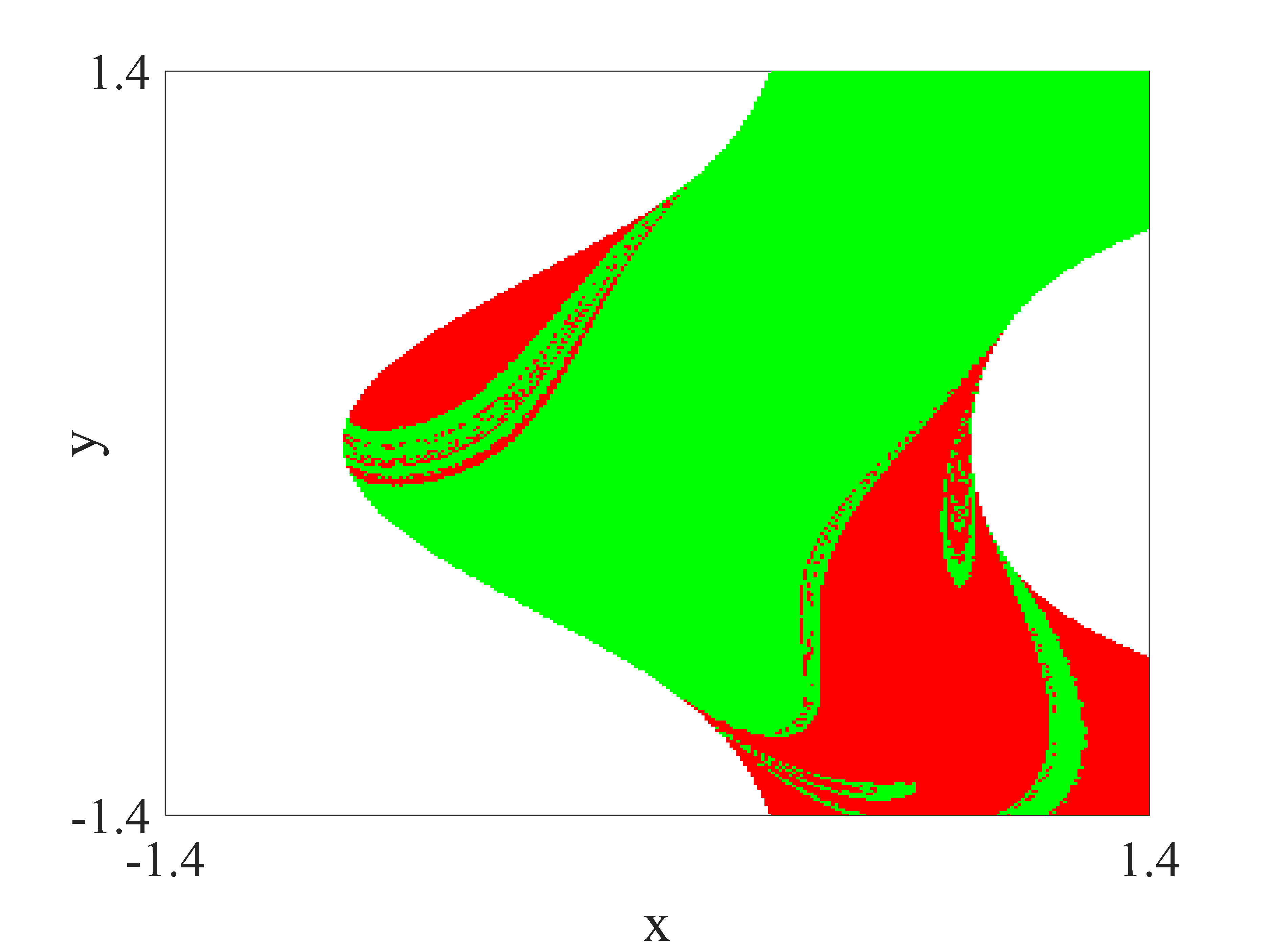}
			\caption{}\label{fig:two-exit-basins:q}
		\end{subfigure}\hfill
		\begin{subfigure}[t]{.31\linewidth}
			\centering
			\includegraphics[width=\linewidth]{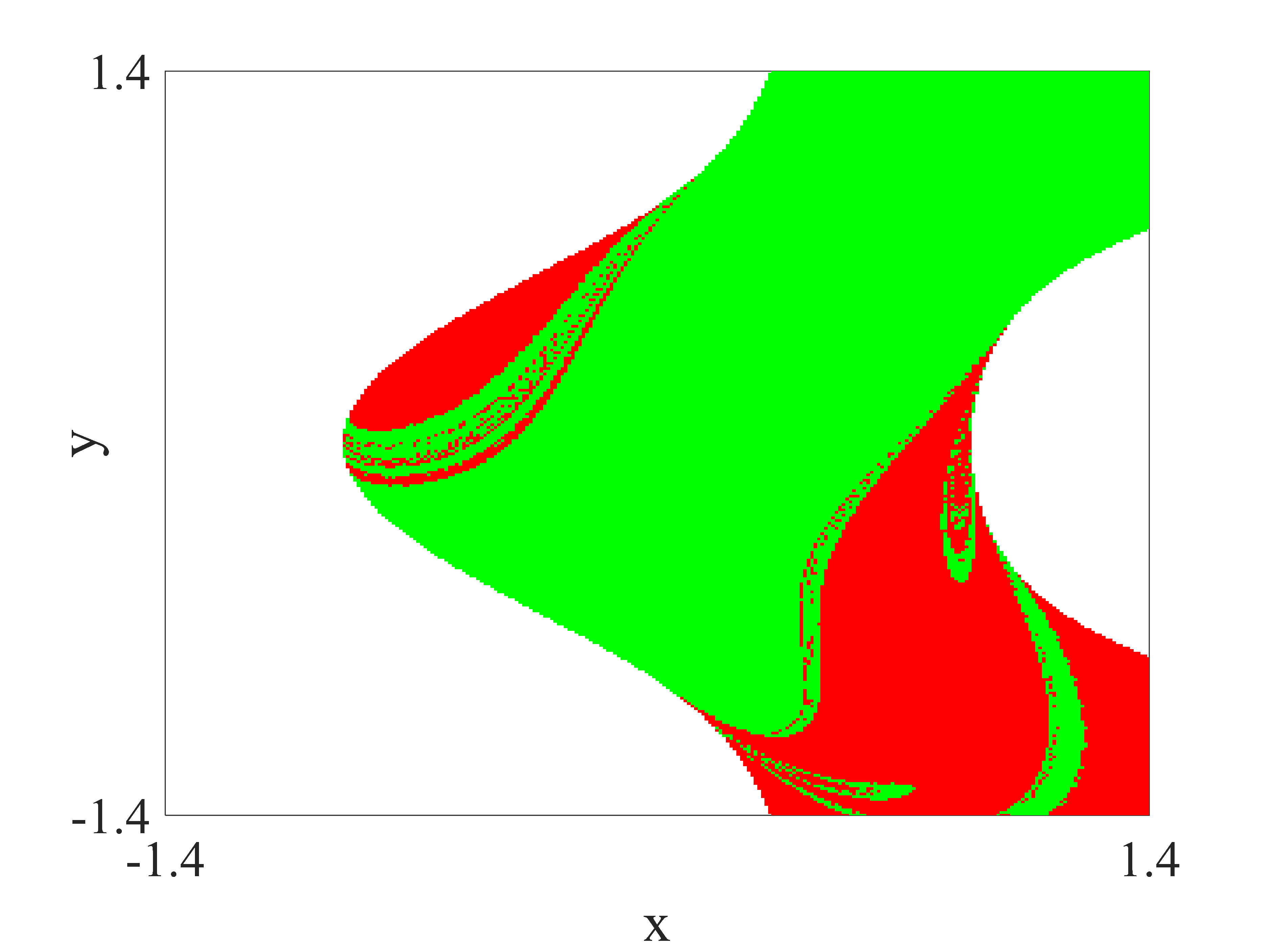}
			\caption{}\label{fig:two-exit-basins:r}
		\end{subfigure}
		\caption[]{\textit{Continued.}}
	\end{figure}
	
\begin{figure}[!htbp]
		\centering
		\includegraphics[width=.64\linewidth]{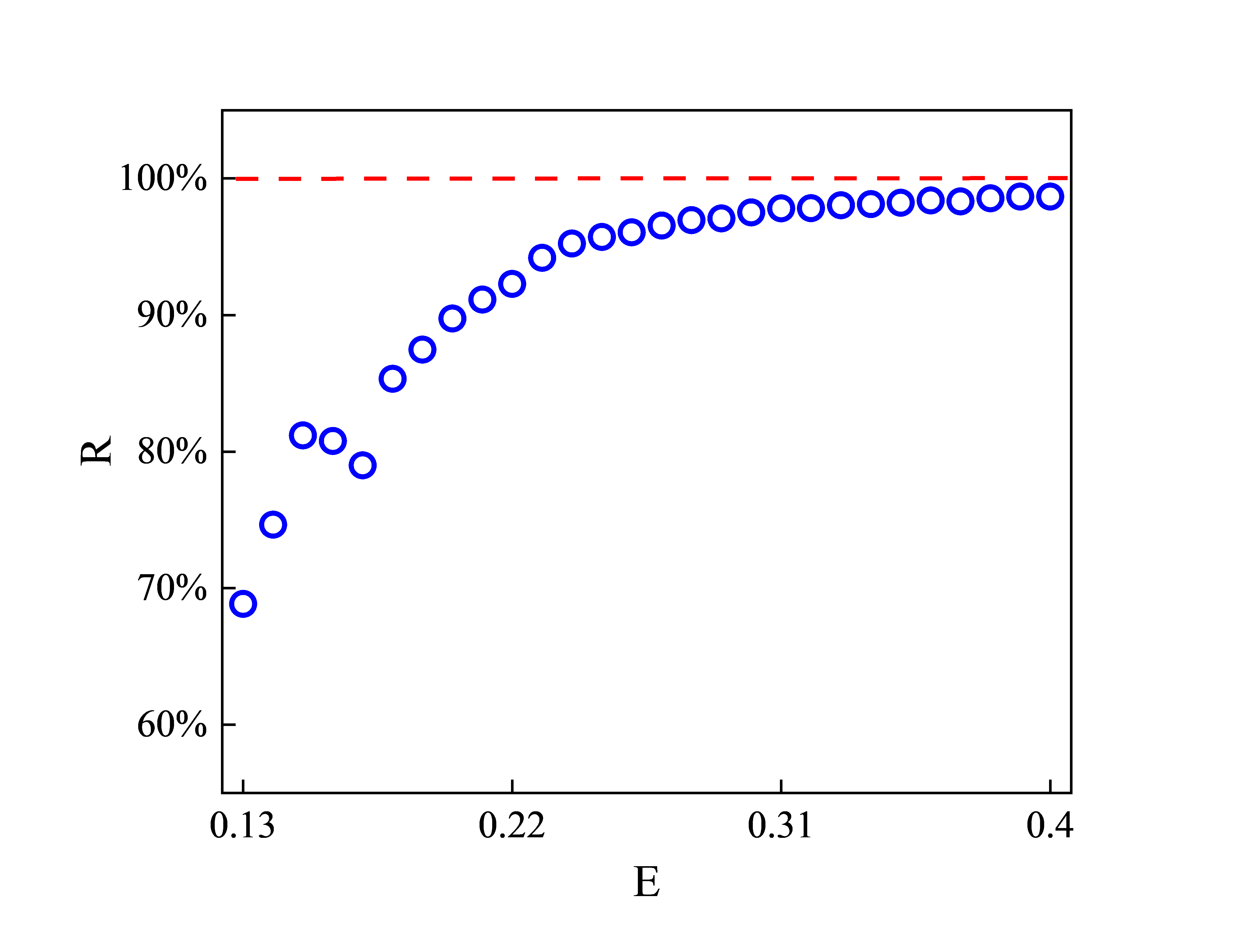}
		\caption{Prediction accuracy of the structure of the escape channels under different energy values.}
		\label{fig:two-exit-accuracy}
	\end{figure}

	Fig.~\ref{fig:two-exit-basins} presents a comparison between the predicted and true structures for energy values $E=0.13,\allowbreak 0.14,\allowbreak 0.15,\allowbreak 0.16,\allowbreak 0.17,\allowbreak 0.18,\allowbreak 0.2,\allowbreak 0.3$, and $0.4$. The true structure of the escape channels is obtained by numerically integrating the original system equations, while the predictions are generated by the trained model. For each initial condition, the model outputs the trajectory. We then classify whether it escapes and, if so, through which exit. The accuracy of the predictions is quantified by the ratio:
	
	\begin{equation*}
		R=\frac{M_{0}}{M},\allowbreak 
	\end{equation*}
	where $M$ denotes the total number of valid initial conditions, i.e., grid points lying within the specified domain and satisfying condition \eqref{eq:two-exit-admissibility}, and $M_{0}$ is the number of those points for which the predicted escape channel (or nonescape) agrees with the true result. Thus, $R$ represents the fraction of correctly predicted escapechannel assignments. Fig.~\ref{fig:two-exit-accuracy} displays the prediction accuracy $R$ at energy values ranging from 0.13 to 0.4 in increments of 0.01. Overall, the accuracy exhibits an upward trend. At the lowest energy considered, $E=0.13$, the accuracy is $68.8\%$. As energy increases, the escape dynamics become more pronounced and the channel boundaries become sharper, which may help explain the improved model performance. At $E=0.4$, the accuracy reaches $98.7\%$. Note that the training data are generated at energy values below the escape threshold, so the model is never exposed to an escaping trajectory during training and learns only the bounded dynamics. Near the threshold (e.g., $E=0.13$), the system is in a critical regime where the transition from trapped to escaping motion occurs. In this region, the escape behavior is extremely sensitive to both the energy value and the initial condition, and the boundaries separating different escape channels, or between escape and non-escape, can be highly intricate or even fractal. Since the model has no direct experience with escape events, its predictions here rely on extrapolation from sub-threshold dynamics, which limits its accuracy. As the energy increases to well above the threshold (e.g., $E=0.4$), the escape dynamics become less sensitive on the sampled grid. The potential barrier is small relative to the total energy, so trajectories escape readily, and the channel structure is governed mainly by the overall geometry of the potential landscape rather than by fine-scale sensitivity. This makes the escape-channel map smoother and more predictable, enabling the model to generalize effectively even though it was trained only on non-escape data.
	
\begin{figure}[!htbp]
		\centering
		\begin{subfigure}[t]{.47\linewidth}
			\centering
			\includegraphics[width=\linewidth]{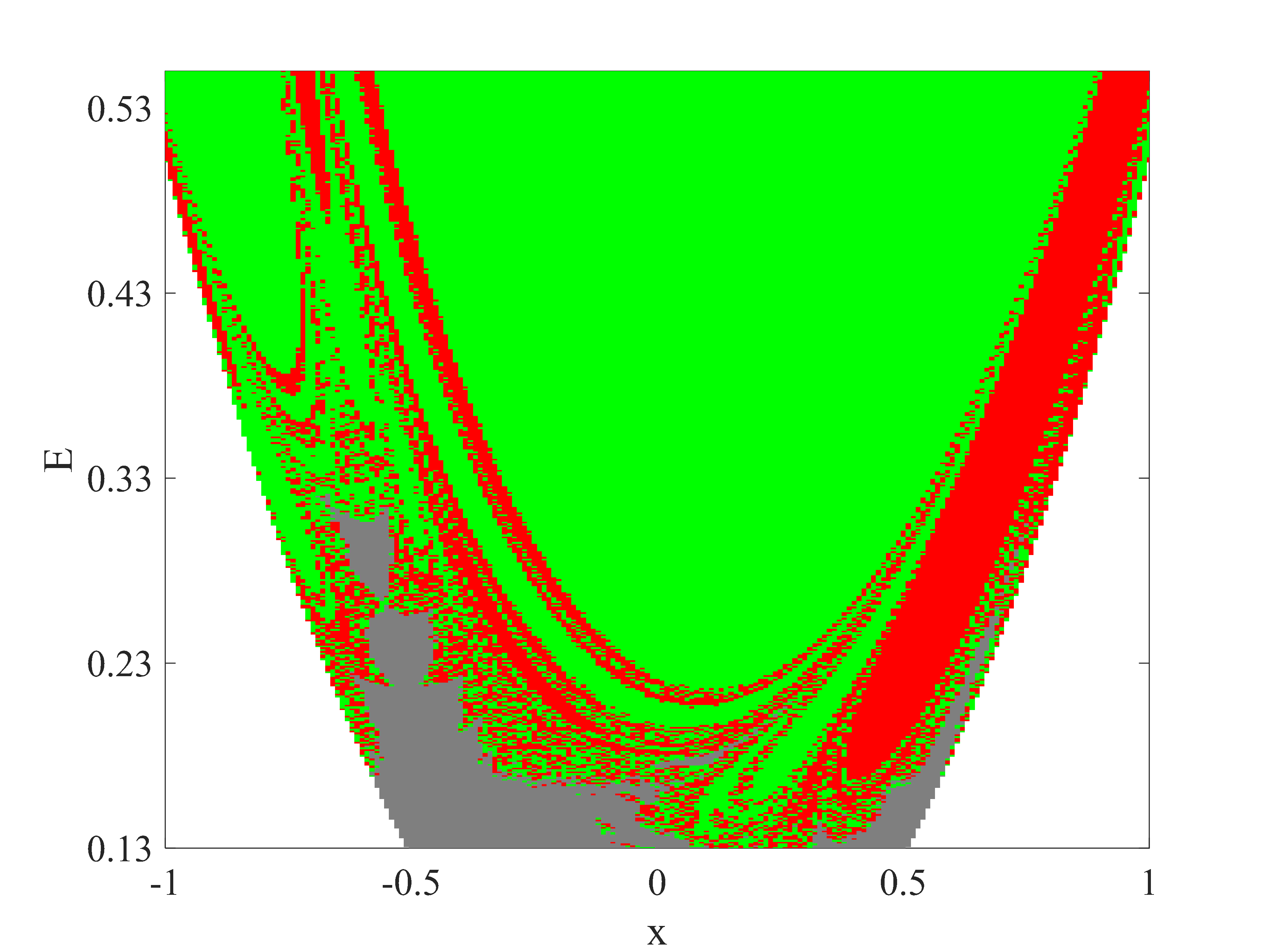}
			\caption{}\label{fig:two-exit-energy-basins:a}
		\end{subfigure}\hfill
		\begin{subfigure}[t]{.47\linewidth}
			\centering
			\includegraphics[width=\linewidth]{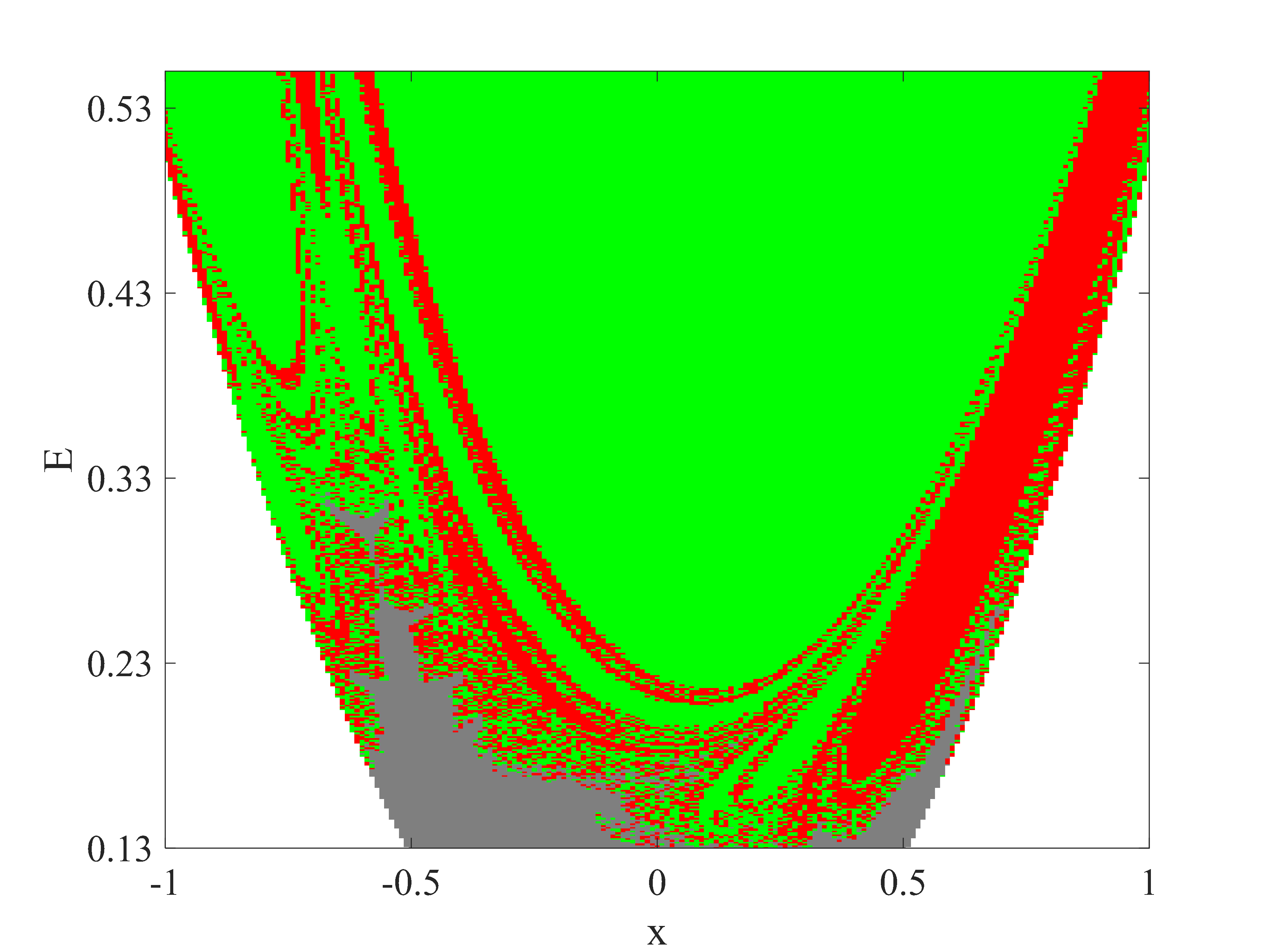}
			\caption{}\label{fig:two-exit-energy-basins:b}
		\end{subfigure}
		\caption{Prediction of the structure of the escape channels in the $(x,\allowbreak E)$ plane. Gray: trapped; green: escape through exit 1; red: escape through exit 2. (\subref{fig:two-exit-energy-basins:a}) True structure. (\subref{fig:two-exit-energy-basins:b}) Predicted structure.}
		\label{fig:two-exit-energy-basins}
	\end{figure}

	When both $y_{0}=0$ and $p_{x0}=0$ are fixed, the initial momentum $p_{y0}$ for a given $x_{0}$ is chosen on the nonnegative branch using the total energy $E$ through the relation $p_{y0}=\sqrt{2E-x_{0}^{2}}$, provided that $2E-x_{0}^{2}\geqslant{}0$. Using this constraint, we construct the structure of escape channels in the $(x,\allowbreak E)$ plane, where the initial coordinate $x_{0}$ varies over $[-1,\allowbreak 1]$ and the energy $E$ ranges from 0.13 to 0.55. Fig.~\ref{fig:two-exit-energy-basins}(\subref{fig:two-exit-energy-basins:a}) displays the true structure obtained from direct numerical integration, while Fig.~\ref{fig:two-exit-energy-basins}(\subref{fig:two-exit-energy-basins:b}) shows the corresponding prediction generated by the model. The prediction accuracy reaches $91.7\%$.
	
	Escape dynamics in Hamiltonian systems is a highly complex phenomenon. Nevertheless, these experimental results demonstrate that the NG-RC model can effectively predict such dynamics. Specifically, the model not only predicts the escape threshold with high accuracy but also reliably reconstructs the structure of the escape channels. The high accuracy achieved over a broad energy range highlights the model's capability to capture the global geometry of the escape-channel map, even though it was trained only on data from a limited set of energy values below the threshold.
	
	\Needspace{4\baselineskip}\subsection{A Hamiltonian system with three exit channels}\label{sec:three-exit}
\begin{figure}[!htbp]
		\centering
		\begin{subfigure}[t]{.49\linewidth}
			\centering
			\includegraphics[width=\linewidth]{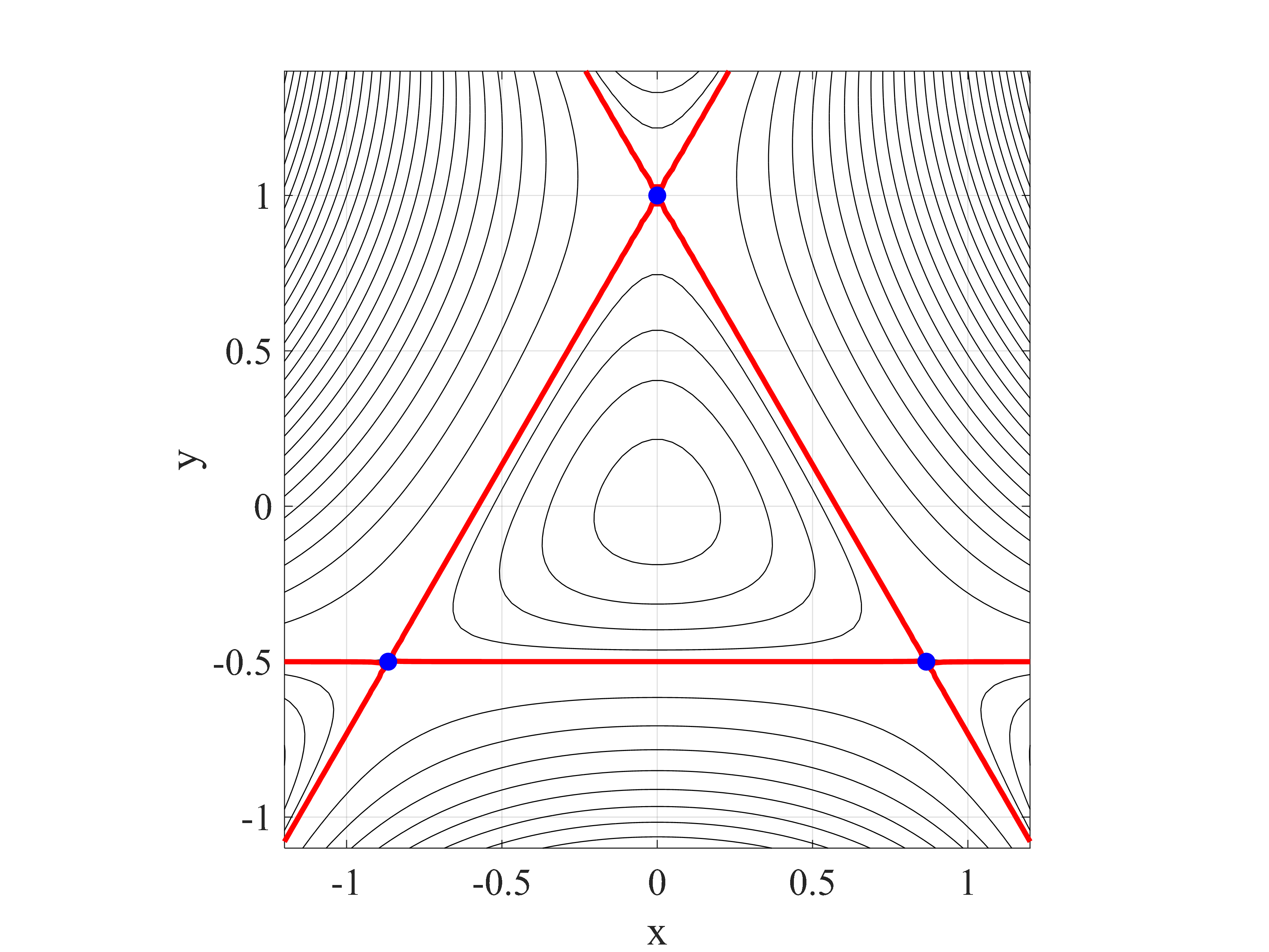}
			\caption{}\label{fig:henon-heiles-geometry:a}
		\end{subfigure}\hfill
		\begin{subfigure}[t]{.49\linewidth}
			\centering
			\includegraphics[width=\linewidth]{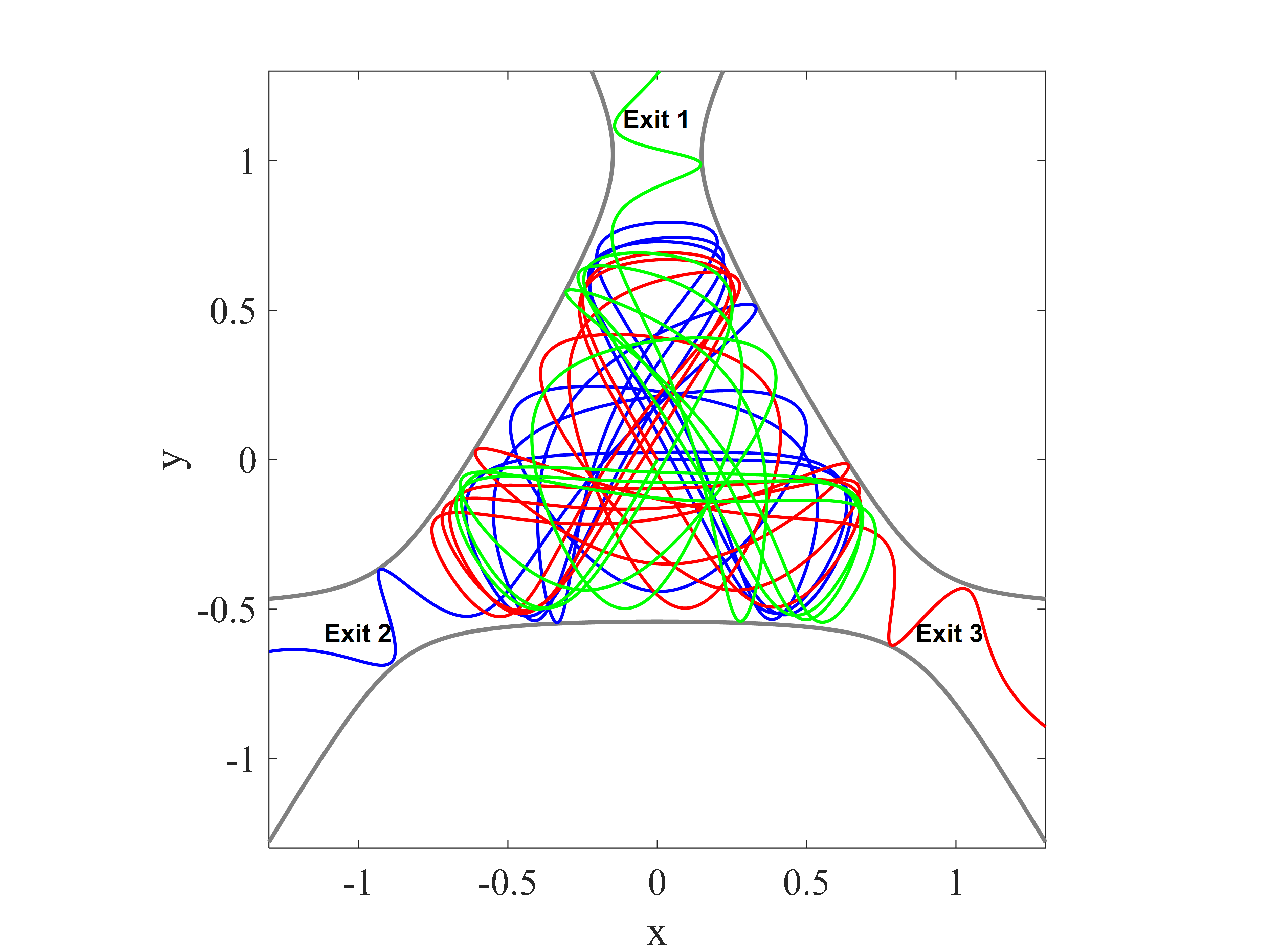}
			\caption{}\label{fig:henon-heiles-geometry:b}
		\end{subfigure}
		\caption{Equipotential curves and exit channels of the H\'enon-Heiles system. (\subref{fig:henon-heiles-geometry:a}) The equipotential curves for different energy values. (\subref{fig:henon-heiles-geometry:b}) The three exit channels and three escaping trajectories. The trajectory escaping through exit 1 is plotted in green, the one escaping through exit 2 in blue, and the one escaping through exit 3 in red.}
		\label{fig:henon-heiles-geometry}
	\end{figure}

	We consider the well-known H\'enon-Heiles system, which was originally introduced to model the motion of a star in a galactic potential. Due to its rich dynamics, including chaotic behavior and escape phenomena, it has been widely studied as a benchmark system in nonlinear dynamics \cite{aguirre2001-1,barrio2009-3,costin2025-14,fordy1991-19,seoane2007-39,seoane2008-40,zotos2015-53}. The potential function is given by
	
	\begin{equation}\notag 
		V(x,\allowbreak y)=\frac{1}{2}(x^{2}+y^{2})+x^{2}y-\frac{1}{3}y^{3}.
	\end{equation}
	Let $p_{x}=\dot{x}$ and $p_{y}=\dot{y}$. With a unit mass, the Hamiltonian takes the standard form
	
	\begin{equation}\notag
		H=\frac{1}{2}(p_{x}^{2}+p_{y}^{2})+\frac{1}{2}(x^{2}+y^{2})+x^{2}y-\frac{1}{3}y^{3}.
	\end{equation}
	From the Hamiltonian $H$, the equations of motion are derived as follows:
	
	\begin{equation}
		\begin{aligned}\dot{x}&=p_{x},\allowbreak  \\ \dot{y}&=p_{y},\allowbreak  \\ \dot{p}_{x}&=-x-2xy,\allowbreak  \\ \dot{p}_{y}&=-y-x^{2}+y^{2}.\end{aligned}
		\label{eq:henon-heiles-motion}
	\end{equation}
	The escape threshold of the H\'enon-Heiles system is $E_{esc}=1/6$. For energy $E<E_{esc}$, the central well is closed, preventing escape from it. In contrast, for energy $E>E_{esc}$, the equipotential curves open and three exit channels appear, through which trajectories may escape to infinity. Fig.~\ref{fig:henon-heiles-geometry}(\subref{fig:henon-heiles-geometry:a}) shows the equipotential curves of the H\'enon-Heiles potential. The equipotential corresponding to the escape energy is plotted in red. Furthermore, the potential has three saddle points at $(0,\allowbreak 1)$, $(-\frac{\sqrt{3}}{2},\allowbreak -\frac{1}{2})$, and $(\frac{\sqrt{3}}{2},\allowbreak -\frac{1}{2})$, which are marked with blue dots. The escape channels and representative trajectories are shown in Fig.~\ref{fig:henon-heiles-geometry}(\subref{fig:henon-heiles-geometry:b}).
	
	With $x_{0}=0$, for initial value $(y_{0},\allowbreak p_{y0})$, the potential energy is
	
	\begin{equation*}
		V(0,\allowbreak y_{0})=\frac{1}{2}y_{0}^{2}-\frac{1}{3}y_{0}^{3}.
	\end{equation*}
	The initial value $(y_{0},\allowbreak p_{y0})$ needs to satisfy the condition
	
	\begin{equation}
		V(0,\allowbreak y_{0})\leqslant{}E.
		\label{eq:henon-heiles-potential-bound}
	\end{equation}
	Define the following function:
	
	\begin{equation*}
		f(y,\allowbreak p_{y},\allowbreak E)=\sqrt{2(E-V(0,\allowbreak y))-p_{y}^{2}}.
	\end{equation*}
	For $x_{0}=0$ and initial value $(y_{0},\allowbreak p_{y0})$, we set $p_{x0}=f(y_{0},\allowbreak p_{y0},\allowbreak E)$ so that the total energy of the system is $E$. To ensure that the initial conditions are real-valued, it is required that
	
	\begin{equation}
		2(E-V(0,\allowbreak y_{0}))\geqslant{}p_{y0}^{2}.
		\label{eq:henon-heiles-admissibility}
	\end{equation}
	We select three energy values, $E=E_{1}=0.14$, $E=E_{2}=0.15$, and $E=E_{3}=0.16$, and then define three initial values:
	
	\begin{equation*}
		\boldsymbol{z}_{1}=(0,\allowbreak  -0.1,\allowbreak  f(-0.1,\allowbreak 0.1,\allowbreak E_{1}),\allowbreak  0.1),
	\end{equation*}
	\begin{equation*}
		\boldsymbol{z}_{2}=(0,\allowbreak  -0.1,\allowbreak  f(-0.1,\allowbreak 0.1,\allowbreak E_{2}),\allowbreak  0.1),
	\end{equation*}
\begin{equation*}
	\boldsymbol{z}_{3}=(0,\allowbreak  -0.1,\allowbreak  f(-0.1,\allowbreak 0.1,\allowbreak E_{3}),\allowbreak  0.1).
\end{equation*}
			We generate data from Eq. \eqref{eq:henon-heiles-motion} with an output sampling interval of approximately $\Delta{}t=0.03$ under $\boldsymbol{z}_{1}$, $\boldsymbol{z}_{2}$, and $\boldsymbol{z}_{3}$. The three time series here, each has a length of 8000, are employed to train the model. Two hyperparameters are set as follows: $k=3$, $\gamma{}=5\times{}10^{-5}$.
			
\begin{figure}[!htbp]
				\centering
				\begin{subfigure}[t]{.47\linewidth}
					\centering
					\includegraphics[width=\linewidth]{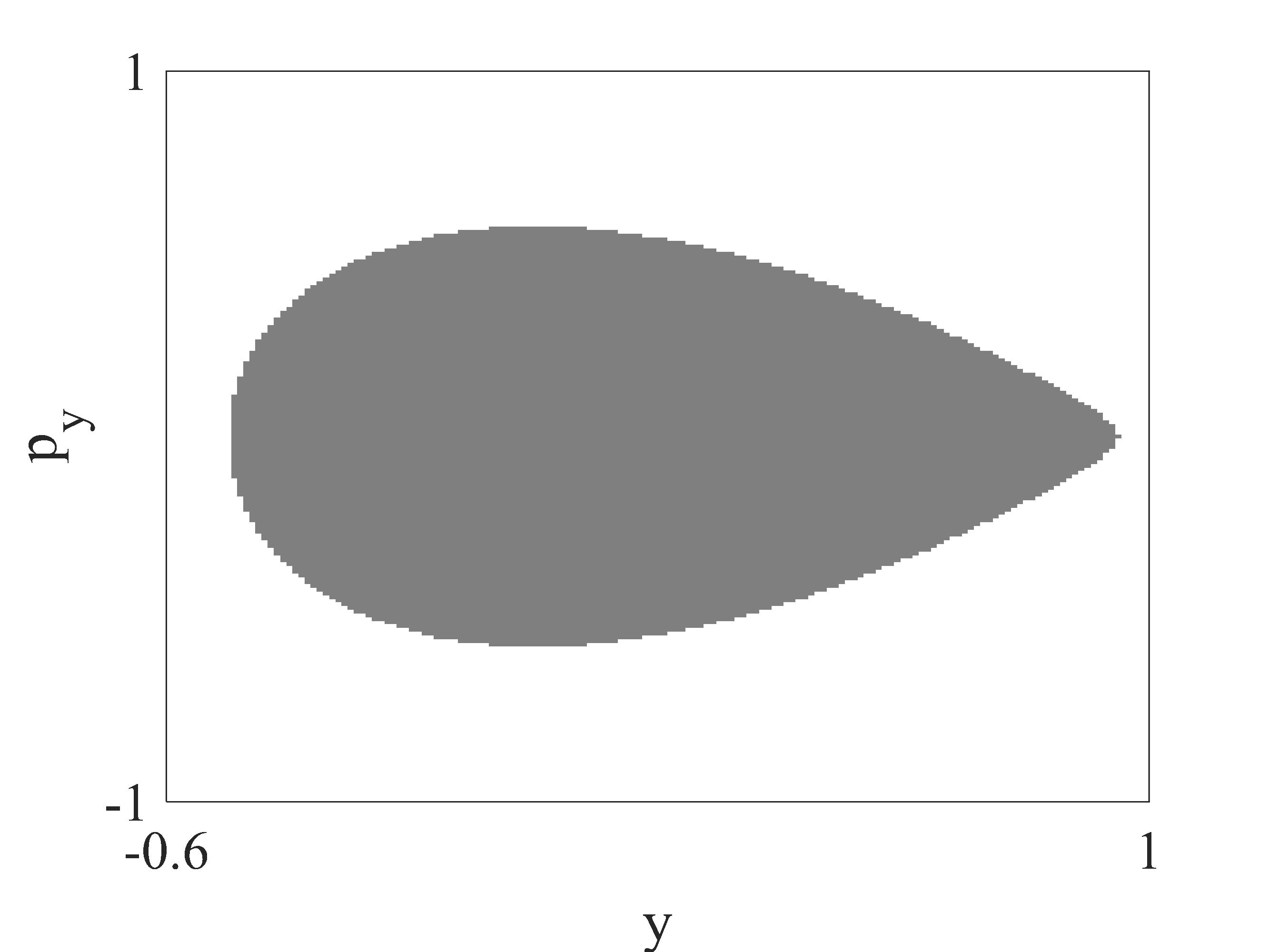}
					\caption{}\label{fig:henon-heiles-threshold-trajectories:a}
				\end{subfigure}\hfill
				\begin{subfigure}[t]{.47\linewidth}
					\centering
					\includegraphics[width=\linewidth]{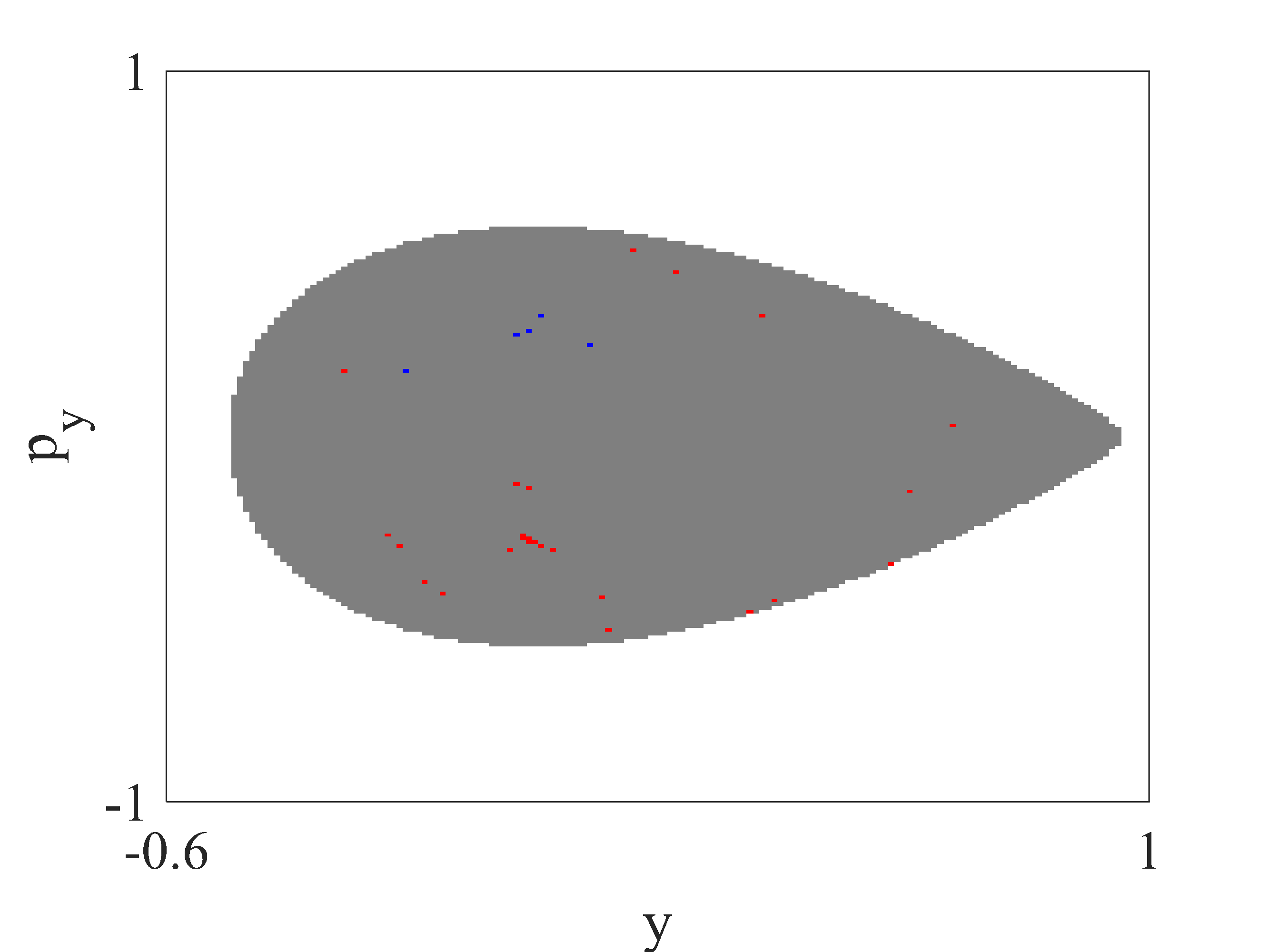}
					\caption{}\label{fig:henon-heiles-threshold-trajectories:b}
				\end{subfigure}\par\smallskip
				\begin{subfigure}[t]{.47\linewidth}
					\centering
					\includegraphics[width=\linewidth]{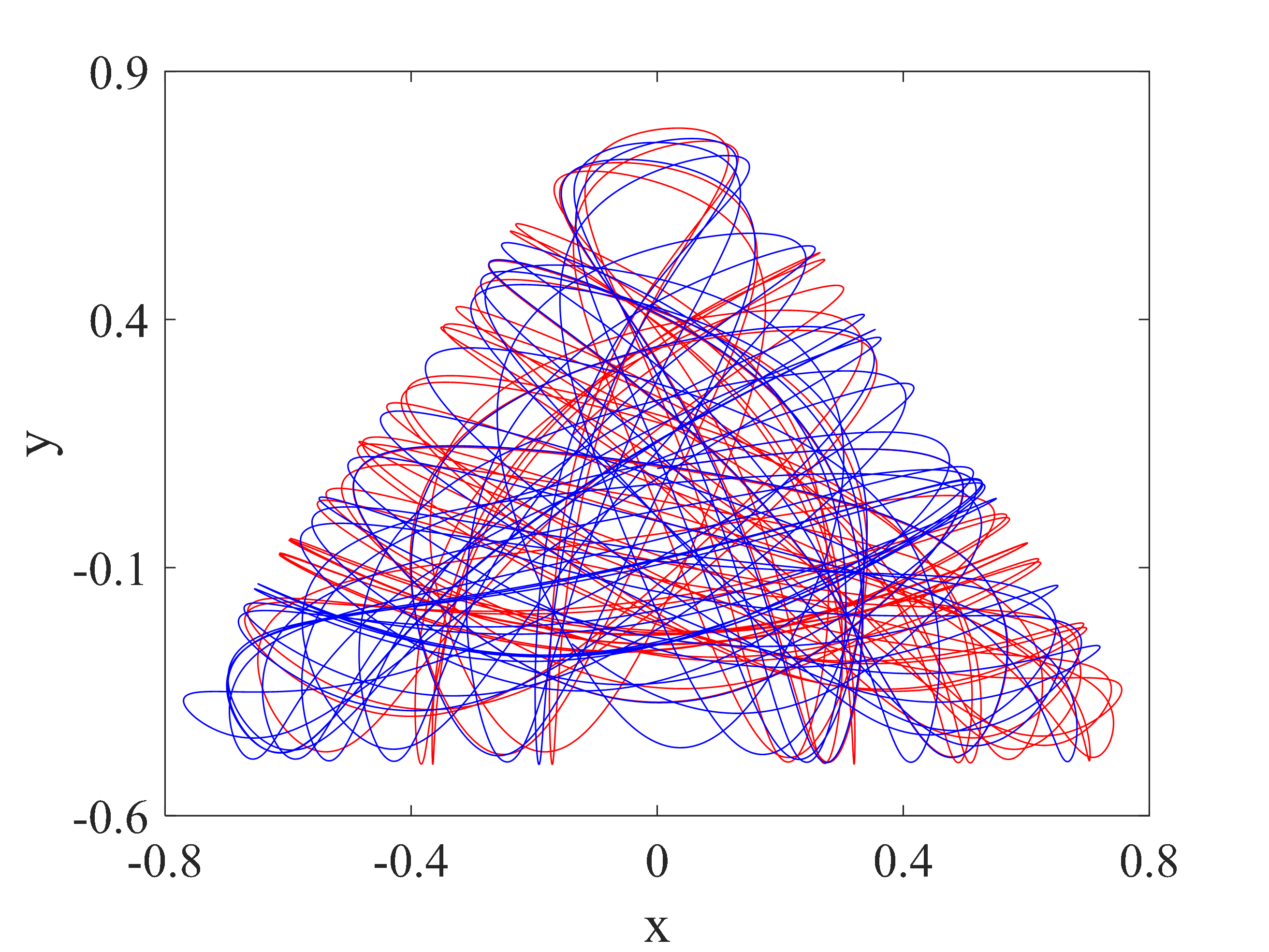}
					\caption{}\label{fig:henon-heiles-threshold-trajectories:c}
				\end{subfigure}\hfill
				\begin{subfigure}[t]{.47\linewidth}
					\centering
					\includegraphics[width=\linewidth]{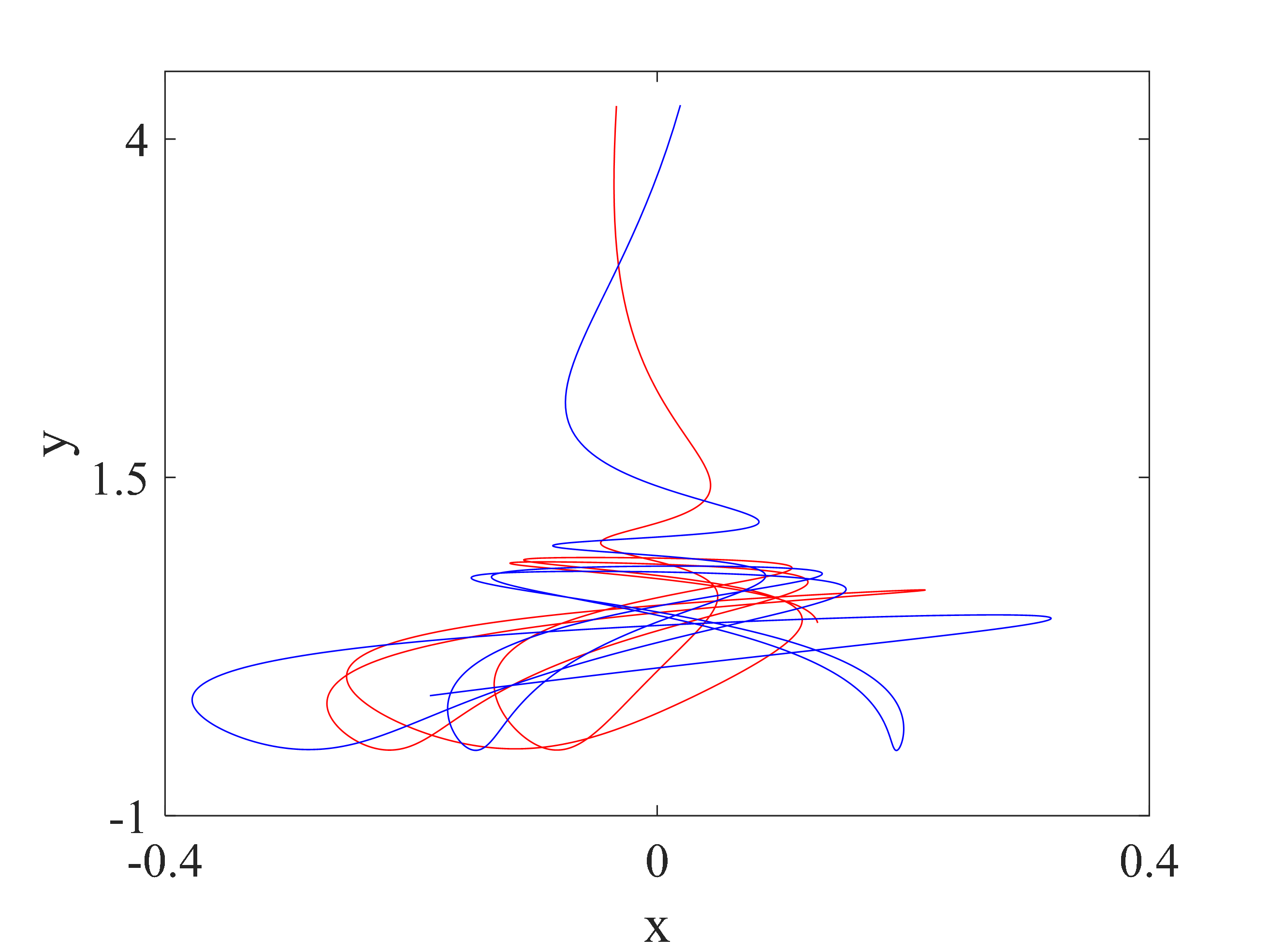}
					\caption{}\label{fig:henon-heiles-threshold-trajectories:d}
				\end{subfigure}\par\smallskip
				\begin{subfigure}[t]{.47\linewidth}
					\centering
					\includegraphics[width=\linewidth]{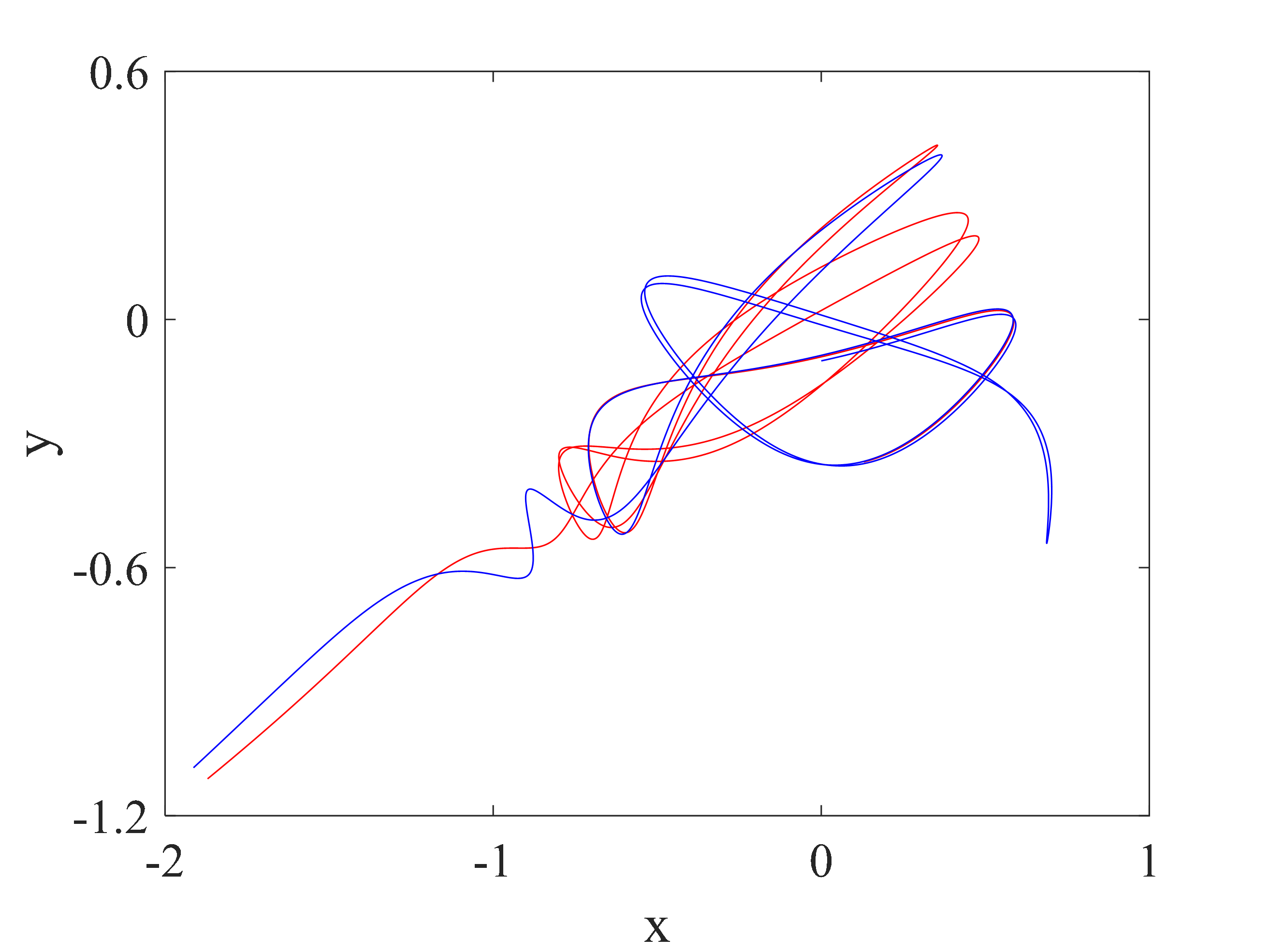}
					\caption{}\label{fig:henon-heiles-threshold-trajectories:e}
				\end{subfigure}\hfill
				\begin{subfigure}[t]{.47\linewidth}
					\centering
					\includegraphics[width=\linewidth]{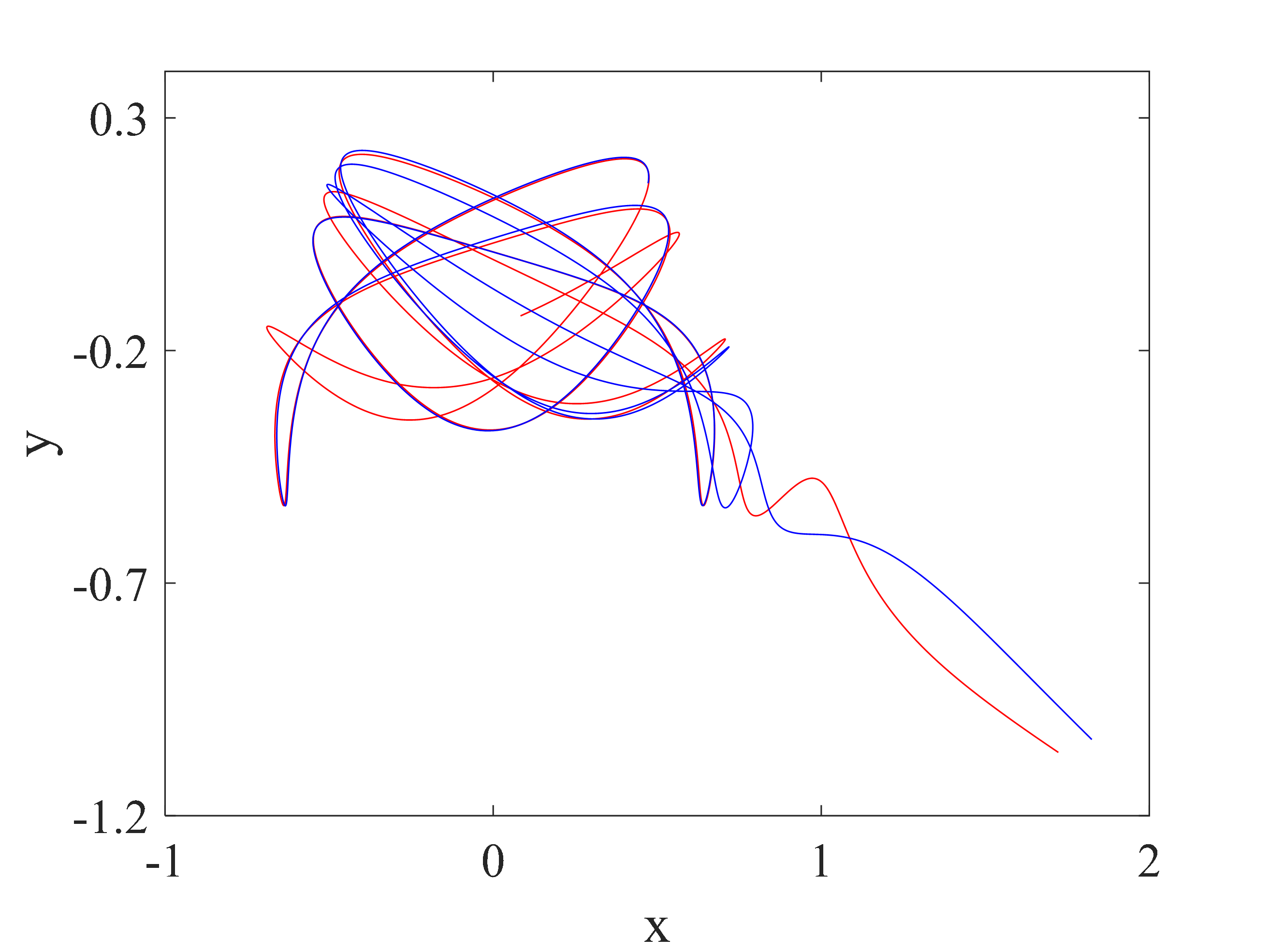}
					\caption{}\label{fig:henon-heiles-threshold-trajectories:f}
				\end{subfigure}
				\caption{Predicted escape and non-escape structures ((\subref{fig:henon-heiles-threshold-trajectories:a})--(\subref{fig:henon-heiles-threshold-trajectories:b}), gray: trapped; green: escape through exit 1; blue: escape through exit 2; red: escape through exit 3) and examples of predicted trajectories under different energy values and initial conditions ((\subref{fig:henon-heiles-threshold-trajectories:c})--(\subref{fig:henon-heiles-threshold-trajectories:f}), blue: true; red: predicted). (\subref{fig:henon-heiles-threshold-trajectories:a}) Predicted escape and non-escape structure at $E=0.1655$. (\subref{fig:henon-heiles-threshold-trajectories:b}) Predicted escape and non-escape structure at $E=0.16555$. (\subref{fig:henon-heiles-threshold-trajectories:c}) $E=E_{4}=0.165$, $\boldsymbol{z}_{4}=(0,\allowbreak  0.1,\allowbreak  f(0.1,\allowbreak 0.1,\allowbreak E_{4}),\allowbreak  0.1)$, trapped. (\subref{fig:henon-heiles-threshold-trajectories:d}) $E= E_{5}=0.18$, $\boldsymbol{z}_{5}=(0,\allowbreak  0.25,\allowbreak  f(0.25,\allowbreak 0.2,\allowbreak E_{5}),\allowbreak 0.2)$, escape through exit 1. (\subref{fig:henon-heiles-threshold-trajectories:e}) $E=E_{5}=0.18$, $\boldsymbol{z}_{6}= (0,\allowbreak -0.1,\allowbreak  f(-0.1,\allowbreak 0.1,\allowbreak E_{5}),\allowbreak  0.1)$, escape through exit 2. (\subref{fig:henon-heiles-threshold-trajectories:f}) $E=E_{5}=0.18$, $\boldsymbol{z}_{7}=(0,\allowbreak  -0.15,\allowbreak  f(-0.15,\allowbreak 0.15,\allowbreak E_{5}),\allowbreak  0.15)$, escape through exit 3.}
				\label{fig:henon-heiles-threshold-trajectories}
			\end{figure}

			We predict the system's evolution at energy values ranging from 0.161 to 0.178 in increments of 0.001. For each energy value, $(y_{0},\allowbreak p_{y0})$ is sampled over a uniform grid covering the domain $[-0.6,\allowbreak 1]\times{}[-1,\allowbreak 1]$ with a step size of 0.01. $x_{0}$ is fixed at 0, while $p_{x0}$ is determined by the function $f(y,\allowbreak p_{y},\allowbreak E)$. Only the points that satisfy conditions \eqref{eq:henon-heiles-potential-bound} and \eqref{eq:henon-heiles-admissibility} are considered. Table~\ref{tab:threshold-3} presents the prediction results. For $E\leqslant{}0.165$, no predicted escape is detected, while for every tested energy $E\geqslant{}0.166$, at least one predicted trajectory is classified as escaping, so the estimated escape threshold is between 0.165 and 0.166. Furthermore, we conduct predictions at energy values ranging from 0.16505 to 0.166, and the results are shown in Table~\ref{tab:threshold-4}. In addition, Fig.~\ref{fig:henon-heiles-threshold-trajectories}(\subref{fig:henon-heiles-threshold-trajectories:a})--(\subref{fig:henon-heiles-threshold-trajectories:b}) display the predicted escape and non-escape structures in the considered domain at $E=0.1655$ and $E=0.16555$, respectively. Therefore, the grid-based threshold estimate is 0.16555, which is very close to the theoretical value.
			
			\begin{table}[!htbp]
				\centering
				\caption{Predictions of escape and non-escape for energy values from 0.161 to 0.178.}
				\label{tab:threshold-3}
				\small
				\setlength{\tabcolsep}{0pt}
				\renewcommand{\arraystretch}{1.3}
				\begin{tabular}{|>{\centering\arraybackslash}m{16mm}|*{9}{>{\centering\arraybackslash}m{15.6mm}|}}
					\hline
					$E$ & 0.161 & 0.162 & 0.163 & 0.164 & 0.165 & 0.166 & 0.167 & 0.168 & 0.169 \\
					\hline
					Result & $\times$ & $\times$ & $\times$ & $\times$ & $\times$ & $\checkmark$ & $\checkmark$ & $\checkmark$ & $\checkmark$ \\
					\hline
					$E$ & 0.17 & 0.171 & 0.172 & 0.173 & 0.174 & 0.175 & 0.176 & 0.177 & 0.178 \\
					\hline
					Result & $\checkmark$ & $\checkmark$ & $\checkmark$ & $\checkmark$ & $\checkmark$ & $\checkmark$ & $\checkmark$ & $\checkmark$ & $\checkmark$ \\
					\hline
				\end{tabular}
			\end{table}
			
			\begin{table}[!htbp]
				\centering
				\caption{Predictions of escape and non-escape for energy values from 0.16505 to 0.166.}
				\label{tab:threshold-4}
				\small
				\setlength{\tabcolsep}{0pt}
				\renewcommand{\arraystretch}{1.3}
				\begin{tabular}{|>{\centering\arraybackslash}m{16mm}|*{8}{>{\centering\arraybackslash}m{17.55mm}|}}
					\hline
					$E$ & 0.16505 & 0.1651 & 0.16515 & 0.1652 & 0.16525 & 0.1653 & 0.16535 & 0.1654 \\
					\hline
					Result & $\times$ & $\times$ & $\times$ & $\times$ & $\times$ & $\times$ & $\times$ & $\times$ \\
					\hline
					$E$ & 0.16545 & 0.1655 & 0.16555 & 0.1656 & 0.16565 & 0.1657 & 0.16575 & 0.1658 \\
					\hline
					Result & $\times$ & $\times$ & $\checkmark$ & $\checkmark$ & $\checkmark$ & $\checkmark$ & $\checkmark$ & $\checkmark$ \\
					\hline
					$E$ & 0.16585 & 0.1659 & 0.16595 & 0.166 &  &  &  &  \\
					\hline
					Result & $\checkmark$ & $\checkmark$ & $\checkmark$ & $\checkmark$ &  &  &  &  \\
					\hline
				\end{tabular}
			\end{table}
			
			Denote $E=E_{4}=0.165$, and $E=E_{5}=0.18$. Define $\boldsymbol{z}_{4}=(0,\allowbreak  0.1,\allowbreak  f(0.1,\allowbreak 0.1,\allowbreak E_{4}),\allowbreak  0.1)$, $\boldsymbol{z}_{5}=(0,\allowbreak  0.25,\allowbreak  f(0.25,\allowbreak 0.2,\allowbreak E_{5}),\allowbreak 0.2)$, $\boldsymbol{z}_{6}=(0,\allowbreak  -0.1,\allowbreak  f(-0.1,\allowbreak 0.1,\allowbreak E_{5}),\allowbreak  0.1)$, and $\boldsymbol{z}_{7}=(0,\allowbreak  -0.15,\allowbreak   f(-0.15,\allowbreak 0.15,\allowbreak E_{5}),\allowbreak  0.15)$. Fig.~\ref{fig:henon-heiles-threshold-trajectories}(\subref{fig:henon-heiles-threshold-trajectories:c})--(\subref{fig:henon-heiles-threshold-trajectories:f}) present the predicted and true trajectories under $\boldsymbol{z}_{4}$, $\boldsymbol{z}_{5}$, $\boldsymbol{z}_{6}$, and $\boldsymbol{z}_{7}$. For initial value $\boldsymbol{z}_{4}$, the system is bounded, and the predicted orbit projection reproduces the shape of the reference orbit projection. For initial values $\boldsymbol{z}_{5}$, $\boldsymbol{z}_{6}$ and $\boldsymbol{z}_{7}$, the system escapes to infinity, specifically through exit 1 for $\boldsymbol{z}_{5}$, through exit 2 for $\boldsymbol{z}_{6}$ and through exit 3 for $\boldsymbol{z}_{7}$. The model correctly predicts the escape channel.
			
\begin{figure}[!htbp]
				\centering
				\includegraphics[width=.64\linewidth]{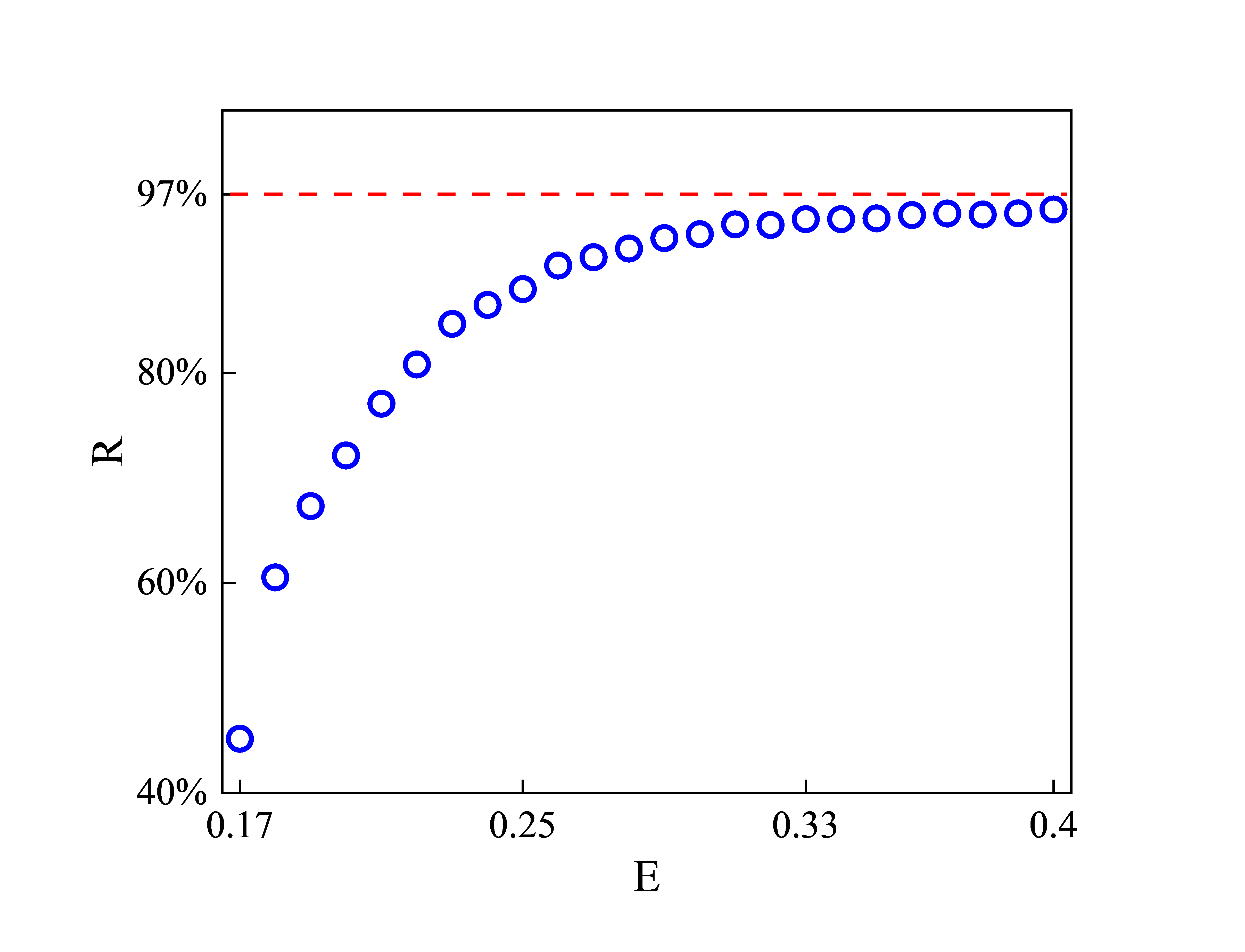}
				\caption{Prediction accuracy of the structure of the escape channels under different energy values.}
				\label{fig:henon-heiles-accuracy}
			\end{figure}

			We predict the structure of the escape channels in the $(y,\allowbreak p_{y})$ plane. For $0.17\leqslant{}E\leqslant{}0.22$, $(y_{0},\allowbreak p_{y0})$ is sampled from the domain $[-0.6,\allowbreak 1.5]\times{}[-1,\allowbreak 1]$; for $0.23\leqslant{}E\leqslant{}0.3$, the sampling domain is $[-0.7,\allowbreak 1.6]\times{}[-1.2,\allowbreak 1.2]$; and for $0.31\leqslant{}E\leqslant{}0.4$, the domain is $[-0.8,\allowbreak 1.8]\times{}[-1.4,\allowbreak 1.4]$. In all cases, the region is discretized into a uniform grid with a step size of 0.01, and every admissible grid point is used as an initial value. $x_{0}$ is fixed at 0, and $p_{x0}$ is determined by the function $f(y,\allowbreak p_{y},\allowbreak E)$. Fig.~\ref{fig:henon-heiles-basins} presents a comparison between the predicted and true structures for energy values $E=0.17,\allowbreak 0.18,\allowbreak 0.19,\allowbreak 0.20,\allowbreak 0.22,\allowbreak 0.26,\allowbreak 0.30,\allowbreak 0.35$, and $0.4$. Fig.~\ref{fig:henon-heiles-accuracy} displays the prediction accuracy. As energy increases, the accuracy exhibits a general upward trend. At $E=0.17$, the accuracy is $45.2\%$, and at $E=0.4$, the accuracy reaches $95.6\%$.
			
			\begin{figure}[!htbp]
				\centering
				\begin{subfigure}[t]{.31\linewidth}
					\centering
					\includegraphics[width=\linewidth]{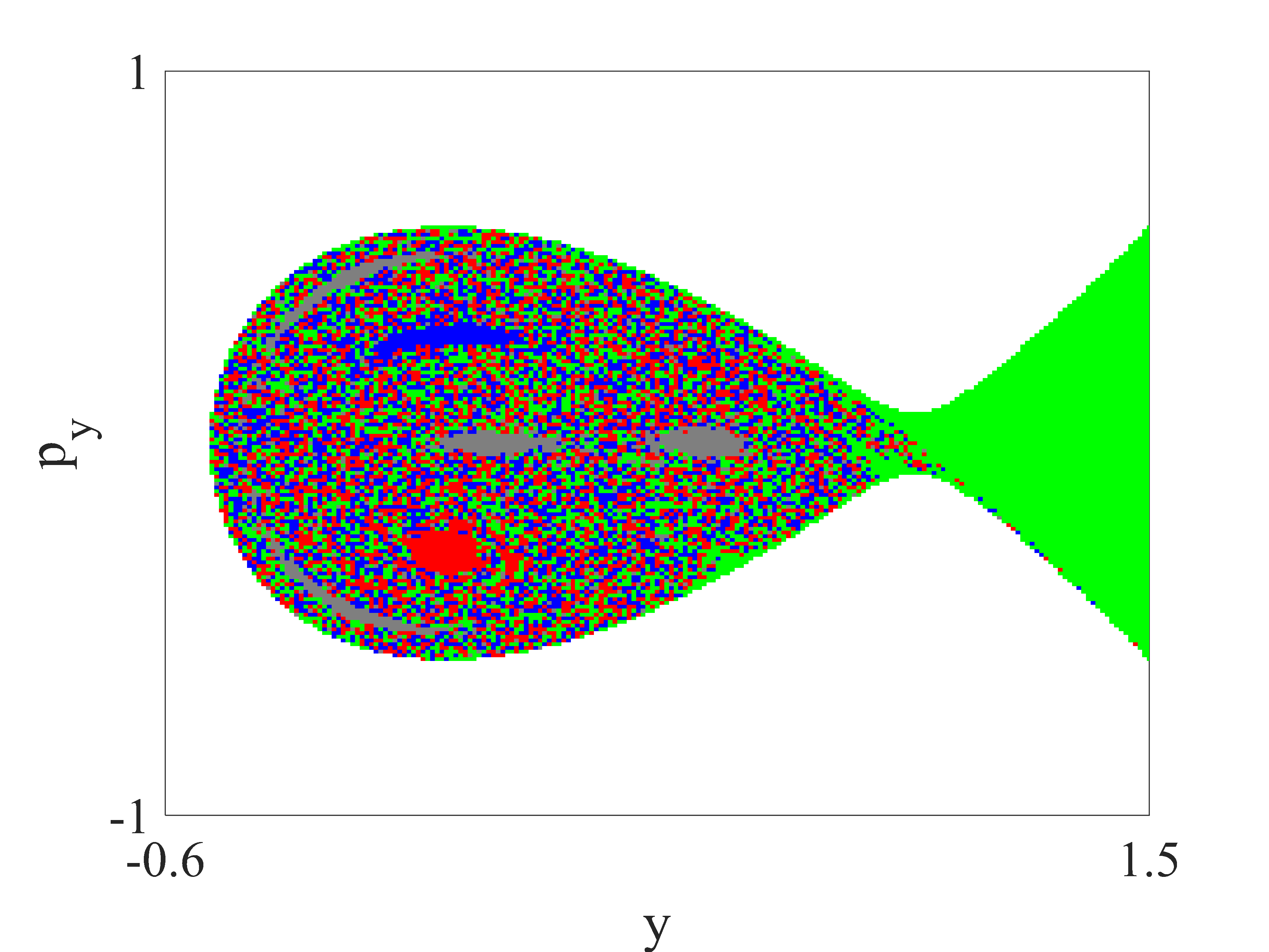}
					\caption{}\label{fig:henon-heiles-basins:a}
				\end{subfigure}\hfill
				\begin{subfigure}[t]{.31\linewidth}
					\centering
					\includegraphics[width=\linewidth]{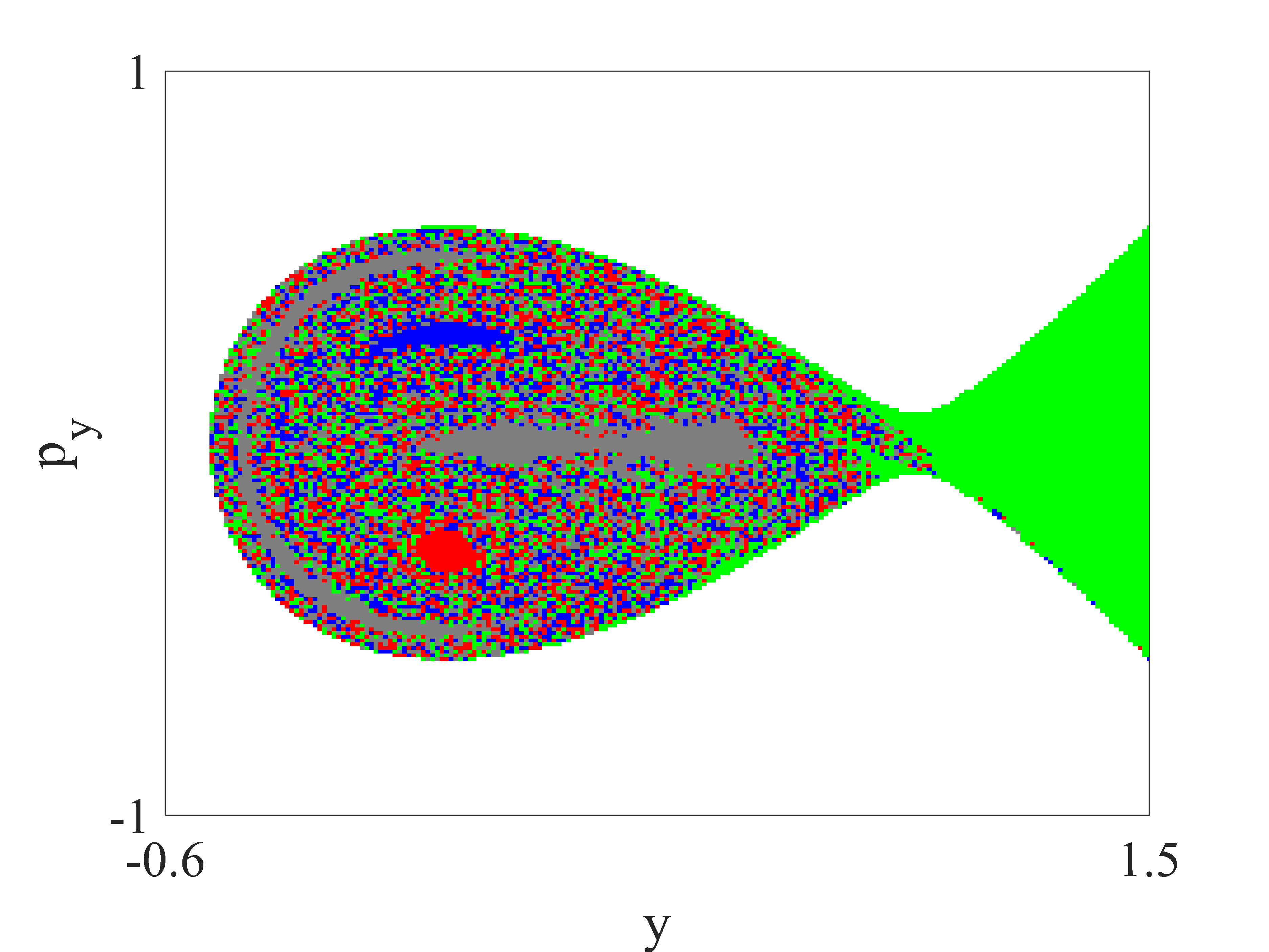}
					\caption{}\label{fig:henon-heiles-basins:b}
				\end{subfigure}\hfill
				\begin{subfigure}[t]{.31\linewidth}
					\centering
					\includegraphics[width=\linewidth]{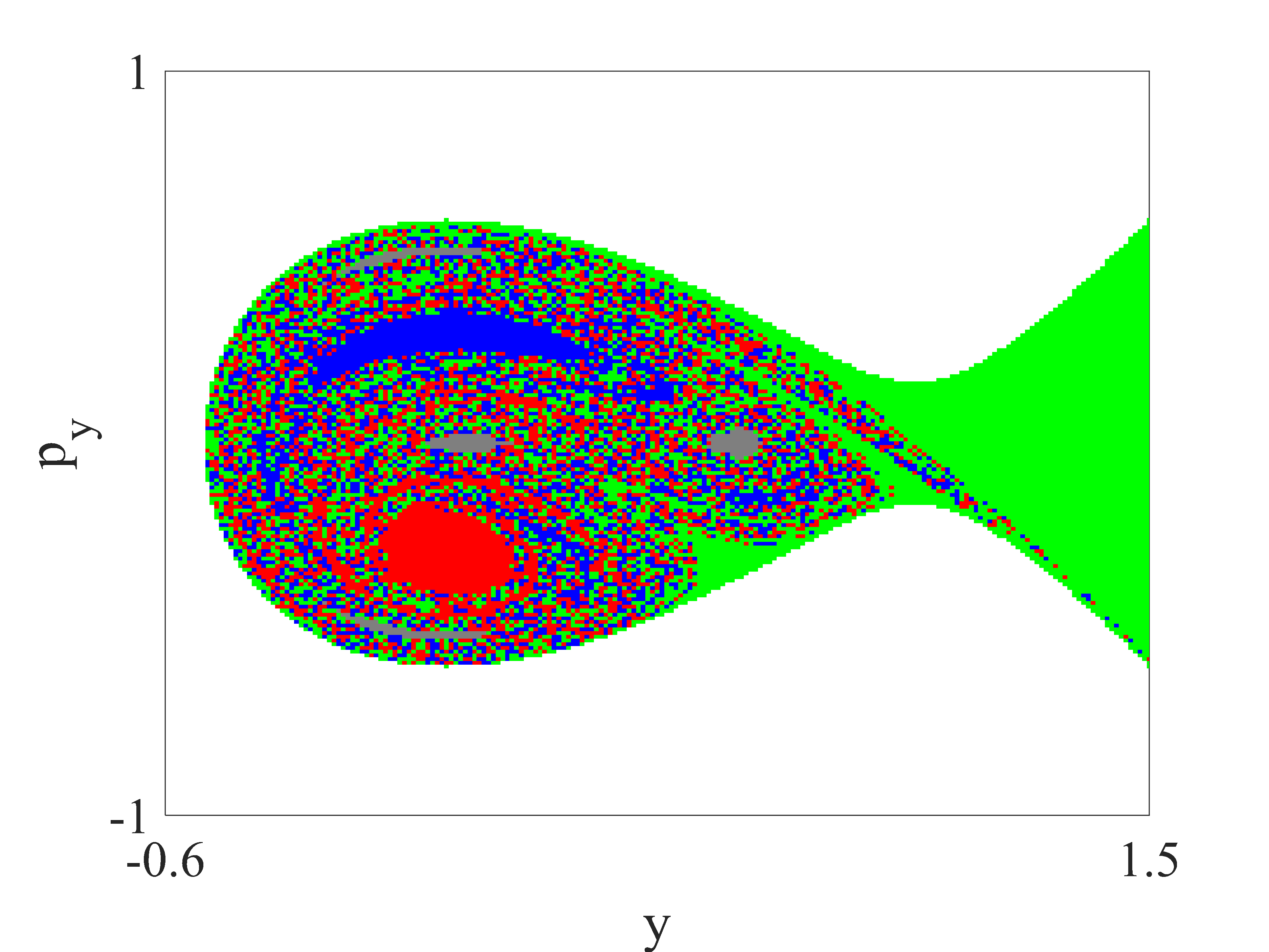}
					\caption{}\label{fig:henon-heiles-basins:c}
				\end{subfigure}\par\smallskip
				\begin{subfigure}[t]{.31\linewidth}
					\centering
					\includegraphics[width=\linewidth]{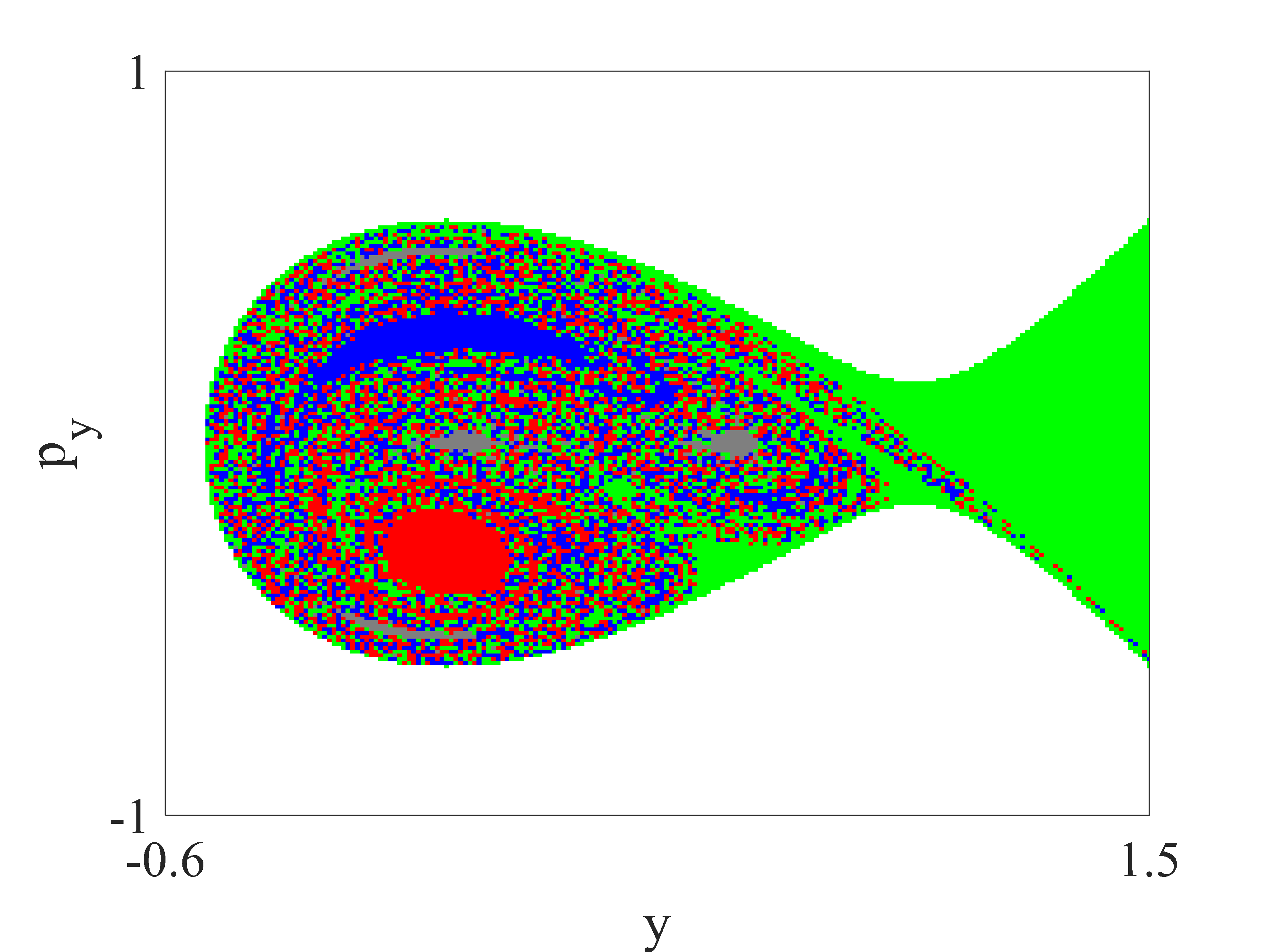}
					\caption{}\label{fig:henon-heiles-basins:d}
				\end{subfigure}\hfill
				\begin{subfigure}[t]{.31\linewidth}
					\centering
					\includegraphics[width=\linewidth]{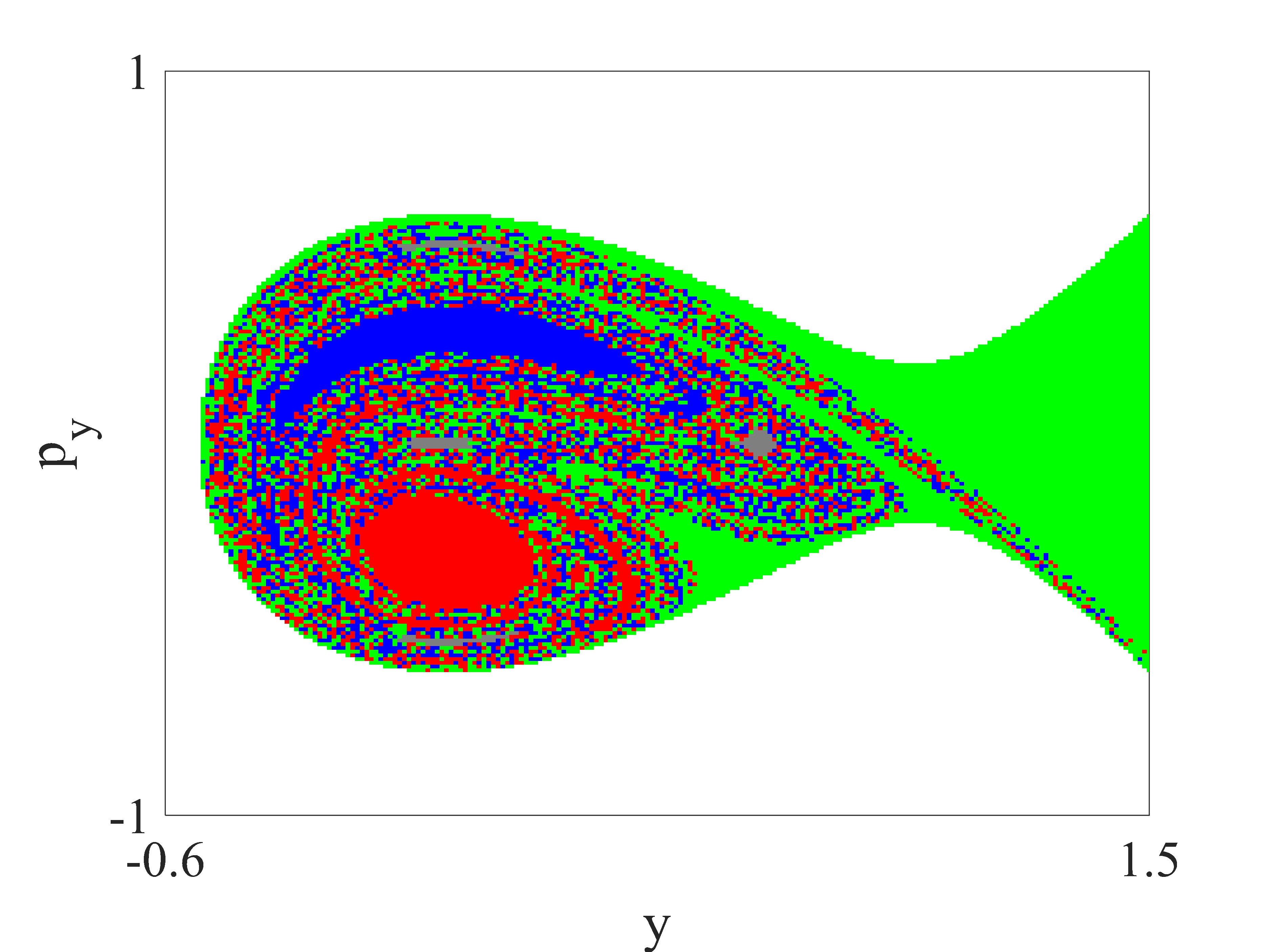}
					\caption{}\label{fig:henon-heiles-basins:e}
				\end{subfigure}\hfill
				\begin{subfigure}[t]{.31\linewidth}
					\centering
					\includegraphics[width=\linewidth]{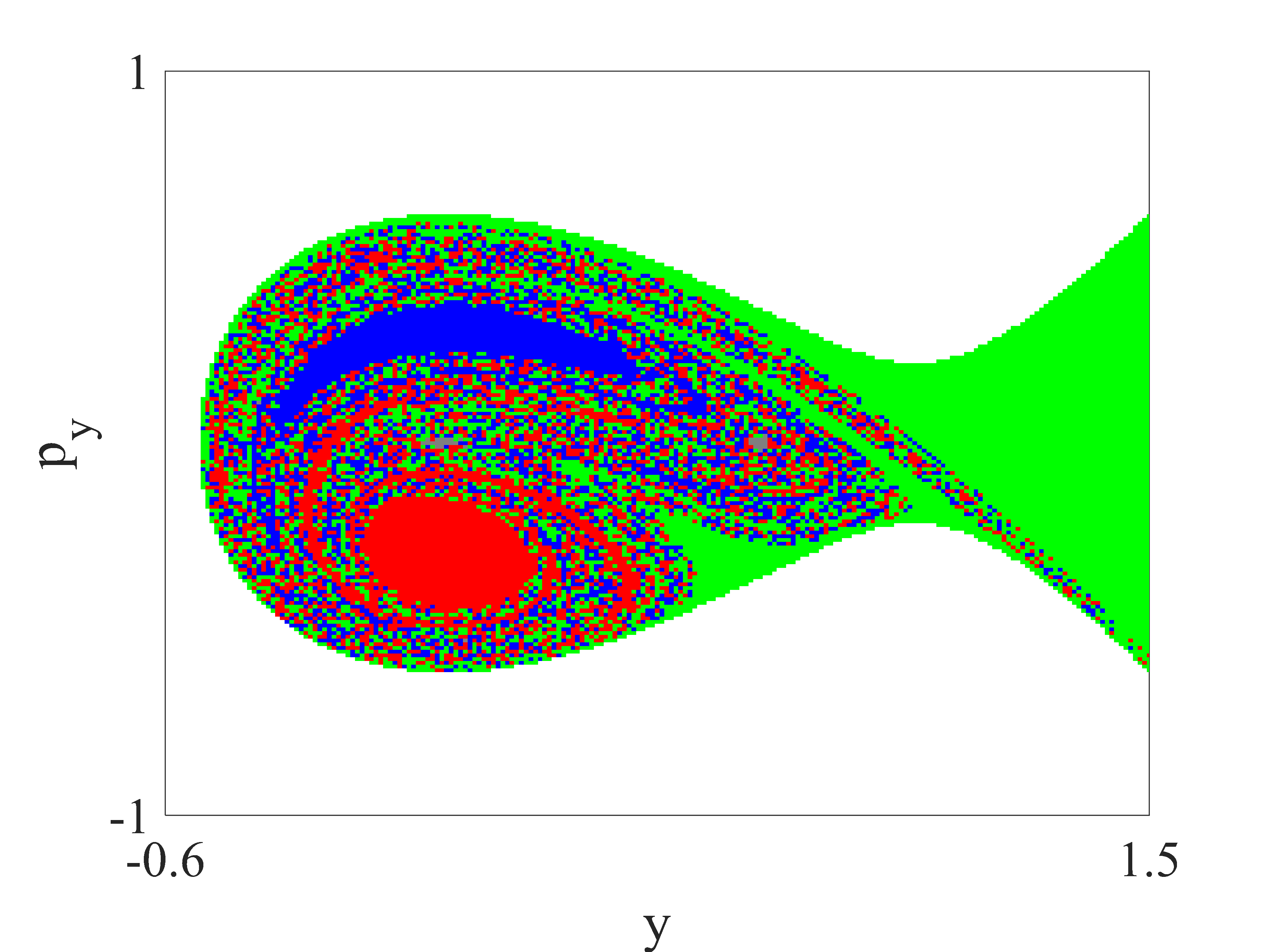}
					\caption{}\label{fig:henon-heiles-basins:f}
				\end{subfigure}
				\caption{Predictions of the structure of escape channels in the $(y,\allowbreak p_{y})$ plane for different energy values. Gray: trapped; green: escape through exit 1; blue: escape through exit 2; red: escape through exit 3. (\subref{fig:henon-heiles-basins:a}) $E=0.17$, true. (\subref{fig:henon-heiles-basins:b}) $E=0.17$, predicted. (\subref{fig:henon-heiles-basins:c}) $E=0.18$, true. (\subref{fig:henon-heiles-basins:d}) $E=0.18$, predicted. (\subref{fig:henon-heiles-basins:e}) $E=0.19$, true. (\subref{fig:henon-heiles-basins:f}) $E=0.19$, predicted. (\subref{fig:henon-heiles-basins:g}) $E=0.2$, true. (\subref{fig:henon-heiles-basins:h}) $E=0.2$, predicted. (\subref{fig:henon-heiles-basins:i}) $E=0.22$, true. (\subref{fig:henon-heiles-basins:j}) $E=0.22$, predicted. (\subref{fig:henon-heiles-basins:k}) $E=0.26$, true. (\subref{fig:henon-heiles-basins:l}) $E=0.26$, predicted. (\subref{fig:henon-heiles-basins:m}) $E=0.3$, true. (\subref{fig:henon-heiles-basins:n}) $E=0.3$, predicted. (\subref{fig:henon-heiles-basins:o}) $E=0.35$, true. (\subref{fig:henon-heiles-basins:p}) $E=0.35$, predicted. (\subref{fig:henon-heiles-basins:q}) $E=0.4$, true. (\subref{fig:henon-heiles-basins:r}) $E=0.4$, predicted.}
				\label{fig:henon-heiles-basins}
			\end{figure}
			
			\begin{figure}[!htbp]
				\centering
				\ContinuedFloat
				\begin{subfigure}[t]{.31\linewidth}
					\centering
					\includegraphics[width=\linewidth]{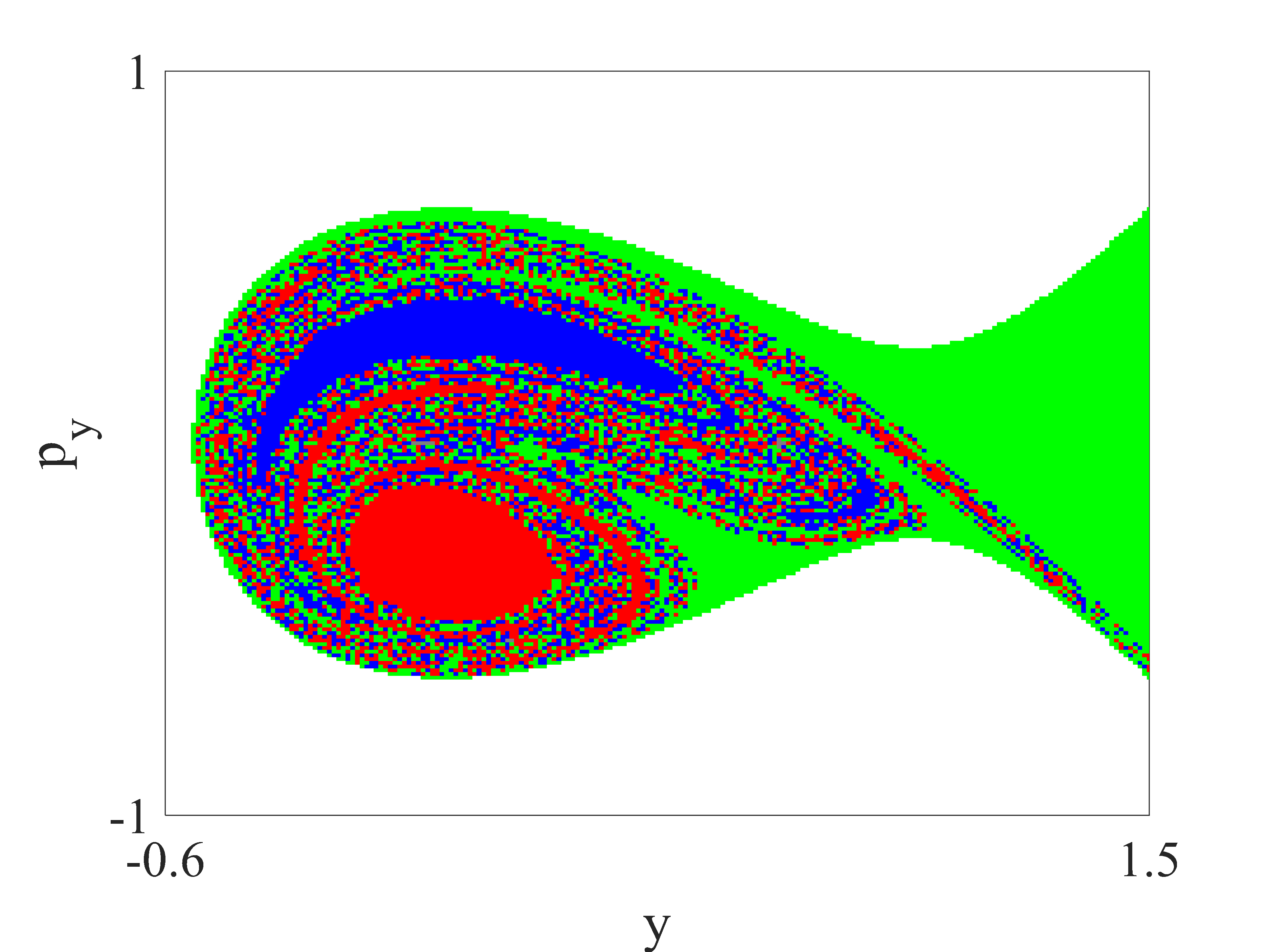}
					\caption{}\label{fig:henon-heiles-basins:g}
				\end{subfigure}\hfill
				\begin{subfigure}[t]{.31\linewidth}
					\centering
					\includegraphics[width=\linewidth]{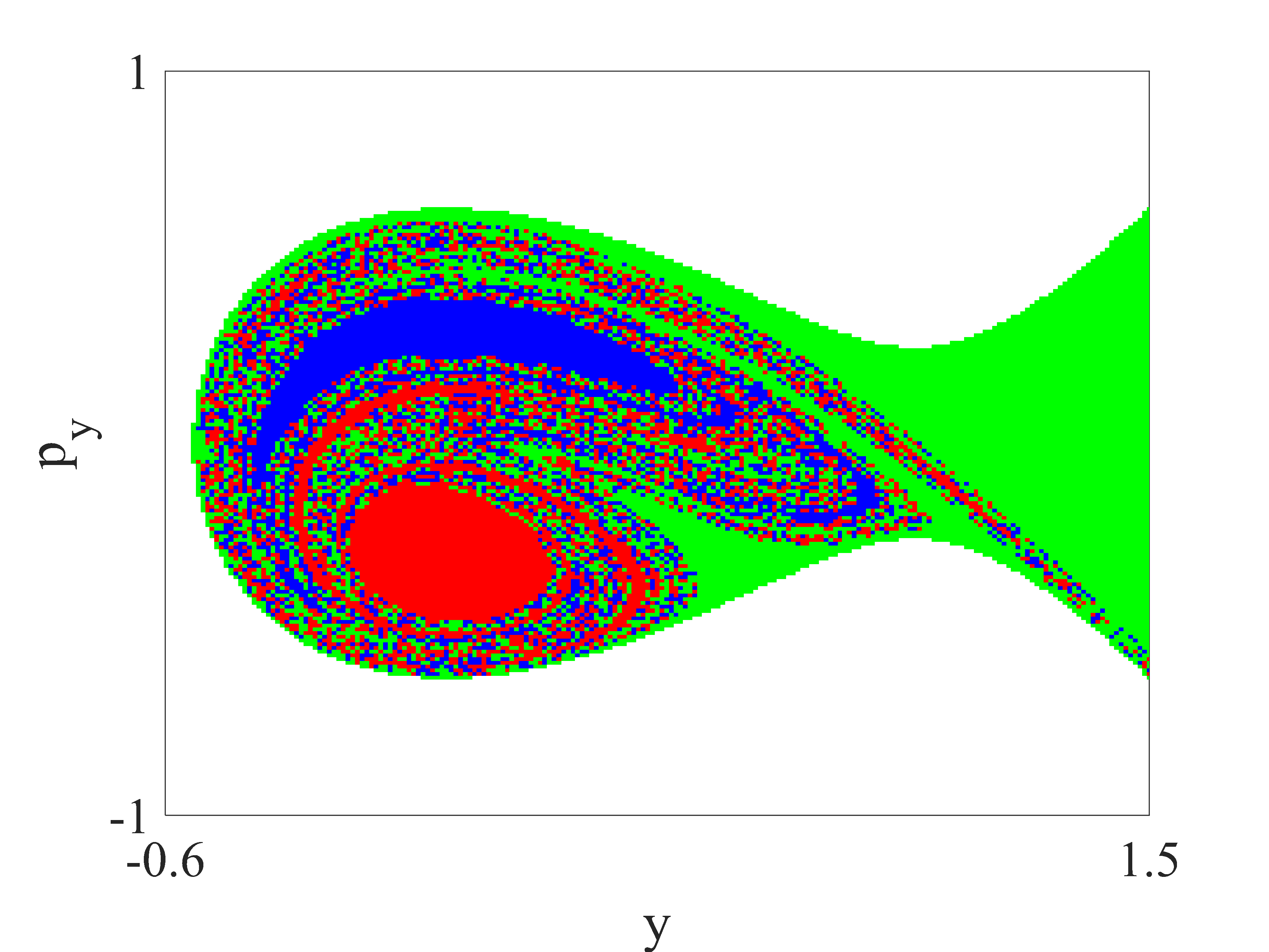}
					\caption{}\label{fig:henon-heiles-basins:h}
				\end{subfigure}\hfill
				\begin{subfigure}[t]{.31\linewidth}
					\centering
					\includegraphics[width=\linewidth]{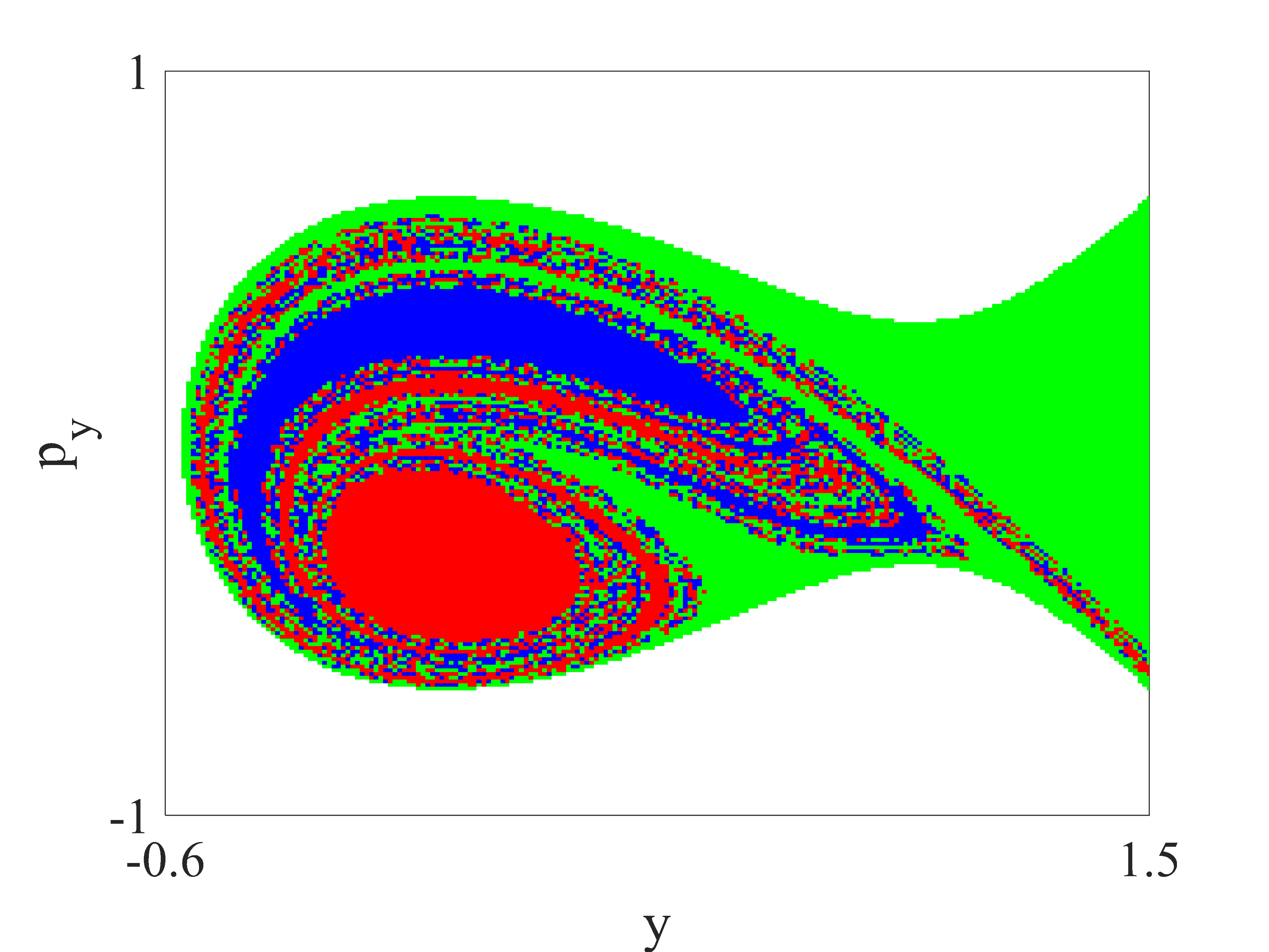}
					\caption{}\label{fig:henon-heiles-basins:i}
				\end{subfigure}\par\smallskip
				\begin{subfigure}[t]{.31\linewidth}
					\centering
					\includegraphics[width=\linewidth]{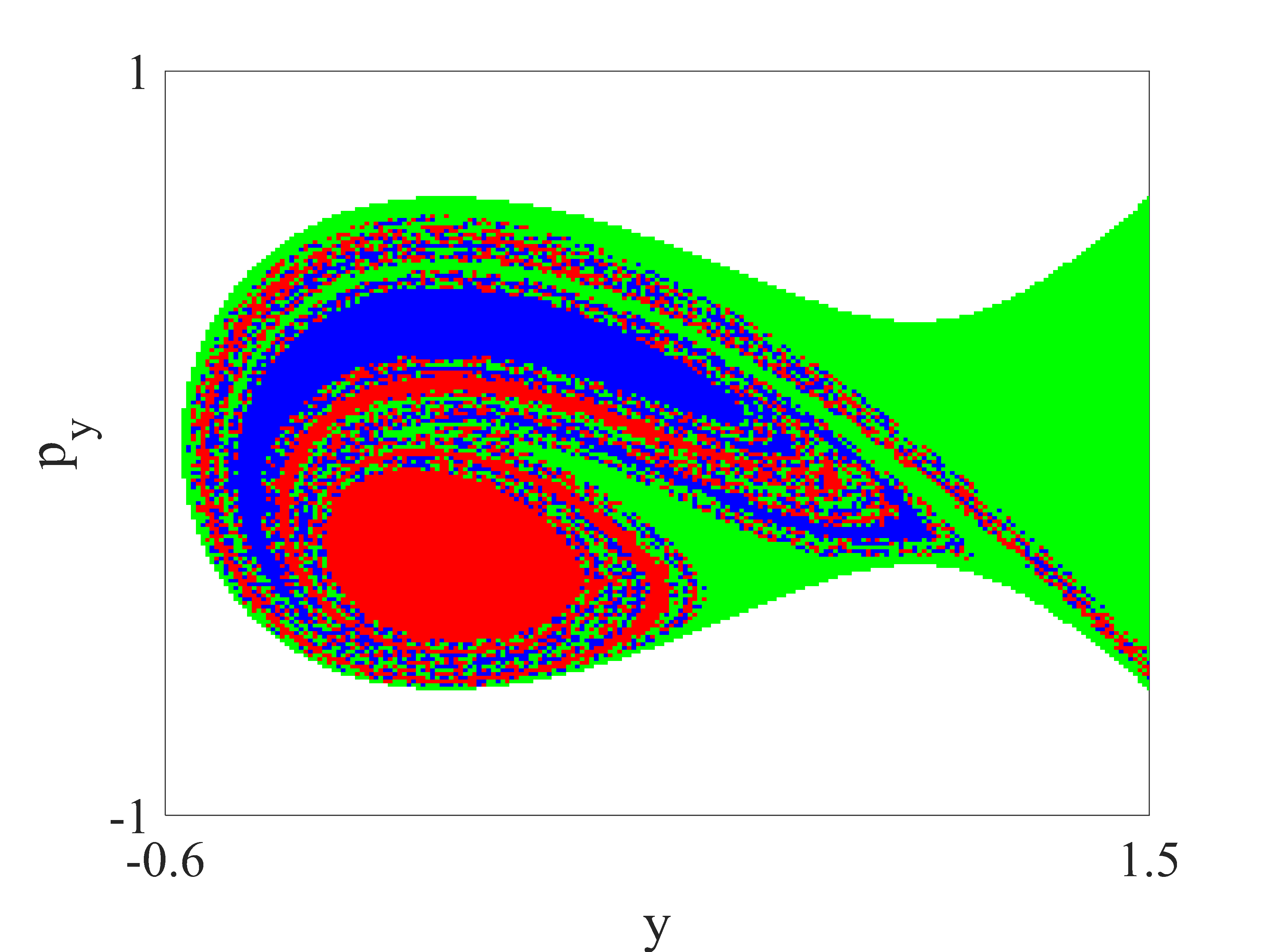}
					\caption{}\label{fig:henon-heiles-basins:j}
				\end{subfigure}\hfill
				\begin{subfigure}[t]{.31\linewidth}
					\centering
					\includegraphics[width=\linewidth]{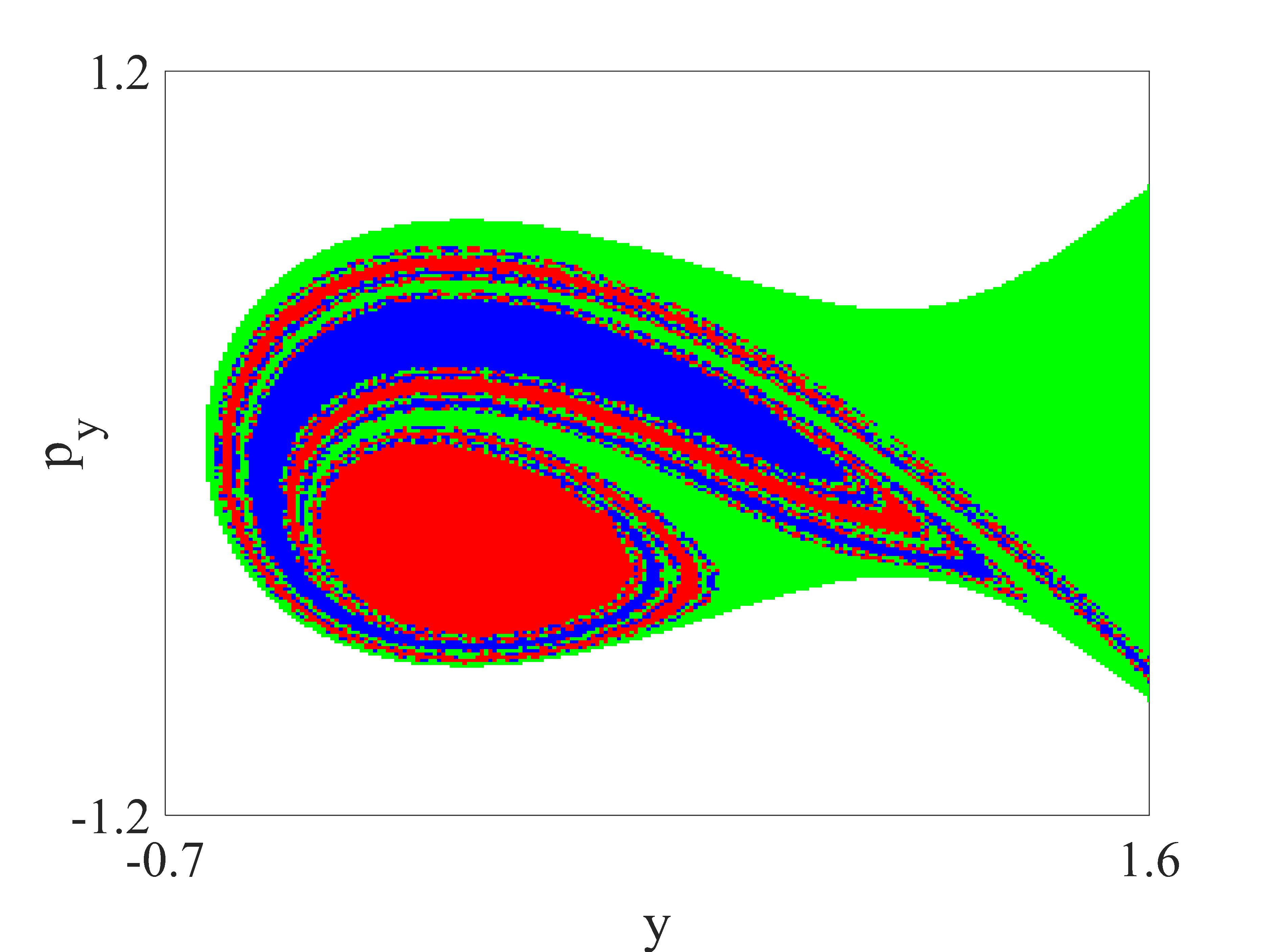}
					\caption{}\label{fig:henon-heiles-basins:k}
				\end{subfigure}\hfill
				\begin{subfigure}[t]{.31\linewidth}
					\centering
					\includegraphics[width=\linewidth]{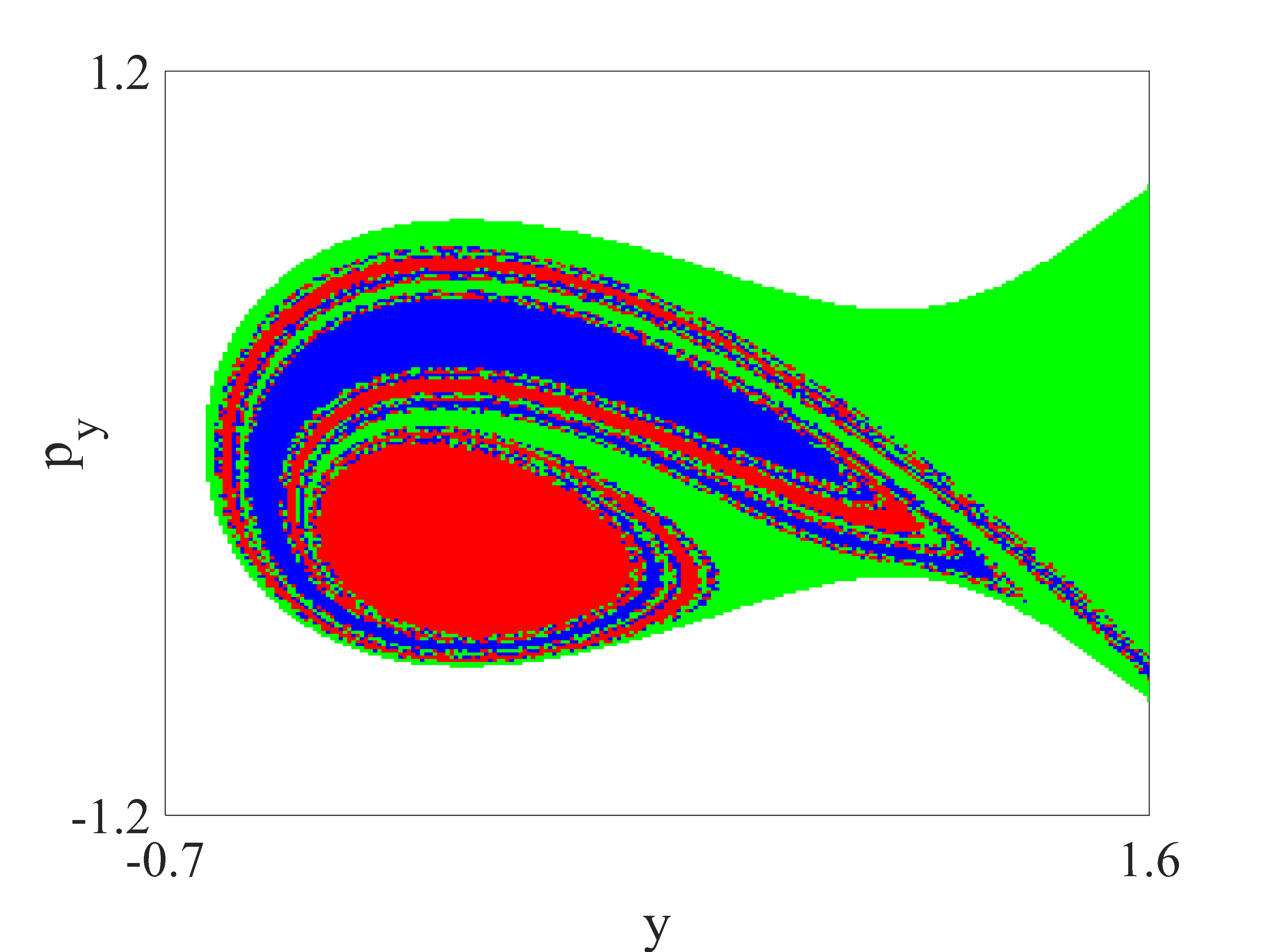}
					\caption{}\label{fig:henon-heiles-basins:l}
				\end{subfigure}\par\smallskip
				\begin{subfigure}[t]{.31\linewidth}
					\centering
					\includegraphics[width=\linewidth]{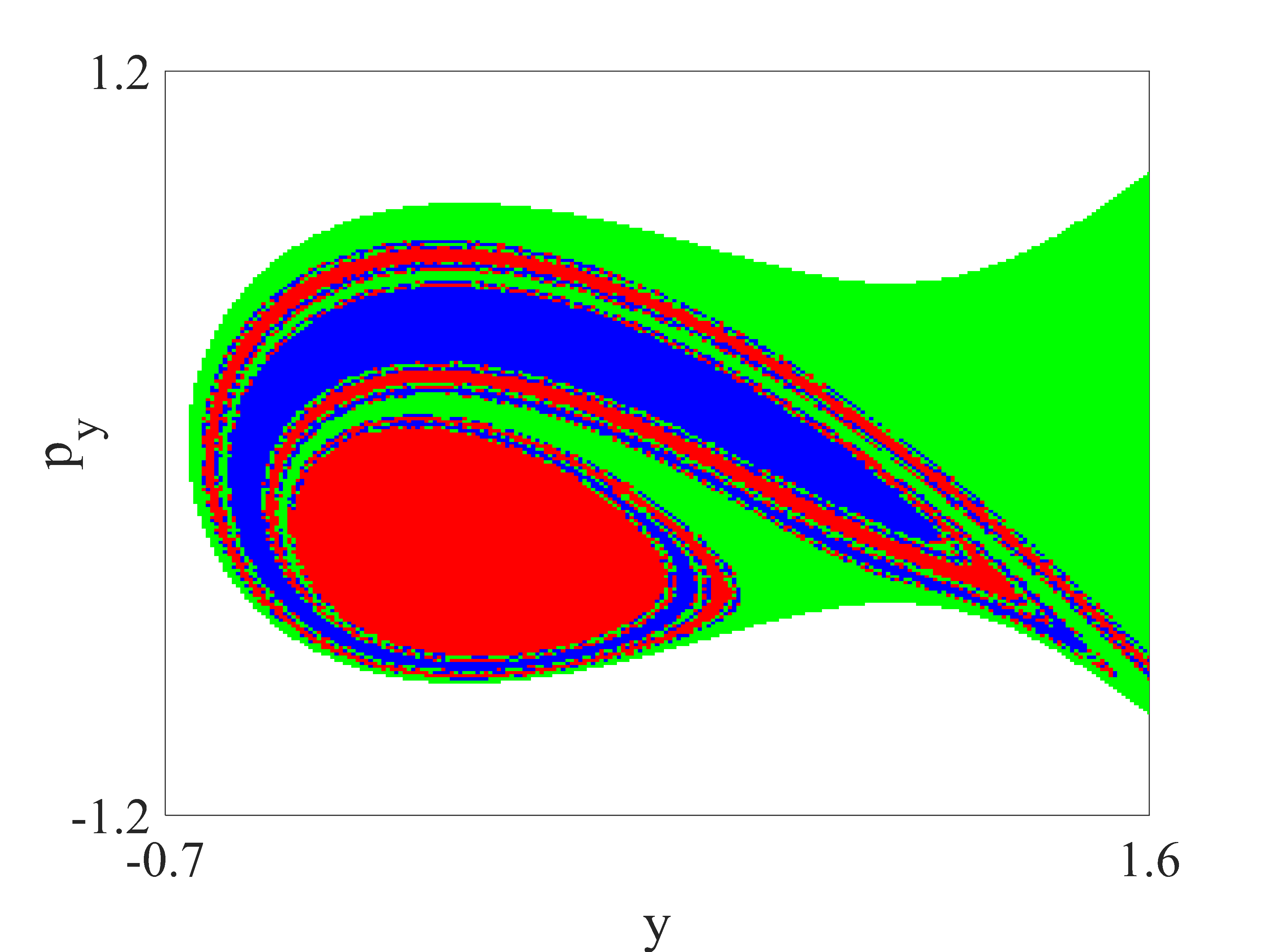}
					\caption{}\label{fig:henon-heiles-basins:m}
				\end{subfigure}\hfill
				\begin{subfigure}[t]{.31\linewidth}
					\centering
					\includegraphics[width=\linewidth]{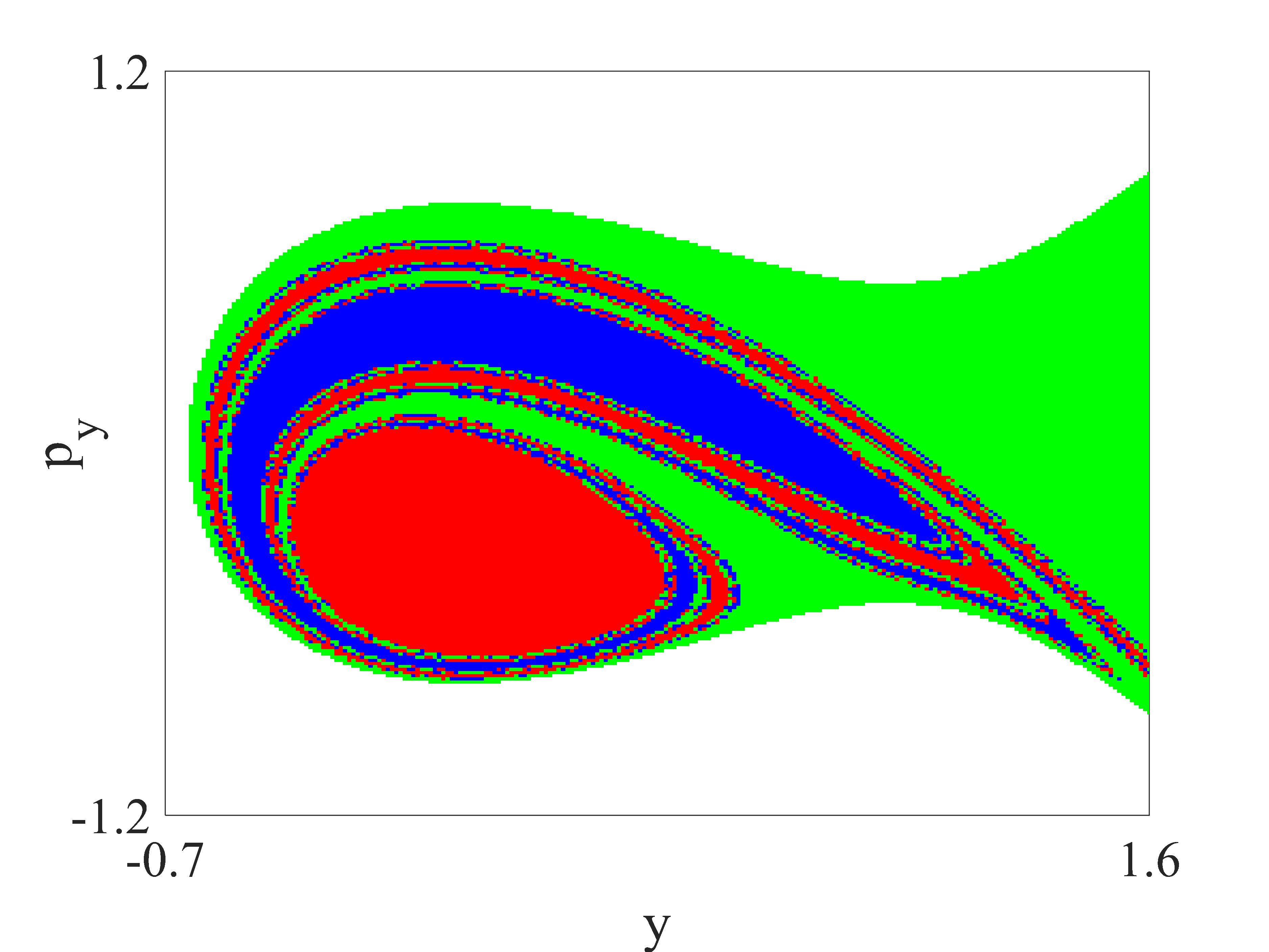}
					\caption{}\label{fig:henon-heiles-basins:n}
				\end{subfigure}\hfill
				\begin{subfigure}[t]{.31\linewidth}
					\centering
					\includegraphics[width=\linewidth]{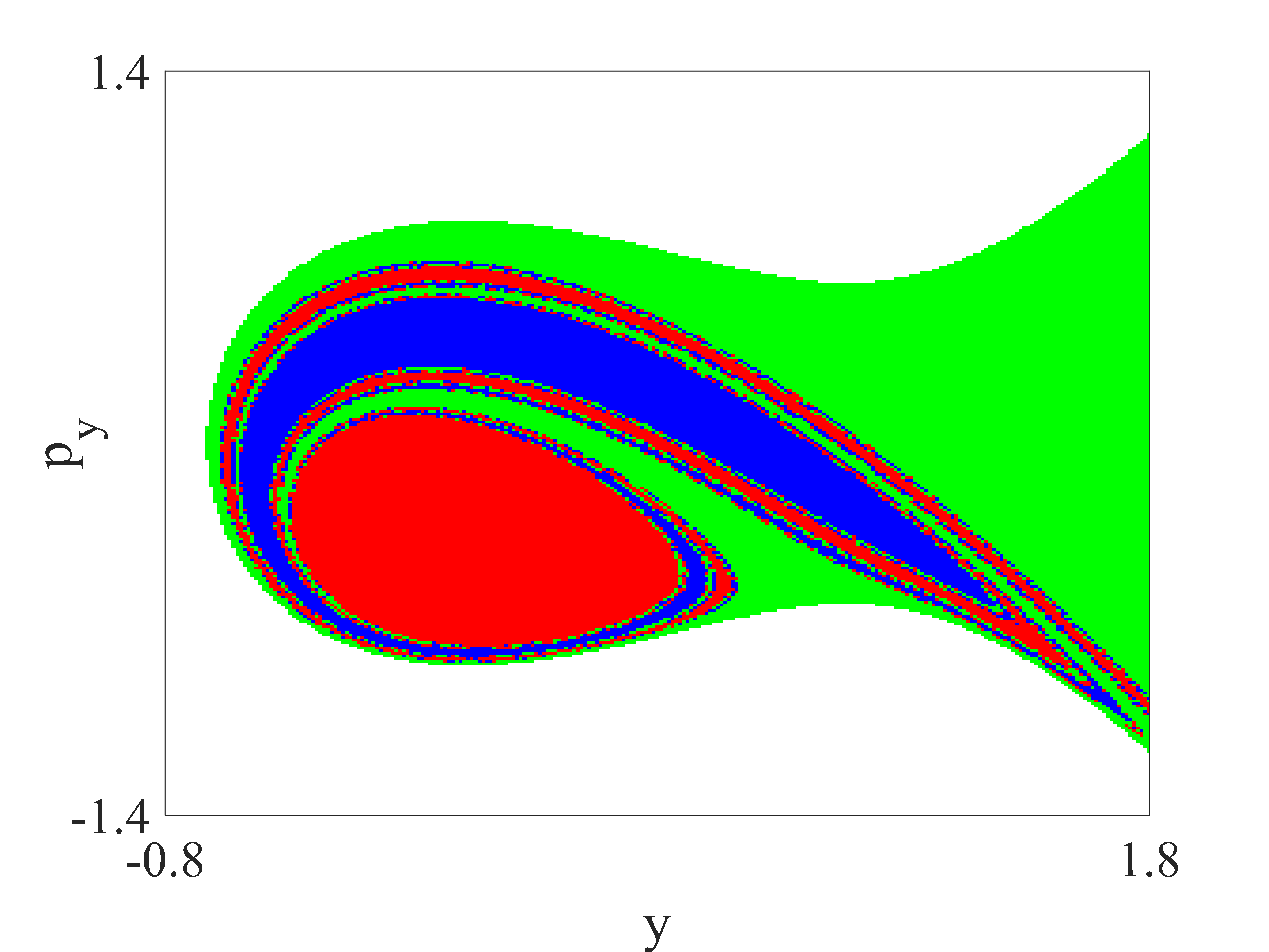}
					\caption{}\label{fig:henon-heiles-basins:o}
				\end{subfigure}\par\smallskip
				\begin{subfigure}[t]{.31\linewidth}
					\centering
					\includegraphics[width=\linewidth]{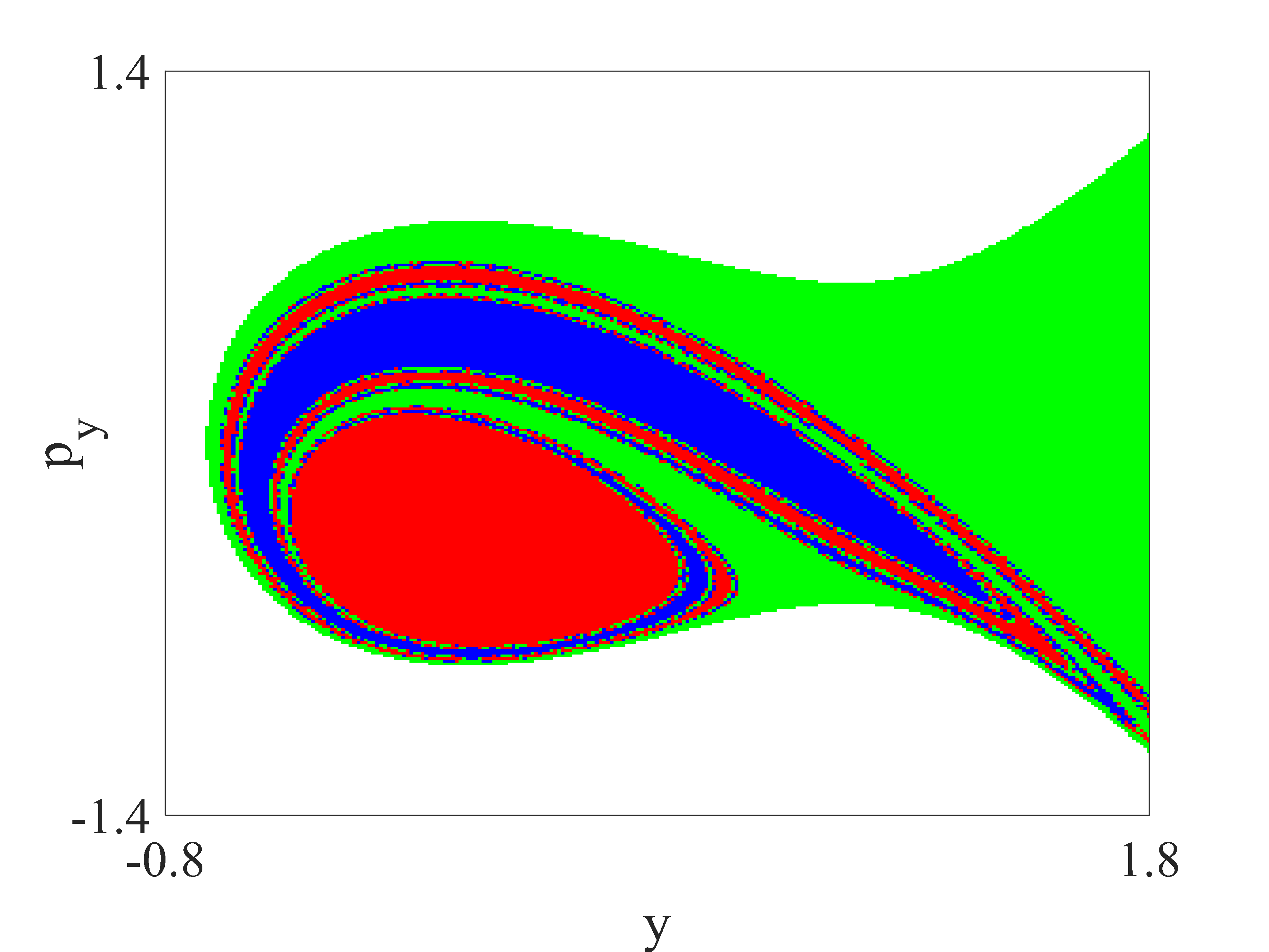}
					\caption{}\label{fig:henon-heiles-basins:p}
				\end{subfigure}\hfill
				\begin{subfigure}[t]{.31\linewidth}
					\centering
					\includegraphics[width=\linewidth]{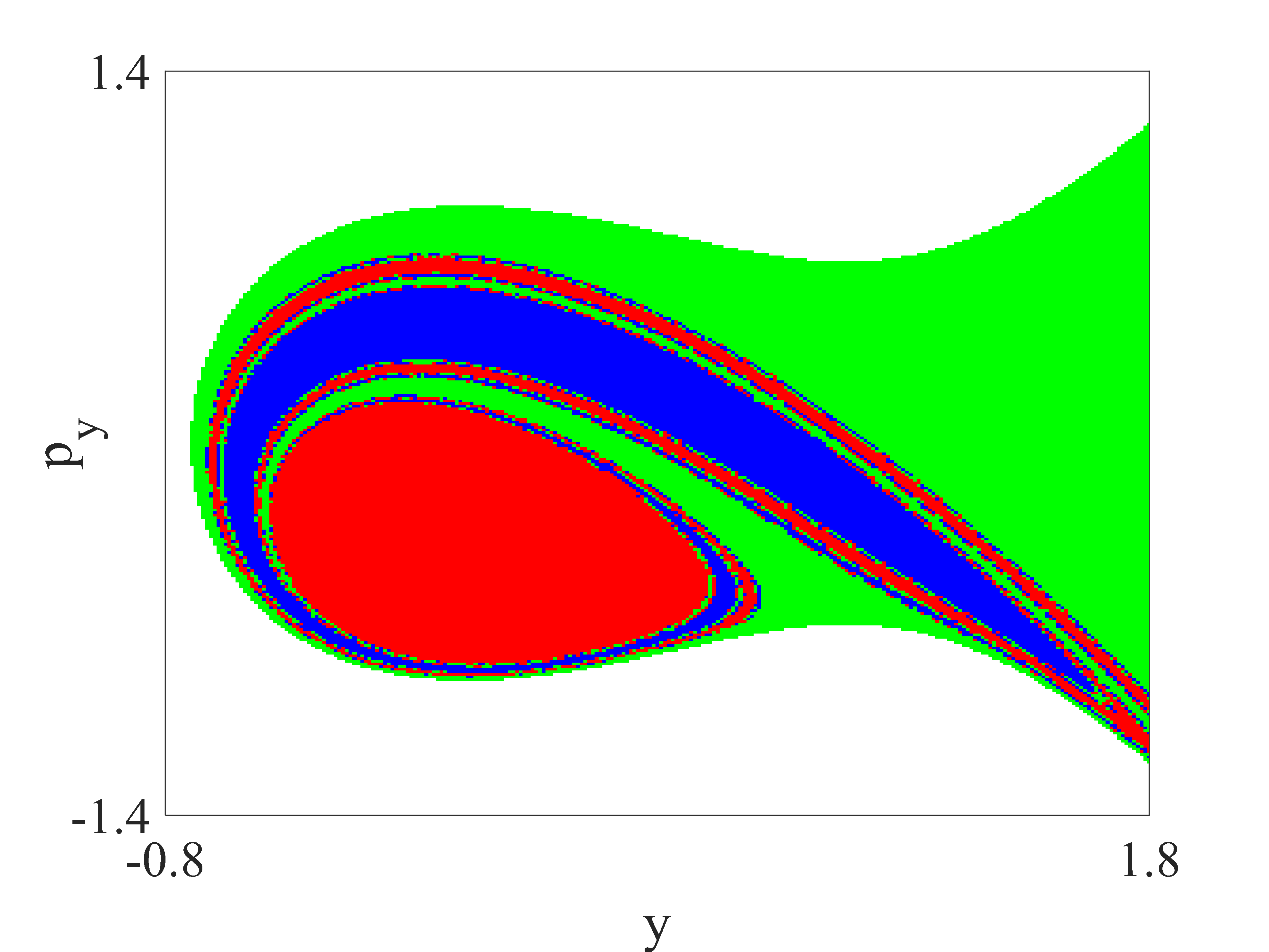}
					\caption{}\label{fig:henon-heiles-basins:q}
				\end{subfigure}\hfill
				\begin{subfigure}[t]{.31\linewidth}
					\centering
					\includegraphics[width=\linewidth]{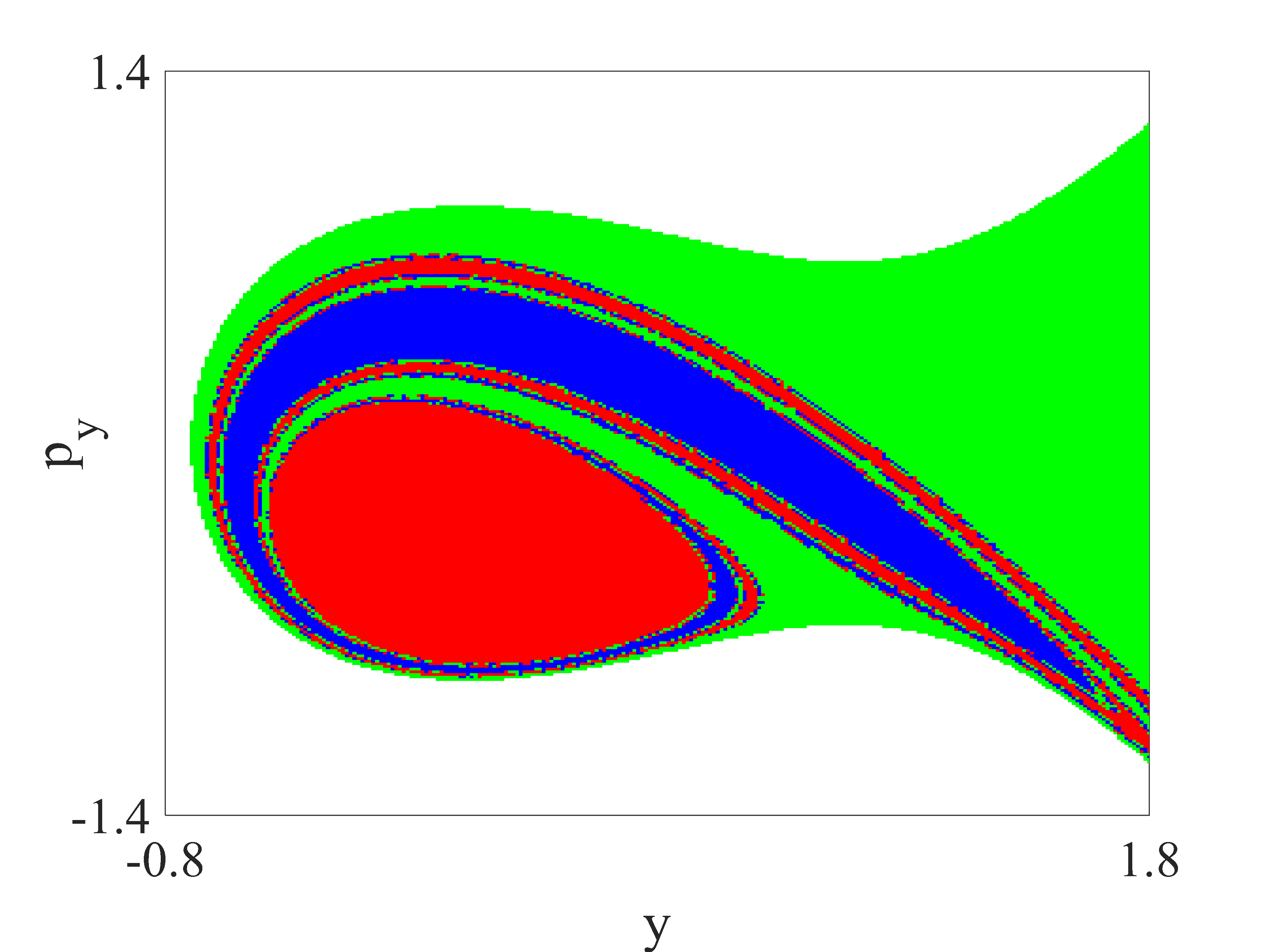}
					\caption{}\label{fig:henon-heiles-basins:r}
				\end{subfigure}
				\caption[]{\textit{Continued.}}
			\end{figure}
			
\begin{figure}[!htbp]
				\centering
				\begin{subfigure}[t]{.47\linewidth}
					\centering
					\includegraphics[width=\linewidth]{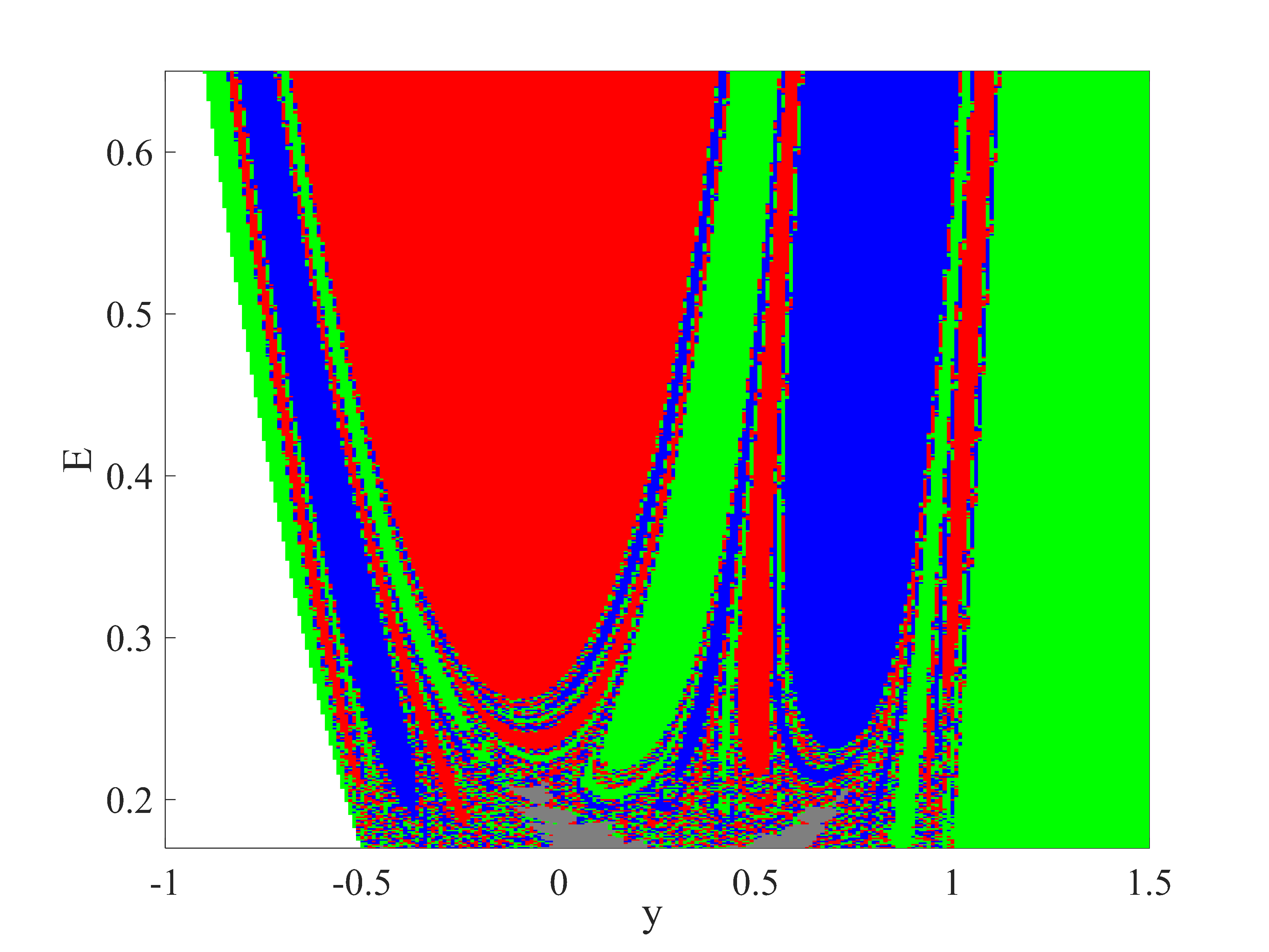}
					\caption{}\label{fig:henon-heiles-energy-basins:a}
				\end{subfigure}\hfill
				\begin{subfigure}[t]{.47\linewidth}
					\centering
					\includegraphics[width=\linewidth]{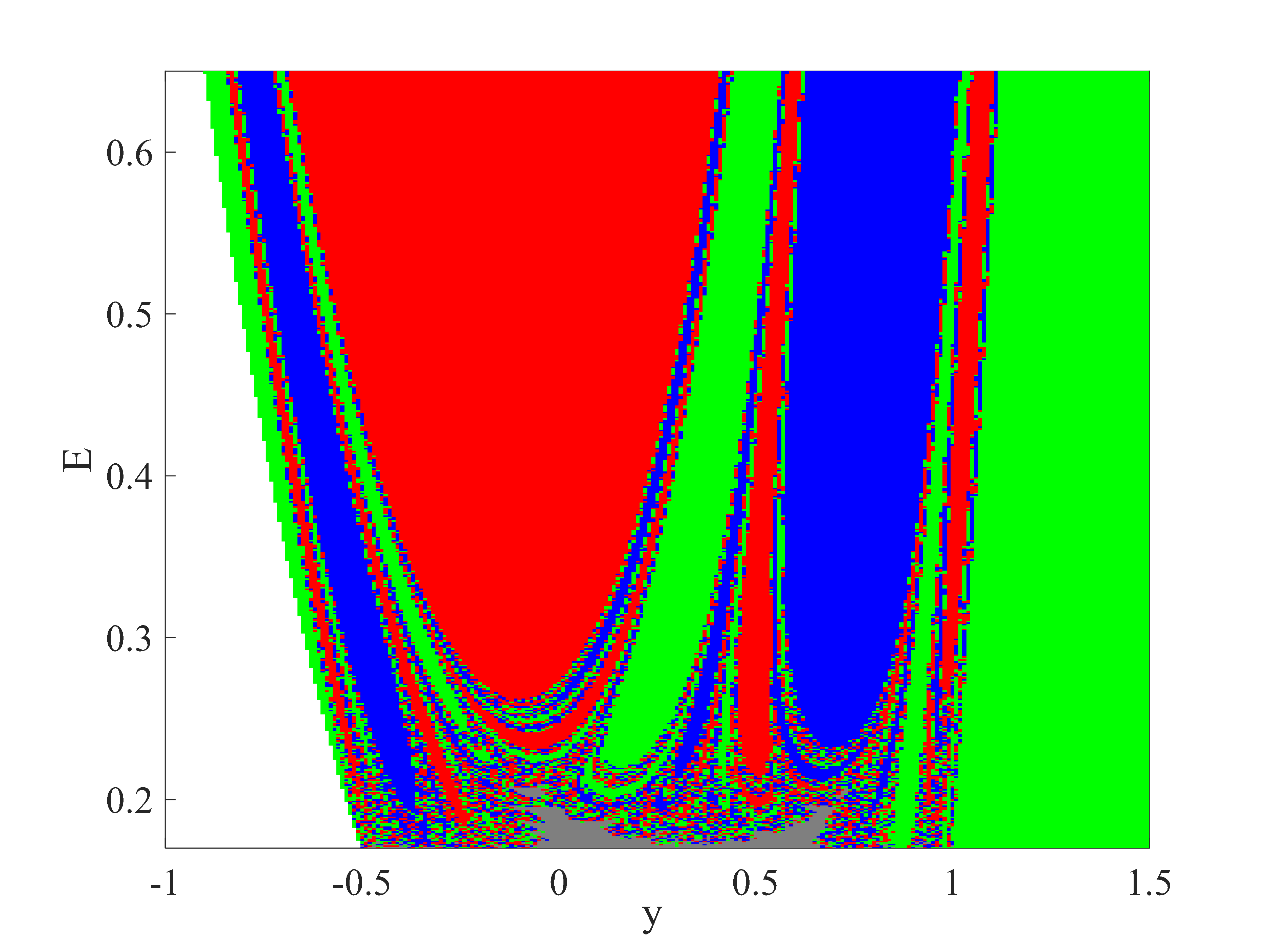}
					\caption{}\label{fig:henon-heiles-energy-basins:b}
				\end{subfigure}
				\caption{Prediction of the structure of the escape channels in the $(y,\allowbreak E)$ plane. Gray: trapped; green: escape through exit 1; blue: escape through exit 2; red: escape through exit 3. (\subref{fig:henon-heiles-energy-basins:a}) True structure. (\subref{fig:henon-heiles-energy-basins:b}) Predicted structure.}
				\label{fig:henon-heiles-energy-basins}
			\end{figure}

			Fix $x_{0}=0$ and $p_{y0}=0$. For a given $y_{0}$, we assign $p_{x0}=\sqrt{2E-y_{0}^{2}+\frac{2}{3}y_{0}^{3}}$, provided that $2E-y_{0}^{2}+\frac{2}{3}y_{0}^{3}\geqslant{}0$; then the total energy is $E$. Fig.~\ref{fig:henon-heiles-energy-basins}(\subref{fig:henon-heiles-energy-basins:a}) shows the true structure of escape channels in the $(y,\allowbreak E)$ plane, where $y_{0}$ takes values in $[-1,\allowbreak 1.5]$ and the energy $E$ takes values in $[0.17,\allowbreak 0.65]$. Fig.~\ref{fig:henon-heiles-energy-basins}(\subref{fig:henon-heiles-energy-basins:b}) presents the prediction result. The shapes and positions of the three exit channels are well reproduced, and deviations occur mainly in narrow bands along the boundaries separating different channels or separating escape from non-escape. These small discrepancies arise from the highly intricate, even fractal, nature of the true boundaries. Detailed calculations show that the prediction accuracy attains $90.4\%$. In fact, extending the energy range may further improve the accuracy.
			
			Compared with the two-exit system studied here, the H\'enon-Heiles system is more complex. Its three exit channels are arranged in a three-fold symmetric manner, and this arrangement gives rise to a larger number of boundaries, including those between adjacent exit channels and those between escape and non-escape regions, which can become highly intricate and, in some cases, even fractal. As a result, predicting the structure of the escape channels for the H\'enon-Heiles system is more challenging, which may contribute to the observed decrease in prediction accuracy relative to the two-exit case. Nevertheless, when we take into account both the prediction of the escape threshold and the reconstruction of the structure of escape channels across different energy values, the results remain satisfactory. This again demonstrates that escape dynamics can be effectively predicted by the NG-RC model.
			
			\Needspace{4\baselineskip}\section{Conclusion}\label{sec:conclusion}
			This study focuses on predicting the escape dynamics of Hamiltonian systems. The findings from numerical experiments on two distinct systems, one with two exit channels and the other with three, are summarized as follows. First, with only three time series as training data, the NG-RC model is able to predict both the escape threshold itself and the structure of the escape channels with high accuracy. The obtained results are based on the predicted trajectories, which, although not exact copies of the true orbits, capture the essential dynamical features of the system. Second, we observe that the prediction accuracy is energy dependent. Near the escape threshold, the system is extremely sensitive to initial conditions due to the presence of chaotic scattering and fractal basin boundaries. Consequently, predicting the correct escape channel becomes more challenging in this regime. As the energy increases, the prediction accuracy shows a gradual improvement overall. For the two-channel system, the accuracy exceeds 98\% at higher energy, and for the three-channel H\'enon-Heiles system, it exceeds 95\%.
			
			Notably, the NG-RC model is trained \textit{solely on non-escaping data, yet it successfully predicts supercritical escape dynamics.} This outcome is somewhat unexpected, as the training set does not explicitly contain information about how trajectories leave the potential well. The fact that the model can infer the global topology of the potential landscape from low-energy data alone highlights its powerful extrapolation capability. Exploring the broader application of machine learning algorithms in data-driven modeling of Hamiltonian systems remains an important direction for future research.
			
			\subsection*{Statements and Declarations}
			The authors declare that they have no known competing financial interests or personal relationships that could have appeared to influence the work reported in this paper.
			
			\subsection*{Acknowledgments}
			Jianming Liu was supported by the National Natural Science Foundation of China (Grant No. 12602043). Zhicheng Tong was supported by the National Natural Science Foundation of China (Grant No. 12601347) and the China Postdoctoral Science Foundation (Grant No. 2025M783102).
			
			\subsection*{CRediT author statement}
			Jianming Liu: Conceptualization, Methodology, Writing - original draft, Software, Writing - review \& editing.\vspace{2mm}
			
			\noindent Zhicheng Tong: Conceptualization, Methodology, Writing - review \& editing, Validation.
			
			\subsection*{Data availability}
			The code for the numerical simulations and the data generated in this study are publicly available on GitHub at \url{https://github.com/LiuJM1026/escape_dynamics.git}
			
			\FloatBarrier
			
		\end{document}